\documentclass[10pt,
twoside]{amsart}
\usepackage{lmodern} 
\usepackage{soul, wrapfig}
\usepackage{graphicx}
\usepackage{color}
\usepackage{amsmath}
\usepackage{amssymb}
\usepackage{amsfonts}
\usepackage{latexsym, amsthm, mathrsfs}
\usepackage{hyperref, breakurl}
\usepackage[OT1]{fontenc}
\usepackage{fancybox}
\usepackage{amscd}
\input{xy}
\usepackage[all]{xy}

\numberwithin{equation}{section}

\newtheorem{theorem}{Theorem}[section]
\newtheorem{lemma}[theorem]{Lemma}

\newtheorem*{theorem1}{Theorem \ref{1}}
\newtheorem*{theorem2}{Theorem \ref{2}}
\newtheorem*{theoremexistence}{Theorem \ref{gmtu_existence}}
\newtheorem*{theoremgowers}{Theorem \ref{gowersown1}}
\newtheorem*{theoremstructure}{Theorem \ref{structure}}
\newtheorem*{theorem4}{Theorem \ref{4}}
\newtheorem*{corollary4}{Corollary \ref{cor4}}

\theoremstyle{definition} 
\newtheorem{definition}[theorem]{Definition}

\theoremstyle{plain} 
\newtheorem{question}[theorem]{Question}

\newtheorem{proposition}[theorem]{Proposition}
\newtheorem{corollary}[theorem]{Corollary}
\newtheorem{fact}[theorem]{Fact}

\newtheorem*{maintheorem*}{Main Theorem}
\newtheorem*{conjecture*}{Conjecture}
\newtheorem*{theorem*}{Theorem}
\newtheorem*{proposition*}{Proposition}
\newtheorem*{corollary*}{Corollary}

\theoremstyle{remark}  
\newtheorem{remark}[theorem]{Remark}
\newtheorem*{remarks*}{Remarks}
\newtheorem*{remark*}{Remark}
\newtheorem*{claim*}{Claim}

\newtheoremstyle{mystyle}
  {3pt}{3pt}{\itshape}{}{\bfseries}{}{.5em}
  {\thmname{#1}\thmnumber{ #2}\thmnote{. #3}}
\theoremstyle{mystyle}
\newtheorem*{maintheoremnew*}{Main Theorem}
\newtheorem*{conjecturenew*}{Conjecture}
\newtheorem*{theoremnew*}{Theorem}
\newtheorem*{propositionnew*}{Proposition}

\newcommand{\nc}{\newcommand}
\nc{\comment}[1]{#1}
\nc{\nothing}[1]{}
\nothing{   
\nc{\dom}{{\rm dom}}
\nc{\card}{{\rm card}}
\nc{\lh}{{\rm lh}}
\nc{\lgg}{{\rm lg}}
\nc{\rge}{\mbox{\rm range}}
\nc{\cf}{{\rm cf}}
nc{\nex}{\mbox{\rm next}}
\nc{\uhr}{\restriction}
\nc{\supt}{{\rm supt}}
\nc{\supp}{{\rm supp}}
\nc{\Lim}{{\rm Lim}}
\nc{\Leb}{{\rm Leb}}
\nc{\modd}{{\rm mod}}
\nc{\RO}{{\rm RO}}
\nc{\prob}{{\rm Prob}}
}
\nc{\minl}{{\min}} 
\nc{\maxl}{{\max}} 
\nc{\dotpl}{{+}}
\nc{\dotunion}{{\cup}} 
\nc{\On}{{\rm On}}
\nc{\Ord}{{\rm On}}
\nc{\nco}{\DeclareMathOperator}
\nco{\Tl}{{T}}
\nco{\lift}{lift}
\nco{\rk}{rk}
\nco{\order}{o}
\nco{\ppower}{pp}
\nco{\pcf}{pcf} 
\nco{\tcf}{tcf} 
\nco{\tlim}{tlim} 
\nco{\limtext}{lim} 
\nco{\prodt}{{\textstyle \prod}}
\nco{\symdiff}{\triangle}
\nco{\dom}{dom}
\nco{\card}{card}
\nco{\lh}{lh}
\nco{\lt}{lt}
\nco{\lgg}{lg}
\nco{\hgt}{ht}
\nco{\rge}{range}
\nco{\otp}{otp}
\nco{\trunk}{tr}
\nco{\cf}{cf}
\nco{\nex}{next}
\nc{\uhr}{\restriction}
\nco{\reduction}{red}
\nco{\supt}{supt}
\nco{\supp}{supp}
\nco{\Lim}{Lim}
\nco{\Leb}{Leb}
\nco{\modd}{mod}
\nco{\invariant}{inv}
\nco{\invalues}{invalues}
\nco{\id}{id}
\nco{\RO}{RO}
\nco{\poss}{pos}
\nco{\Inc}{Inc} 
\nco{\Ge}{Ge}
\nco{\hdrop}{\hat{drop}}
\nco{\Gen}{Gen}
\nco{\Property}{Pr}
\nco{\GT}{GT}
\nco{\refl}{refl}
\nco{\nmg}{nmg}
\nco{\en}{en}
\nco{\truth}{truth} 
\nco{\sw}{sw} 

\nc{\potom}{\ensuremath{{\cal P}(\omega)}}
\nc{\potinf}{\ensuremath{[\omega]^\omega}}
\nc{\pfin}{\ensuremath{{\cal P}(\omega)/{\rm fin}}}

\nc{\potfin}{\ensuremath{[\omega]^{<\omega}}}
\nc{\inn}{\ensuremath{{\omega^{\uparrow \omega}}}}
\nc{\baire}{{}^\omega \omega}
\nc{\bair}{\omega^\omega}
\nc{\hoch}{^{<\omega}}
\nc{\hocho}{^{\omega}}
\nc{\tree}[1]{{[} #1 {]}_0}
\nc{\tre}[2]{ {#1}_{#2}}

\nc{\prooff}[1]{{\bf Proof} of #1:}
\nc{\proofend}{\makebox{} \hfill ${\bf \square}$ \\}
\nc{\proofendof}[1]{\makebox{} \hfill ${\boldmath{\square}}_{\rm #1}$ \\}
\nc{\beq}{\begin{eqnarray*}}
\nc{\eeq}{\end{eqnarray*}}
\nc{\bde}{\begin{list}}
\nc{\ede}{\end{list}}

\newenvironment{myrules}
{\begin{list}{}
{
 \setlength{\leftmargin}{0.8cm}
 \setlength{\labelwidth}{0.8cm}
 \setlength{\labelsep}{0.2cm}
 \setlength{\parsep}{0.5ex plus 0.2ex minus 0.1 ex}
 \setlength{\itemsep}{0.3ex plus 0.2 ex minus 0ex}
}}{\end{list}}

\newenvironment{myrules1}
{\begin{list}{}
{
 \setlength{\leftmargin}{1.0cm}
 \setlength{\labelwidth}{0.8cm}
 \setlength{\labelsep}{0.3cm}
 \setlength{\parsep}{0.5ex plus 0.2ex minus 0.1 ex}
 \setlength{\itemsep}{0.5ex plus 0.2 ex minus 0ex}
}}{\end{list}}

{\end{list}}

\newcounter{subalph}
{\end{list}}

\newcommand{\greek}[1]{\ifthenelse{\value{#1}=1}{\mbox{$\alpha$}}%
  {\ifthenelse{\value{#1}=2}{\mbox{$\beta$}}{%
   \ifthenelse{\value{#1}=3}{\mbox{$\gamma$}}{%
   \ifthenelse{\value{#1}=4}{\mbox{$\delta$}}{%
   \ifthenelse{\value{#1}=5}{\mbox{$\varepsilon$}}{%
   \ifthenelse{\value{#1}=6}{\mbox{$\zeta$}}{%
   \ifthenelse{\value{#1}=7}{\mbox{$\eta$}}{%
   \ifthenelse{\value{#1}=8}{\mbox{$\theta$}}{%
   \ifthenelse{\value{#1}=9}{\mbox{$\iota$}}{%
   \ifthenelse{\value{#1}=10}{\mbox{$\kappa$}}{%
   \ifthenelse{\value{#1}=11}{\mbox{$\lambda$}}{%
   \ifthenelse{\value{#1}=12}{\mbox{$\mu$}}{%
   \ifthenelse{\value{#1}=13}{\mbox{$\nu$}}{%
   \ifthenelse{\value{#1}=14}{\mbox{$\xi$}}{%
   \ifthenelse{\value{#1}=15}{\mbox{$\rm o$}}{%
   \ifthenelse{\value{#1}=16}{\mbox{$\pi$}}{%
   \ifthenelse{\value{#1}=17}{\mbox{$\varrho$}}{%
   \ifthenelse{\value{#1}=18}{\mbox{$\sigma$}}{%
   \ifthenelse{\value{#1}=19}{\mbox{$\tau$}}{%
   \ifthenelse{\value{#1}=20}{\mbox{$\upsilon$}}{%
   \ifthenelse{\value{#1}=21}{\mbox{$\varphi$}}{%
   \ifthenelse{\value{#1}=22}{\mbox{$\chi$}}{%
   \ifthenelse{\value{#1}=23}{\mbox{$\psi$}}{\mbox{$\omega$}%
  }}}}}}}}}}}}}}}}}}}}}}}}

\newcounter{subgreek}
{\end{list}}

\newcounter{subarabic}
{\end{list}}

\newcounter{subroman}
{\end{list}}

\newcount\skewfactor

\def\mathunderaccent#1#2 {\let\theaccent#1\skewfactor#2
\mathpalette\putaccentunder}
\def\putaccentunder#1#2{\oalign{$#1#2$\crcr\hidewidth
\vbox to.2ex{\hbox{$#1\skew\skewfactor\theaccent{}$}\vss}\hidewidth}}
\def\name{\mathunderaccent\tilde-3 }

\nc{\nname}{\name}

\nc{\even}{\ensuremath{\rm Even}}
\nc{\odd}{\ensuremath{\rm Odd}}

\nc{\al}{$\alpha$\  }
\nc{\om}{\omega}
\nc{\omm}{\ensuremath{\omega_1}}
\nc{\ep}{\varepsilon}
\nc{\tk}{\tilde{K}}
\nc{\concat}{\,\hat{} \,}   
\nc{\force}{\Vdash}
\nc{\fb}{f_{\overline{M}}}
\nc{\such}{\, : \,}
\newcommand{\la}{\langle}
\newcommand{\ra}{\rangle}

\nc{\meager}{\ensuremath{{\cal M}}}
\nc{\lebesgue}{\ensuremath{{\cal N}}}
\nc{\nulll}{\ensuremath{{\cal N}}}
\nc{\ksigma}{\ensuremath{{\bf K}_\sigma}}
\nc{\ideal}{\ensuremath{{\cal I}}}
\nc{\ga}{\ensuremath{\frak a}}
\nc{\AAA}{{\cal A}}   
\nc{\gc}{\ensuremath{\frak c}}
\nc{\gs}{\ensuremath{\frak s}}
\nc{\gh}{\ensuremath{\frak h}}
\nc{\gd}{\ensuremath{\frak d}}
\nc{\gb}{\ensuremath{\frak b}}
\nc{\gro}{\ensuremath{\frak g}}
\nc{\gu}{\ensuremath{\frak u}}
\nc{\gr}{\ensuremath{\frak r}}
\nc{\gt}{\ensuremath{\frak t}}
\nc{\fff}{\ensuremath{\frak f}}
\nc{\gm}{\ensuremath{\mathfrak{mcf}}}
\nc{\gge}{\ensuremath{\mathfrak e}}
\nc{\cfupro}{\ensuremath{\cf(\upro)}}
\nc{\cfvpro}{\ensuremath{\cf(\vpro)}}
\nc{\gp}{\ensuremath{\frak p}}
\nc{\gk}{\ensuremath{\frak k}}

\nc{\add}{\mbox{\ensuremath{{\rm add}}}}
\nc{\cov}[1]{\mbox{\ensuremath{{\rm cov}(#1)}}}
\nc{\unif}[1]{\mbox{\ensuremath{{\rm unif}(#1)}}}
\nc{\cof}[1]{{\mbox{\ensuremath{\rm cof}(#1)}}}

\nc{\addd}[2]{\mbox{\ensuremath{{\rm add}^{#1}(#2)}}}   
\nc{\covv}[2]{\mbox{\ensuremath{{\rm cov}^{#1}(#2)}}}   
\nc{\uniff}[2]{\mbox{\ensuremath{{\rm unif}^{#1}(#2)}}} 
\nc{\coff}[2]{{\mbox{\ensuremath{\rm cof}^{#1}(#2)}}}

\nc{\cd}{Cicho\'n's Diagram}

\nc{\MA}{\mbox{\sf MA}}
\nc{\PFA}{\mbox{\sf PFA}}
\nc{\OCA}{\mbox{\sf OCA}}
\nc{\GCH}{\mbox{\sf GCH}}
\nc{\CH}{\mbox{\sf CH}}
\nc{\zfc}{\mbox{\sf ZFC}}
\nc{\sch}{\mbox{\sf SCH}}
\nc{\ZF}{\mbox{\sf ZF}}
\nc{\NCF}{\mbox{\sf NCF}} 
\nc{\FD}{\mbox{\sf FD}}   
\nc{\SFT}{\mbox{\sf SFT}}
\nc{\fourG}{\mbox{\rm 4G}}
\nc{\fourI}{\mbox{\rm 4I}}
\nc{\past}{\ ;{\rm past}\;}
\nc{\Borelhood}{Borel measurability} 
\nc{\Pieinseins}{\mbox{${\bf \Pi}^1_1$}}
\nc{\seinseins}{\mbox{${\bf\Sigma}^1_1$}}
\nc{\seinszwei}{\mbox{${\bf\Sigma}^1_2$}}
\nc{\seinsdrei}{\mbox{${\bf\Sigma}^1_3$}}
\nc{\Deleinszwei}{\mbox{${\bf\Delta}^1_2$}}

\nc{\up}{\ensuremath{{\cal U}\mbox{\ensuremath{\rm -prod}}\,\omega}}
\nc{\upp}{\ensuremath{{\cal U}'\mbox{\ensuremath{\rm -prod}}\,\omega}}
\nc{\upro}{\ensuremath{{\cal U}\mbox{\ensuremath{\rm -prod}}\,\om}}
\nc{\fupro}{\ensuremath{f({\cal U})\mbox{\ensuremath{\rm -prod}}\,\om}}
\nc{\vpro}{\ensuremath{{\cal V}\mbox{\ensuremath{\rm -prod}}\,\om}}
\nc{\fpro}{\ensuremath{{\cal F}\mbox{\ensuremath{\rm -prod}}\,\om}}

\nc{\cff}[1]{{\text{cf}\,(#1)}}           
\nc{\cu}{\ensuremath{\cal U}}             
\nc{\ai}{\ensuremath{\forall^\infty}}     
\nc{\ei}{\ensuremath{\exists^\infty}}     
\nc{\ww}{\ensuremath{\omega^\omega}}      

\nc{\gw}{groupwise dense}

\nc{\kk}{car\-dinal cha\-rac\-teris\-tic}
\nc{\joker}{\ast}
\nc{\gtc}{Galois-Tukey connection} 
\nc{\ntc}{not nearly coherent} 
\nc{\nnc}{non-nearly coherent} 

\nc{\av}[1]{{\rm Av}_{#1}}
\nc{\eps}{\varepsilon}
\renewcommand{\epsilon}{\varepsilon}

\nc{\n}{{\bf n}}                 
\nc{\m}{{\bf m}}

\nc{\marginparr}[1]{}
\nc{\footnoteee}{} 
\nc{\footnotee}{}  

\newcommand{\cal}{\mathcal}

\nc{\divs}{{c_0 \setminus \ell^1}}
\nc{\divser}{(\divs, \leq^*)/\thickapproy}
\nc{\bfin}{\RO(\pfin \setminus\{0\},\subseteq^*)}
\nc{\bdivser}{\RO(\divser)}
\nc{\inc}{{\rm INC}}
\nc{\com}{{\rm COM}}
\nc{\thickapproy}{\makebox{}\!\!\thickapprox}
\nc{\approy}{\makebox{}\!\!\approx}
\nc{\lessi}{\leqslant}
\nc{\gessi}{\geqslant}
\nc{\interior}[1]{{\rm int}(#1)}
\nc{\closure}[1]{{\rm cl}(#1)}
\nc{\Vo}{Vojt\'a\v{s}}

\nc{\precedeseq}{\leq^*} 
\nc{\precedes}{\prec}
\nc{\stronger}{\leqslant_{\bf P}}
\nc{\underlline}[1]{\hat{#1}}
\nc{\PO}{{\bf P}}
\nc{\charak}{\text{ch}}
\nc{\symom}{{\rm{Sym}(\omega)}}

\nc{\needed}{needed\ }
\nc{\neededc}{needed}
\nc{\Needed}{Needed\ }
\nc{\wneeded}{weakly needed\ }
\nc{\Wneeded}{Weakly needed\ }
\nc{\wneededc}{weakly needed}

\nc{\mup}{m_{\rm up}}
\nc{\mdn}{m_{\rm dn}}
\nco{\may}{may}
\nco{\aver}{av} 
\nco{\norm}{nor} 
\nco{\val}{val} 
\nco{\dis}{dis} 
\nco{\basis}{basis}
\nco{\pos}{pos}
\nco{\spec}{spec}
\nc{\err}{\mbox{err}}
\nc{\eee}{\mbox{e}}
\nco{\Expect}{Exp}
\nco{\rt}{rt}
\nco{\pr}{pr}
\nco{\suc}{suc}
\nco{\splitt}{sp}
\nco{\halv}{h}
\nco{\Add}{Add}
\nco{\Cov}{Cov}
\nco{\Unif}{Unif}
\nco{\Cof}{Cof}
\nco{\htt}{ht}
\nco{\cl}{cl}
\nc{\bbforcing}{\mathbb A}
\nc{\itername}{\mathfrak q}
\nc{\iterp}{\mathfrak p}
\nc{\iterq}{\mathfrak q}
\nc{\invcm}{\rm inv_{cm}}
\nc{\invcf}{\rm inv_{cf}}
\nc{\invgm}{\rm inv_{gm}}
\renewcommand{\P}{\mathbb P}
\nc{\Q}{\mathbb Q}
\nc{\subsetsim}{\underset{\raise0.6em\hbox{$\sim$}}{\subset}}
\newcommand{\subsim}{\underset{\raise20pt\hbox{$\rightarrow$}}{\rightarrow}}
\newcommand{\ssim}{\overset{\raise-40pt\hbox{$\leftarrow$}}{\subsim}}
\nc{\rest}{\restriction}
\nc{\bF}{\mathbb F}
\nc{\F}{{\rm{Fin}}} 
\nc{\Fk}{{\rm{Fin}}_k}
\nc{\Ftwo}[1]{{\rm{Fin}}_{#1}}
\nc{\R}{\mathbb R}
\nc{\M}{\mathbb M}
\nc{\GM}{\mathbb{GM}} 
\nc{\strk}{\sqsubseteq_k^{\rm Tetr}} 
\nc{\strkast}{\sqsubseteq_k^{\rm Tetr,\ast}}
\nc{\strtwo}[1]{\sqsubseteq_{#1}^{\rm Tetr}}
\nc{\bQ}{\mathbb Q}
\nc{\bP}{\mathbb P}
\nc{\bG}{\mathbf G}

\nc{\cJ}{\mathcal J}
\nc{\cE}{\mathcal E}
\nc{\cP}{\mathcal P}
\nc{\cU}{\mathcal U}
\nc{\cV}{\mathcal V}
\nc{\cW}{\mathcal W}
\nc{\cC}{\mathcal C}
\nc{\cT}{\mathcal T}
\nc{\cX}{\mathcal X}
\nc{\cD}{\mathcal D}
\nc{\cS}{\mathcal S}

\nc{\cG}{\mathcal G}
\nc{\cI}{\mathcal I}
\nc{\cF}{\mathcal F}
\nc{\cH}{\mathcal H}
\nc{\cM}{\mathcal M}
\nc{\cA}{\mathcal A}
\nc{\cB}{\mathcal B}
\nc{\cR}{\mathcal R}

\nc{\god}{\mathfrak{od}} 
\nc{\mcf}{\mathfrak{mcf}}
\nc{\roth}{[\omega]^{\omega}}
\nco{\last}{last} 
\nco{\filter}{fil}
\nco{\semifilter}{semi}
\nco{\sfil}{sfil}
\nco{\ssfil}{ssfil}
\nco{\supersets}{cl} 

\nco{\upwards}{up}
\nco{\inter}{inter}
\nco{\FU}{FU} 
\nco{\FUU}{FU} 
\nco{\TFU}{TFU} 
\nco{\set}{set}
\nco{\seq}{seq}
\nc{\forks}{\underset{\raise0.6em\hbox{$\smile$}}{\mid}}
\renewcommand{\setminus}{\smallsetminus}
\nc{\bV}{{\bf V}}\nc{\bW}{{\bf W}}
\nco{\fil}{fil}
\nco{\CFF}{CFF}
\nco{\rfl}{rfl}

\nc{\fin}{\emptyset}
\nc{\ba}{\bar{a}}\nc{\bs}{\bar{s}}
\nc{\bb}{\bar{b}}
\nc{\bc}{\bar{c}}
\nc{\bd}{\bar{d}}
\nc{\be}{\bar{e}}
\nc{\kor}{\mbox{kor}}
\nc{\bx}{\bar{x}}
\nc{\bg}{\bar{g }}
\nc{\mtu}{Milliken--Taylor ultrafilter} 
\nc{\kmtu}{Milliken--Taylor\- ultrafilter}
\nc{\gmtu}{Gowers--Milliken--Taylor\- ultrafilter}
\renewcommand{\rho}{\varrho}
\begin{document}


\title[Finitely Many Near-Coherence Classes \today]{Local Ramsey Spaces in Matet Forcing Extensions and Finitely Many Near-Coherence Classes}

\author{Heike Mildenberger}

\address{Heike Mildenberger, Albert-Ludwigs-Universit\"at Freiburg,
Mathematisches Institut, Abteilung f\"ur math. Logik, Ernst--Zermelo--Stra\ss e~1, 
79104 Freiburg im Breisgau, Germany}


\begin{abstract}
 We introduce Gowers--Matet forcing with a finite sequence of pairwise non-isomorphic Ramsey ultrafilters over $\omega$, and with this forcing we settle the long-standing problem about
the spectrum of numbers of near-coherence classes. We prove that for any finite $n \geq 1$, there is a forcing extension with exactly $n$ near-coherence classes of non-principal ultrafilters.

For evaluating the new forcing, we prove a strengthening of Gowers' theorem on colourings of $\Fk$.
\end{abstract}

\nothing{
\thanks{}
}

\email{heike.mildenberger@math.uni-freiburg.de}
\date{July 25, 2019}
\maketitle
\tableofcontents

\section{Introduction}

We connect two lines of research: The topic of whether special kinds of ultrafilters
from the ground model have extensions of the same kind in a forcing extension
by a Ramsey-theoretic forcing and the investigation of the possible numbers of
near-coherence classes. We answer Banakh's and Blass' question
\cite[Question 31]{BanakhBlass} on the finite part of 
near-coherence spectrum.

On the existence of special ultrafilters:
For ${}^\omega \omega$-bounding forcings and $P$-points or even Ramsey ultrafilters the following is known: Kunen \cite{kunen:some-pts} proved that
no Ramsey ultrafilter can be extended to a $P$-point after addition of any number of random reals at once, Shelah \cite{Sh:f} constructed a model with no $P$-points, recently Chodounsky and Guzm\'{a}n \cite{chodounskyguzman} proved that there are no $P$-points in the Silver model and that no $P$-point from the ground model can be extended in a Silver extension. On the other hand \cite{Zheng2017} and \cite{FernandezBreton_Hrusak} proved that there are \mtu s in the Sacks model, indeed, any \mtu\ from the ground model is preserved. 
Blass \cite{blass-topap2009} proved that the minimum and the maximum projection of a \mtu\ are Ramsey ultrafilters.

In the case of forcings adding an unbounded real, Ketonen proved that $\gd = \gc$ implies the existence of a $P$-point, Canjar \cite{canjar:gen-ex} proved the generic existence of Ramsey ultrafilters under $\cov{\meager} = \gc$ and Eisworth \cite{Eisworth}
proved the generic existence of a \mtu\  with the  Galvin--Glazer technique (see, e.g. \cite{HindmanStrauss}) under the same condition.
 We refer the reader to \cite{Blasshandbook} for
the definitions of the cardinal characteristics $\cov{\meager}$, $\gu$, $\gd$, etc. We write $\gc$ for $2^{\aleph_0}$.

Here we work with variants of Matet forcing \cite{Matet} that come from various constraints
on the reservoir of the pure components of the conditions. The full Matet forcing preserves any $P$-point from the ground model \cite[Theorem~4]{Eisworth} and destroys any Ramsey ultrafilter, since it adds an unbounded real. The (non-complete) subforcings
with pure parts from a \mtu\ have specific preservation properties; they destroy some $P$-points and preserve others, see \cite[Theorem~2.5]{Eisworth}.
We show that the reservoir \mtu\ $\cU$ can be extended to a new \mtu\ after forcing. A new technical ingredient is the work with names for diagonal constructions.

We prove a preservation theorem for countable support iterations
and show that there is a model of $\aleph_1 = \gu < \gd=\gc=\aleph_2$ with at least three
names of different near coherence classes of ultrafilters and a \mtu\ of character $\gc$.

We fix some $k \in \omega \setminus \{0\}$.
We introduce \mtu s that are related to Hindman's theorem
for colourings of $\Fk$ as \mtu s are related to Hindman's theorem, and
investigate Matet forcing with \kmtu s. We show that this generalisation gives
nothing new in point of near coherence classes.

Then we change the condensation order into the Tetris condensation order
$\strk$ (see Def.~\ref{collection_on_Tetris}) and introduce a variant of Matet forcing in which
particular projections of pure conditions are taken from prescribed
pairwise nnc Ramsey ultrafilters $\cR_{i,x}$, $i =1,\dots, k$, $x ={\min},{\max}$. To get the pure decision property and hence properness, we introduce a space
\begin{equation}\label{gist}
  \begin{split}
    (\Fk)^\omega(\bar{\cR}) = & \{\ba \in (\Fk)^\omega \such \forall (X_{i} \in \cR_{i,{\min}}, Y_i \in \cR_{i,{\max}})_{1 \leq i \leq k}\\
    &(\exists^{\rm min-unb} s \in \TFU_k(\ba))\\
    &
    (\forall i \in\{1,\dots, k\})(\min_i(s) \in X_i \wedge \max_i(s) \in Y_i)\}
\end{split}\end{equation}
and investigate
$((\Fk)^\omega(\bar{\cR}), \strk)$ (see Def.~\ref{PP})  and prove a new
generalisation of Gowers' theorem and of Blass' theorem \cite[Theorem 2.3]{Blass:ufs-hindman}.

\begin{theoremgowers} Let $k \geq 1$ and let $\bar{\cR}=\la \cR_{i,\min}, \cR_{i,\max}
  \such i = 1,\dots k\ra$ be a sequence  of pairwise non nearly coherent Ramsey ultrafilters.
  \begin{myrules}
\item[(1)]
  Any $\strk$-descending $\omega$-sequence of elements of $(\Fk)^\omega(\bar{\cR})$
  has a $\strkast$-lower bound in $(\Fk)^\omega(\bar{\cR})$.
\item[(2)]
Let $n \in \omega \setminus \{0\}$ and $\ba \in (\Fk)^\omega(\bar{\cR})$ and let $c $ be a colouring of
  $[\TFU_k(\ba)]^n_{<}$ into finitely many colours. Then there is a $\bb \sqsubseteq_k \ba$, $\bb \in (\Fk)^\omega(\bar{\cR})$ such that $[\TFU_k(\bb)]^n_{<}$ is $c$-monochromatic.
\end{myrules}
\end{theoremgowers}

As a corollary we get the following theorems about existence and
additional  structure:

\begin{theoremexistence}
For any pairwise nnc Ramsey ultrafilters $\cR_{i,x}$, $i = 1, \dots, k$, $
x = {\min}, {\max}$ under \CH\ or {\sf MA}($\sigma$-centred) there is a \gmtu\ $\cU$ such that for $i = 1\dots, k$,
$\minl_i(\cU) = \cR_{i,\min}$ and $\maxl_i(\cU) = \cR_{i,\max}$.
\end{theoremexistence}

Pairwise non-near-coherence and Ramseyness is necessary, however, more structure 
emerges in the cases $k \geq 2$: The statement about the
higher core filters $\Phi_{\geq i+1}$ is novel. For the definition of $\Phi_{\geq i+1}$ see Def.~\ref{higher_core}.

\begin{theoremstructure}
   For any \gmtu\ $\cU$ over $\Fk$ the following holds:
   \begin{myrules}
   \item[(1)]
     For each $1 \leq i  \leq k$, $\cR_{i,{\min}}$ and $\cR_{i,{\max}}$ are nnc and
     for $i< k$, $\cR_{i,{\min}}$ and $\cR_{i,{\max}}$  are nnc $\Phi_{\geq i+1}(\cU)$.
\item[(2)]
The projections $\minl_1(\cU)$, \dots , $\minl_k(\cU)$, $\maxl_k(\cU)$, \dots, $\maxl_1(\cU)$ are  pairwise non-nearly coherent Ramsey ultrafilters over $\omega$.
 \item[(3)] All cores $\Phi_i(\cU)$ are nearly coherent.
        \end{myrules}
\end{theoremstructure}

Now we turn to the second line of research: the number of near-coherence classes of (non-principal) ultrafilters over $\omega$, see Def.~\ref{near_coherence}.
In \cite{BsSh:242} a model with one near-coherence class is given.
Blass \cite{ncf1} showed that under $\gd \leq \gu$ there are $2^{\gu}$ near-coherence classes.
Banakh and Blass \cite{BanakhBlass}
showed: If the number of near-coherence classes is infinite then it is $2^{(2^\omega)}$.

One of our main results about one iterand is:

\begin{theorem1} Let $\cE$ be a $P$-point and $\cU$ be a \mtu\ with $\Phi(\cU) \not\leq_{\rm RB} \cE$. Then in the forcing extension by $\M(\cU)$
  the \mtu\ $\cU$ is destroyed and can be completed to a \mtu\ $\cU^{\rm ext} \supseteq \cU$ with
  $\Phi(\cU^{\rm ext}) \not\leq_{\rm RB} \cE$. 
\end{theorem1}

The notion of a $P$-point will be explained in this section, the core $\Phi(\cU)$ is defined in Def.~\ref{core}, the Rudin--Blass order $\leq_{\rm RB}$ is defined in Def.~\ref{RudinBlassordering}. \mtu s are defined in Def.~\ref{past_etc}(6), the forcing $\M(\cU)$ is defined in Def.~\ref{Matet_centred}. 
By \cite[Theorem 2.5]{Eisworth}, $\cE$ generates a $P$-point in the extension.

Theorem \ref{1} serves as a successor step in the forcing that is used
in the following theorem:

\begin{theorem2} Assume \CH.
  Then there is a countable support iteration of proper iterands $\bP = \la \bP_\alpha, \M(\cU_\beta) \such \beta < \omega_2, \alpha \leq \omega_2 \ra$ such that
  in the extension there at least three near-coherence classes of ultrafilters
   and there is a \mtu\ of character $\gc$.
  \end{theorem2}

We generalise from the set of blocks to $\Fk$ for some $k \geq 1$
and show that this does hardly change the situation.

We use the new Ramsey space $((\Fk)^\omega(\bar{\cR}), \strk)$ to define
Gowers--Matet forcing  $\GM_k(\bar{\cR})$ in a localisation $\bar{\cR}$ to a
$2k$-sequence $\bar{\cR}$ of non-isomorphic Ramsey ultrafilters and prove our main result:

\begin{theorem4}
  Assume \CH\ and let $k \geq 1$ and fix $2k$ pairwise nnc Ramsey ultrafilters
  $\bar{\cR}_0=\la \cR_{i,x,0} \such i =1, \dots, k, x = {\min},{\max}\ra$.\footnote{The last index in the ultrafilters indicates the iteration
    stage  $\bV^{\bP_{\alpha}}$. The ultrafilters grow, i.e., $\cR_{i,x,\alpha} = \cR_{i,x,\beta} \cap \bV^{\bP_{\alpha}}$ for $\alpha < \beta \leq \omega_2$.}
  Then there is a countable support iteration of proper iterands $\bP = \la \bP_\alpha, \GM_k(\bar{\cR_{\beta}}) \such \beta < \omega_2, \alpha \leq \omega_2 \ra$ such that
  in the extension there \emph{exactly} $2k+1$ near-coherence classes of ultrafilters. Namely, one class is represented by a $P$-point of character $\omega_1$ and
  $2k$ classes represented by Ramsey ultrafilters $\cR_{i,x,\omega_2}$, $i = 1, \dots, k$,
  $x = {\min}, {\max}$.
\end{theorem4}

A slight variation of the model gives exactly $2k$ near coherence classes.
We thus get the full finite near-coherence spectrum.

\begin{corollary4} For any $n\in \omega$, the statement ``there are exactly $n+1$ near-coherence classes of ultrafilters'' is consistent relative to {\sf ZFC}.
\end{corollary4}

By work of Mioduszewski our result has 
applications to analysis,
namely the number of composants of $\beta({\mathbb R}^+)-{\mathbb R}^+$ corresponds by
\cite{Mio1,Mio2} to the number of near-coherence classes of ultrafilters.
Blass \cite{ncf1} gives applications to cofinality classes of short 
non standard models of arithmetic, and to the decomposition of the 
ideal of compact linear operators on a Hilbert space into proper subideals.
His results on equivalent characterisations of indecomposability can be
translated to: For $n \geq 1$, There is a decomposition of the ideal of compact operators into
$n$ proper subideals such that the union of any two different of them
is the whole ideal, if and only if there
are exactly $n$ near-coherence classes.
The correspondence is defined in \cite{BlassWeiss}.

\medskip

In the remainder of the introduction we recall some definitions from the realm of near coherence and special ultrafilters over $\omega$.
For the cardinal invariants $\gd$ and $\cov{\meager}$ we refer the reader to \cite{Blasshandbook}.

Let $S$ be a countable set.
By a {\em filter over $S$} we mean a non-empty subset of ${\mathcal P}(S)$
that is closed under supersets and under finite intersections
and that does not contain the empty set.
 We call a filter over $S$ {\em non-principal} if it
contains all cofinite subsets of $S$. A $\subseteq$-maximal filter is an ultrafilter.          

For $B \subseteq \omega$ and $f \colon \omega\to\omega$,
we let $f[B] = \{f(b)\such b \in B\}$ and $f^{-1}[B] 
= \{n \such f(n) \in B\}$.
The set of all infinite subsets of $\omega$ is denoted by $\roth$.
For $\cB \subseteq {\mathcal P}(\omega)$ we let
$f(\cB)= \{ X \subseteq \omega \such f^{-1}[X]
\in \cB\}$. \footnote{If $f$ is surjective, then $f(\cB)$ can be written as
  $\{f[X] \such X \in \cB\}$. In any case,  $f(\cB)$ is contained in the set of supersets of
  members of the latter set.}
This double lifting is an important function from 
${\mathcal P}{\mathcal P}(\omega)$ into itself. In analysis the special case of 
$f$ being finite-to-one  (that means that the preimage of each natural number is finite) is particularly useful, see e.g., \cite{ncf2}.

From now on all filters over $\omega$ will be non-principal filters over
$\omega$, though we write only ``filter'' over $\omega$.
If $f\colon \omega\to \omega$ is finite-to-one,
then also $f(\cF)$ is a non-principal filter. It is the filter generated
by $\{f[X] \such X\in \cF\}$.
\begin{definition}\label{near_coherence}
  \begin{myrules}
  \item[(1)] A non-empty family $\cG \subseteq \roth$ is called a \emph{filter subbase (over $\omega$)} if
    any intersection of finitely many elements of $\cG$ is infinite.
    We write
    \[
    \filter(\cG) = \{X \in \roth \such (\exists n \in \omega)(\exists G_0 \dots \exists G_n \in \cG)( X \supseteq G_0 \cap \dots \cap G_n)\}
    \]
    for the filter generated by $\cG$. The \emph{character} of a filter $\cF$ is the smallest size of a generating subbase.
    \item[(2)]
Two filters ${\cF}, \cG \subseteq \roth$ are {\em nearly coherent}, if there is
some finite-to-one $f\colon \omega \to \omega$ such that $f({\cF} ) \cup
f({\cG})$ generates a filter. \nothing{We also say to this situation that
$f({\cF})$ and $f({\cG})$ are coherent.}

\item[(3)]
  On the set of non-principal ultrafilters near-coherence is an equivalence relation (for a proof see \cite{ncf1}, e.g.)
whose equivalence classes are called \emph{near-coherence classes}.

\item[(4)]
Two subsets $\cH_1$, $\cH_2$ of $\roth$ are called \emph{nnc}, if
for any $X_i \in \cH_i$, $i =1,2$ and any finite-to-one $h$
there is $Y_i \subseteq X_i$, $Y_i \in \cH_i$, $i = 1,2$ such that
$h[Y_1] \cap h[Y_2] =\emptyset$.

\item[(5)] $\cH \subseteq \roth$ is called \emph{nowhere almost a filter}
  if for any $X \in \cH$, and any finite-to-one $h$ the
  set $\{h[Y] \such Y \in \cH \wedge Y \subseteq X \}$ does not have
  the finite intersection property.
  \end{myrules}
\end{definition}
For filters, nnc is the negation of  near coherence.
Near coherence is witnessed by a weakly increasing surjective finite-to-one
function. $f$ is weakly increasing if $x<y \rightarrow f(x)\leq f(y)$.
A coideal $\cH$ is nowhere almost a filter if and only if it is nowhere
almost an ultrafilter. We will use the property for
coideals.

We say ``$A$ is almost a subset of $B$'' 
and write $A \subseteq^* B$ if $A \smallsetminus B$ is finite.
Similarly, the symbol  $=^*$ denotes equality
up to finitely many exceptions between elements of 
$[S]^\omega$ for a set $S$.

Let $\kappa $ be a regular uncountable cardinal.
An ultrafilter $\cW$ is called a {\em $P_\kappa$-point} if for every
$\gamma <\kappa$, for every $A_i \in \cU$, $i <\gamma$,
there is some $A \in \cW$ such that  for any $i<\gamma$,
$A \subseteq^* A_i$; such an
$A$ is called a {\em pseudo-intersection} 
of  $\{A_i \such i<\gamma\}$.
A $P_{\aleph_1}$-point is called a $P$-point.

Let $\bP$ be a notion of forcing.
We say that \emph{$\bP$ preserves an ultrafilter $\cW$ over $I$} if 
\[\Vdash_\bP
\mbox{``}(\forall X \subseteq I)(\exists Y \in \cW) (Y \subseteq X
\vee Y \subseteq I \smallsetminus X)\mbox{''}
\] 
and in the
contrary case we
say ``$\bP$ destroys $\cW$''.
In the first case $\{ X\in \roth \cap \bV[G]\such
(\exists Y\in \cW) X\supseteq Y\}$ is an ultrafilter in $\bV[G]$
and $\cW$ generates an ultrafilter in $\bV[G]$. We just say:
$\cW$ is an ultrafilter in $\bV[G]$.
If $\bP$ is proper and preserves $\cW$ and $\cW$ is a $P$-point in the ground model, then $\cW$
stays a $P$-point in the forcing extension by \cite[Lemma~3.2]{BsSh:242}.

An ultrafilter $\cW$ over $\omega$ is called a \emph{$Q$-point}
if for every strictly increasing sequence $\la n_i \such i \in \omega\ra$ of natural numbers
there is $X \in \cW$ such that for every $i$,
$|X\cap [n_i,n_{i+1})|\leq 1$.
  Any $Q$-point from the ground model ceases to be a $Q$-point after adding an unbounded real.

An ultrafilter $\cR$ is called
\emph{selective} (or \emph{Ramsey ultrafilter}) if it is a $P$-point and a $Q$-point.
We use the von Neumann natural numbers $n=\{0,\dots, n-1\}$.
We often use the following, equivalent (see \cite[Theorem 4.9]{booth:thesis}  characterisation of selectivity:
\begin{myrules}
\item[($\ast$)]
For any $\subseteq$-descending sequence  $\la A_n \such n \in \omega\ra$
of sets $A_n \in \cR$
there is $A\in \cR$ such that $A \subseteq A_0$ and
$(\forall n \in A)( A\setminus(n+1) \subseteq A_{n})$.
\end{myrules}

\nothing{
We remark that the weaker property of semiselectivity (see \cite[Def.~7.9]{Todorcevic:Ramsey}) would be sufficient for many of our colouring theorems. Since we are localising to ultrafilters and since semiselectivity
there coincides with selectivity we dispense with this generality.
{\sf either add more or leave it out}
}
\smallskip

Kunen \cite{kunen:some-pts} contains more information on Ramsey  ultrafilters. Two Ramsey ultrafilters are nnc iff they are not isomorphic \cite{Blass:ppts}.

\nothing{ 
A useful measure for growth is the $\leq^*$-relation:
Let $\baire$ denote the set of functions
from $\omega$ to $\omega$, and let
$f,g\in \baire$. We say $g$ eventually dominates $f$ and 
write $f\leq^*g$ if $(\exists n)(\forall k \geq n)
f(k)\leq g(k)$.  A family $\cD \subseteq \baire$ is called
a dominating family if $\forall f \in \baire \exists g\in \cD f \leq^*g$.

\begin{definition}\label{nextfunction}
Let $E\in \roth$.
The function $\nex(\cdot,E)\colon \omega\to\omega$
is defined by $\nex(n, E) = \min(E\cap [n,\infty))$.
\end{definition}

Near coherence is connected to smallness:

\begin{fact}\label{factaboutdom}(Proof of \cite[Theorem 3.2]{ncf2})
Let $\cV$, $\cW$ be non-principal filters over $\omega$.
$\cV$ is nearly coherent to $\cW$ iff
\[ \{\max(\nex(\cdot,X),\nex(\cdot,Y)) \such X \in \cV, Y\in \cW\}\]
is not a $\leq^*$-dominating family.
\end{fact}

Since this fact is crucial we give a proof.
Assume that $\cV$ and $\cW$ are nearly coherent. Let $f$ be a finite-to-one function such that  $f(\cV) \cup f(\cW)$ is a filter base.
W.l.o.g.\ we assume that there is a strictly increasing sequence of natural numbers $\la \pi_i \such i < \omega\ra$ such that $\pi_0 = 0$ and $f(n) =
i$ for $n \in [\pi_i,\pi_{i+1})$. We let $h(i) = \pi_{n+1}$ for $i \in [\pi_n,\pi_{n+1})$. We show that
 no member of $F$ dominates $h$. 
 Let $X \in \cV$, $Y \in \cV$. Then $f[X] \cap f[Y]$ is infinite.
 Let $i = f(x) = f(y)$ for some $x \in X$, $y \in Y$. Then $x,y \in [\pi_i,\pi_{i+1})$, and we assume that $x,y$ are minimal in $X \cap [\pi_i,\pi_{i+1})$
     and in $Y \cap [\pi_i,\pi_{i+1})$ respectively. Thus $h(\pi_i) = \pi_{i+1} >
       \max(x,y) = \max(\nex(\pi_i,X),\nex(\pi_i, Y))$.
       Hence we have $h \not\leq^* \max(\nex(\cdot,X),\nex(\cdot,Y))$.
       
       Now assume that $\cV$ and $\cW$ are not nearly coherent. Let $h \in {}^\omega \omega$ be given, w.l.o.g.\ we assume that $h$ is strictly increasing and $h(0)>0$.
       We let $\tilde{h}$ be the iterate of $h$: $\tilde{h}(0) = 0$,
       $\tilde{h}(n+1) = h(\tilde{h}(n))$.
       We let $f_e(n) = i$  for $n \in [\tilde{h}(2i), \tilde{h}(2i+2))$ and we let
         $\tilde{h}(-1) =0$ and
         $f_o(n)= i $ for  $n \in [\tilde{h}(2i-1), \tilde{h}(2i+1))$.
           Then there are $V_e, V_o \in \cV$ and $W_e, W_0 \in \cW$ such that
           $f_e[V_e] \cap f_e[W_e] = \emptyset$ and
           $f_o[V_o] \cap f_o[W_o] = \emptyset$. Since $\cV$ and $\cW$ are filters, we can assume $V_e = V_o$ and $W_e = W_o$.
           Then for any $n  \in \omega$ we
           have $\max(\nex(n,V_e),\nex(n,W_e))> h(n)$. \proofend

\smallskip
         }
         
In forcing, we follow the Kunen style that the stronger condition is the 
{\em smaller} one. This corresponds  to the close relationship
between the $\leq$-relation in Matet forcing  and the condensation relation $\sqsubseteq$  on the
second components, the so-called pure parts, of a condition in Matet forcing.

Readers who want to focus on our extension of Gowers' theorem in Theorem~\ref{gowersown1}, that
does not contain any forcing,
can just read the notational part from Def.~\ref{ksetting1} to Def.~\ref{ksetting2}  and then
the new Ramsey-theoretic work from the beginning of Section~\ref{S6}  through
Question~\ref{baumgartnerstyle}.

Biggest thanks go to the referee for careful reading and detecting
some gaps and errors in an earlier version.
I am grateful for Andreas Blass for explaining me his proof of \cite[Theorem~2.1]{Blass:ufs-hindman}.

\section{Matet forcing with \mtu s}
\label{S2}

In this section we review results of Blass, Eisworth, and Hindman and carry it
a bit further.
Our nomenclature follows Blass \cite{Blass:ufs-hindman},
Eisworth \cite{Eisworth} and Todor\v{c}evi\'{c}~\cite{Todorcevic:Ramsey}.

\begin{definition}\label{2.1}\hfill
\begin{myrules}
\item[(1)]
We let $\F$ denote the set of finite non-empty subsets of $\omega$.
\item[(2)]
An element $a\in \F$ is called a block.
\item[(3)]
For $a, b \in \F$ we write $a< b$ if $(\forall n \in a) (\forall m \in b) (n <m)$.
\item[(4)]
  We define a well-order (of type $\omega$) $\leq_{\rm lex, \F}$ on the set $\F$ via  
  $a <_{\rm lex,\F} b$ if $\max(a) < \max(b)$ or ($\max(a) = \max(b)$ and
  $\min(a \triangle b) \in a$).\footnote{Note that this is not the usual lexicographic well-order; e.g., $\{0,1\} <_{{\rm lex}, \F} \{1\}$. The aim of this well-order is to define the $\sqsubseteq$-largest common lower bound of two $\sqsubseteq$-compatible elements of $(\F)^\omega$ by induction on the blocks. The relation $\sqsubseteq$ is defined in Def.~\ref{condensation}.}
\item[(5)]
A sequence $\bar{a}=\langle a_n \such n\in\omega\rangle$ of members of $\F$ is
called \emph{unmeshed} if for all $n$, $a_n<a_{n+1}$.
\item[(6)]
By $(\F)^\omega$ we denote the set of unmeshed sequences of members in $\F$.
\item[(7)] Let $a$, $b$ be blocks. We let $a\cup b$ be undefined 
unless $a<b$. Otherwise,  $a \cup b$ is defined as the union.

\item[(8)]

  A set $X \subseteq \F$ is called \emph{min-unbounded} if
  for any $n \in \omega$ there is some $x \in X$ with $\min(x) \geq n$.

  \item[(9)]
  $(\F,\cup)$  is a partial semigroup.
The associative partial binary operation $\cup$ lifts to $\gamma(\F)$, the space of min-unbounded ultrafilters over $\F$, as follows (and we write $\dotunion$ for the lifted operation):
\begin{equation*}
\begin{split}
\cU_1 \dotunion\, \cU_2 = 
\{X \subseteq \F & \such 
\{s \such \{t \such s \cup t \in X \} \in \cU_2\} \in \cU_1\}.
\end{split}
\end{equation*}
For details and history see \cite[Section~4.1]{HindmanStrauss}.
\item[(10)]

  If $X$ is a subset of  $\F$, we write $\FU(X)$ for the set of all
unions of finitely many  members of $X$. We write $\FU(\bar{a})$ instead of 
$\FU(\{a_n \such n \in \omega\})$. We call $X$ an $\FU$-set if there is some unmeshed $\ba$ such that  $X = \FU(\ba)$.

 \item[(11)]
   For $\emptyset \neq X \subseteq \F$ we let $\min_{<_{{\rm lex},\F}}(X)$ be the $\leq_{\rm lex,\F}$-element of $X$.  For $\ba \in (\F)^\omega$ we let
   $\min_{\F}(\ba) =
  \min_{<_{{\rm lex}, \F}}\{a_n \such n < \omega\}$,
  which is $a_0$.
\item[(12)]
A \emph{filter over  $\F$} is a non-empty subset of ${\mathcal P}(\F)$ that
is closed under intersections and supersets and does not contain the empty set.
\item[(13)]
 For $X\subseteq \F$,
 the set $(\FU(X))^\omega$ denotes the collection of all infinite unmeshed sequences in  $\FU(X)$. 
For $\ba \in (\F)^\omega$,
the set $(\FU(\ba))^\omega$ denotes the collection of all infinite unmeshed
sequences in  $\FU(\ba)$ (recall item (9)).
\item[(14)] For $Y \subseteq \F$ and $s \in \F$ we write $(Y \past s) $ for
  $\{u \in Y \such \max(s) < \min(u)\}$.
  For $\ba \in (\F)^\omega$ and $s \in \F$ we write $(\ba \past s)$ for
  $\la a_n \such n \geq n_0\ra$, where 
  $n_0 = \min\{n \such \max(s) < \min(a_n)\}$.
\end{myrules}
\end{definition}

Now the set of min-unbounded elements ${\mathcal P}(\F)$ is equipped with a partial order $\sqsubseteq$
that makes it into a topological Ramsey space in the sense of \cite{Todorcevic:Ramsey}.
We will work with the (closed) subspace of sets of the form $\FU(\ba)$. \footnote{Differences in behaviour of union ultrafilters
  that contain only min-unbounded sets, and
ordered-union ultrafilters that have a bases of unmeshed sets, are not yet known. Any known proof of Hindman's theorem results in an unmeshed sequence.}

\begin{definition}\label{condensation}
Given $X$ and $Y \subseteq \F$, we say that 
\emph{$Y$ is a condensation of $X$} and we write 
$Y \sqsubseteq X$ if $Y \subseteq \FU(X)$.   
We say $Y$ is \emph{almost a condensation of} $X$ and we write 
$Y \sqsubseteq^* X$ if there is an $n \in \omega $ such that 
$(Y\past \{n\})$ is a condensation of $X$.
\end{definition}

The ``blurred'' order  $\sqsubseteq^*$ is
an $<\omega_1$-complete preorder.

We use the relation $\sqsubseteq$ mainly for $Y =\rge(\bb)$ with
$\bb\in (\F)^\omega$, $X=\rge(\ba)$ for $\ba \in (\F)^\omega$, and then we write
$\bb \sqsubseteq \ba$ for $\rge(\bb) \sqsubseteq \rge(\ba)$, and analogously for $\sqsubseteq^*$.
\nothing{
We also call $\bb$ with  $\bb \sqsubseteq \ba$ a strengthening $\ba$ and we call
 $\bb$ with  $\bb \sqsubseteq^* \ba$ an almost strengthening $\ba$.
Strengthening can be imagined as the following procedure:
First we drop (possibly infinitely many) blocks from $\ba$
such that infinitely many blocks remain. Then we merge 
finite groups of adjacent blocks to one block of $\bb$. 
We use the verb ``to strengthen $\ba$''  as an abbreviation for
``to replace $\ba$  by an appropriate (almost) strengthening and call that
strengthening again $\ba$''. It  will be clear from the context which strengthening is meant.
}
\begin{definition}\label{compatible}
  Let $\ba, \bb \in (\F)^\omega$. We say $\ba$ and $\bb$ are compatible, if there is a $\bc \sqsubseteq\ba,\bb$. In the contrary case we write $\ba \perp \bb$. 
\end{definition}

\begin{lemma}\label{luxury1}
  If $\ba$ and $\bb$ are compatible, there is a weakest $\bc \in (\F)^\omega$ such that $\bc=\la c_n \such n < \omega \ra \sqsubseteq \bb, \ba$. Moreover,
 $\bc$ is given by the following procedure. 
  By induction on $n$ we define $c_n$ as follows
$c_{0} = \min_{<_{{\rm lex}, \F}}(\FU(\ba) \cap \FU(\bb))$, $c_{n+1} = \min_{<_{{\rm lex}, \F}}\{s \in \FU(\ba) \cap \FU(\bb) \such c_n < s\}$. 
\end{lemma}

\begin{proof}
  By definition $\bc \sqsubseteq \ba,\bb$.
  To see that $\bc$ is the largest witness, for a contradiction, we
  suppose that there is $\bd \sqsubseteq \ba, \bb$ such that
  $\bd \not\sqsubseteq \bc$. We take the
  first $n$ such that $d_n \neq c_n$.

  Since $\ba$ and $\bb$ are
  unmeshed sequences, there are $i_1< \dots < i_k$, and $j_1< \dots < j_\ell$  such that
  \begin{equation}\label{fixed}\begin{split}
    &c_n = a_{i_1} \cup \dots \cup a_{i_k}= b_{j_1} \cup \dots \cup b_{j_\ell}\end{split}\end{equation}
  and there are no common maxima of support
  in the sequence $(a_{i_1},a_{i_1+1}, \dots, a_{i_{k}-1})$ and in the sequence
  $(b_{j_1},b_{j_1+1} \dots, b_{j_{\ell}-1})
  $
  in at any pair $(a_\varrho, b_\sigma)$ with $\varrho \in [j_1, \dots, j_{\ell})$ and $\sigma \in [j_1,\dots, j_{l})$ and there are no $i<i_1$. $j< j_1$,
      and $a_i$, $b_j$ past $\max(\supp(c_{n-1})$ ($-1$ in the case $n=1$) such
      that $\min(\supp(a_i)) = \min(\supp(b_j))$.

  First case: $\max(\supp(d_n)) > \max(\supp(c_n))$.
  However, then there is a block $d'_n$ with $\max(\supp(d'_n)) = \max(\supp(c_n))$
  and there is $d{''}_n \in \FUU(\ba \past c_{n-1}) \cap \FUU(\bb \past c_{n})$ such that $d_n = d'_n + d''_n$ and $\bd$  could be replaced the
strictly  larger upper bound 
 \[\bd'= \la d_0, \dots, d_{n-1}, d'_n, d{''}_n , d_{n+1}, \dots \ra \sqsupset \bd.
 \]
 In the second and the  third case, applied to $d'_n$, we show that such a $\bd'\not\sqsubseteq \bc$ cannot exist
 and thus we get a contradiction.
 
 Second case: $\max(\supp(d_n)) < \max(\supp(c_n))$.
  Since $\ba$ and $\bb$ are
  unmeshed sequences, then $\bc$ could be broken up, similar to $\bd$ is the
  previous case. This possibility to break up contradicts the
  $<_{{\rm lex},\F}$-minimality of $c_n$.

Third case: $\max(\supp(d_n)) = \max(\supp(c_n))$.
First subcase $d_n <_{{\rm lex},\F} c_n$. This contradicts the choice of $c_n$.

Second subcase: $d_n >_{{\rm lex,\F}} c_n$. Then $d_n$ uses in its sum
not all the summands $a_{i_1}$, \dots, $a_{i_\ell}$ that combine to $c_n$ and
hence $d_n = 0$, because leaving off summands on one side or on both sides of
\eqref{fixed} is not possible by
the choice of $c_n$ and  since $\max(\supp(d_n)) = \max(\supp(c_n))$.
However, $d_n$ is a block and hence not zero. Thus also this case
ends in a contradiction.
\end{proof}

We write $\bc = \ba \wedge \bb$.

\begin{definition}\label{2.4}
A set non-empty subset $\cC \subseteq (\F)^\omega$ is called \emph{centred}, if 
for any finite $C \subseteq \cC$ there is  $\ba \in \cC$/ $\ba \in (\F)^\omega$ 
that is a condensation of any $\bar{c}\in C$
and if $\cC$ is closed under finite alterations 
i.e., if $\bd\in\cC$ and there are $n,m \in \omega$ such that
$\la d_m \such m \geq m_0\ra = \la e_n \such n \geq n_0\ra$ then $\be\in \cC$.
\end{definition}

We specialise $\cC$ further.

\begin{definition}\label{past_etc}\makebox{}
\begin{myrules}
\item[(1)] Let $\cF$ be a filter. A \emph{basis} $\cB$ of $\cF$ is a subset of $\cF$ such that $(\forall F \in \cF)(\exists B \in \cB)(B \subseteq F)$.
\item[(2)]
A non-principal filter $\cF$ over $\F$ is said to be a {\em union} filter if it 
has a basis of sets of the
form $\FU(X)$ for $X\subseteq \F$ such that the elements of $X$ are pairwise disjoint. Note $\FU(X)$ need not be an $\FU$-set.
\item[(3)]
  A non-principal filter $\cF$ over $\F$ is said to be an {\em min-unbounded} filter if $(\forall n \in \omega)(\exists X \in \cF)(\forall s \in X)(\min(s) > n)$.
  \item[(4)]
  A non-principal filter $\cF$ over $\F$ is said to be an {\em ordered-union} filter if it has a basis of sets of the form $\FU(\bar{d})$ for $\bar{d}\in (\F)^\omega$.

\item[(5)]
Let $\mu$ be an uncountable cardinal. A union filter 
is said to be {\em $(<\mu)$-stable} if,
 whenever it contains $\FU(X_\alpha)$ for $X_\alpha \subseteq \F$, 
$\alpha<\kappa$, for some $\kappa<\mu$, then it also contains 
 some $\FU(Y)$ for some $Y$ such that for each $\alpha <\kappa$
 there is $n_\alpha \in \omega$ with $(Y \past \{n_\alpha\}) \subseteq \FU(X_\alpha)$. Such an $Y$ is called a \emph{lower bound} of
$\{X_\alpha \such \alpha<\kappa\}$. For ``$<\omega_1$-stable'' we say ``stable''.
\item[(6)]
A stable ordered-union ultrafilter is also called a {\em \mtu}. 
\item[(7)] An ultrafilter is called \emph{idempotent} if
  $\cU\dotunion \cU=\cU$.
\end{myrules}
\end{definition}

Ordered-union ultrafilters need not exist, as their existence implies the existence of $Q$-points \cite{Blass:ufs-hindman} and there are models without 
$Q$-points \cite{Miller:q-pts}.
Even union ultrafilters need not exist:
 Blass  \cite[Theorem~38]{blass-topap2009}
showed that the existence of a union ultrafilter implies the existence of at
least two near-coherence classes of ultrafilters.
In \cite{BsSh:242} Blass and Shelah show that it is consistent
relative to {\sf ZFC} to have exactly one near-coherence class of non-principal ultrafilters.
Union ultrafilters are idempotent. Idempotent ultrafilters 
exist by the Ellis--Numakura Lemma \cite{Ellis:lemma, Namakura}.
With the help of Hindman's theorem one  shows that 
\CH\ or  Martin's Axiom for $\sigma$-centred posets and
$<2^\omega$ dense sets implies that (even $< 2^\omega$-) stable \mtu s 
 exist \cite{Blass:ufs-hindman}. We recall Hindman's theorem:

\begin{theorem}\label{hindmantheorem}
 (Hindman, \cite[Corollary 3.3]{hindman:sums})
If the set $\F$ is partitioned into finitely many pieces then there is a set 
$\bar{d} \in (\F)^\omega$ such that  $\FU(\bar{d})$ is included in one piece.
\end{theorem}

The theorem also holds if instead of $\F$ we partition some $\FU(\bar{c})$
for a $\bar{c} \in (\F)^\omega$ and search for a homogeneous sequence $\bar{d}
\sqsubseteq \bar{c}$, see \cite[p.~92]{Blass:ufs-hindman}.

\begin{corollary}\label{2.7}(See \cite[p.~93]{Blass:ufs-hindman}.)
Under \CH\ or ${\sf MA}_{<2^\omega}(\sigma\mbox{-centred})$, for every $\bar{a}\in (\F)^\omega$ there is a \mtu\  $\cU$ such that $\FU(\bar{a}) \in \cU$.
\end{corollary}

We let for $X\subseteq\F$, $[X]_<^n$ be the set of increasing
unmeshed $n$-sequences of members of $X$.
For the evaluation of our forcings, Taylor's theorem \cite{taylor} is 
utilised:

\begin{theorem}\label{taylor}(Taylor \cite{taylor}.)
Let $\cU$ be a \mtu, $n\in \omega$. Let $[\F]_<^n$ be partitioned into finitely many sets.
Then there is $A\in \cU$ such that $[\FU(A)]_<^n$ is monochromatic.
\end{theorem}

\begin{corollary}\label{diagonal1}(\cite[Cor.~1.3]{Eisworth})
  Existence of diagonal lower bounds in \\ $((\F)^\omega, \sqsubseteq)$.
  Let $\cU$ be a \mtu, and let  $\la A_n \such n \in \omega\ra$ be
  a $\sqsubseteq$-descending sequence of members of $\cU$.
  Then there is a $B \in \cU$ such $B \sqsubseteq A_0$ and
  \begin{equation}\label{def_diagonal}
    (\forall s \in B) ((B \past s) \mbox{ is a condensation of } A_{\max(s) +1}).
  \end{equation}
  Such a $B$ is called a \emph{diagonal lower bound of $\la A_n \such n \in \omega \ra$.}
  \end{corollary}

In Equation~\eqref{def_diagonal} we can equivalently let $s$ range over $\FU(B)$. If $B$ is a diagonal lower bound and $B' \sqsubseteq B$ then $B'$ is a diagonal lower bound as well.

\begin{definition}\label{core}
Let $\cH$ be a subset of ${\mathcal P}(\F)$.
\begin{myrules}
\item[(1)] 
The {\em core of $\cH$} is the set
$\Phi(\cH) \subseteq \roth$ such that 
\[
X \in \Phi(\cH) \mbox{ iff }  (\exists Y \in \cH)
( \bigcup Y \subseteq X).
\]
\item[(2)] 
The {\em minimum projection of $\cH$} is the set
\[\minl(\cH)=
\{{\rm min}[Y] \such Y\in \cH\},\] 
where 
\[
{\rm min}[Y]= \{{\rm min}(y) \such y \in Y\},
\]
and analogously we define the \emph{maximum projection} $\maxl(\cH)$.
\item[(3)] For $\ba\in (\F)^\omega$ we write
  $\min[\ba]= \{ \min(a_i) \such i \in \omega\}$.
  \item[(4)]  For $\ba\in (\F)^\omega$ we write
  $\set(\ba)= \bigcup\{a_i \such i \in \omega\}$.
\item[(5)]
  For $\cH \subseteq (\F)^\omega$ we let $\Phi(\cH) =
  \Phi(\{\rge(\ba) \such \ba \in \cH\})$ and
  $\minl(\cH) = \{\min[\ba] \such \ba \in \cH\}$.
\end{myrules}
\end{definition}

Note that $\min[\ba] = \min[\FU(\ba)]$ and $\bigcup \{a_n \such n \in \omega\} = \bigcup\FU(\ba)$.

If $\cC \subseteq(\F)^\omega$ is centred, then the core of $\cC$ is a filter over $\omega$.
For centred $\cC \subseteq (\F)^\omega$ and for filters $\cC$ over $\F$
the projection  $\minl(\cC)$ is a filter and 
$\minl(\cC), \maxl(\cC) \supseteq \Phi(\cC)$.
Blass  (\cite[3.6--3.9]{Blass:ufs-hindman} together with
\cite[Theorem~38]{blass-topap2009}) showed that for a \mtu\ $\cU$,
$\minl(\cU)$ and $\maxl(\cU)$ are nnc Ramsey ultrafilters.

The cores of centred systems are just filters over $\omega$.
Even if $\cU$ is an ultrafilter over $\F$,  for any finite-to-one $f$,
$f(\Phi(\cU))$ need not be an ultrafilter over $\omega$.
Blass  \cite[Theorem~38]{blass-topap2009}) showed that for union-ultrafilters
$\cU$,  for any finite-to-one $f$,
$f(\Phi(\cU))$ is not an ultrafilter because among its supersets there two nnc ultrafilters generated
by the $f$-images of the minimum projection and the maximum projection..

If $\cU$ is a  \mtu, then $\Phi(\cU)$  does not have a pseudointersection
(see \cite[Prop.~2.3]{Eisworth}) and also any finite-to-one image 
of $\Phi(\cU)$ does not have a pseudointersection by the same proof.
Hence, by Talagrand \cite{Talagrand}
 $\Phi(\cU)$ is not meagre.
Thus the filter dichotomy principle (see, e.g., \cite{Blasshandbook}) precludes the existence of a \mtu.

\begin{definition}\label{RudinBlassordering}
The {\em Rudin--Blass ordering} for filters over $\omega$ 
is defined as follows:
Let $\cF \leq_{RB} \cG$ if there is a finite-to-one $f$ such that  
$f(\cF) \subseteq f(\cG)$.\footnote{Also the definition $f(\cF) \subseteq \cG$ is used in the literature. If $\cG$ is a $P$-point both definitions
 are closely related.} 
\end{definition}

For filters $\cF$, $\cG$,  the relation $\cF \leq_{\rm RB} \cG$ implies that $\cF$ is nearly coherent to $\cG$.
If $\cG$ is an ultrafilter, also the converse holds.

Now we turn to forcing:

\begin{definition}\label{Matetforcing}
Conditions in {\em Matet forcing}, $\M$, are
pairs $(s, \bar{c})$ such that
$s\in \F$ and $\bar{c} \in (\F)^\omega$ and $s<c_0$.
The forcing order is $(t, \bar{d}) \leq (s,\bar{c})$ 
(recall the stronger condition is the smaller one) 
if $s \subseteq t$ and $t$
is the union of $s$ and finitely many  of the $c_n$ and
$\bar{d}$ is a condensation of $\bar{c}$.
\end{definition}

\begin{definition}\label{Matet_centred}
  For a family $\cH\subseteq (\F)^\omega$,
the notion of forcing $\M(\cH)$ 
consists of all pairs $(s,\bar{a})$
 such that $\bar{a}\in \cH$.
The forcing order is the same as in the Matet forcing.
\nothing{
  If we have that $s=t$, then we call 
$(t,B)$ a pure extension of $(s,A)$ and write $(t,B) \leq_{pr} (s,A)$.
In the special case that $\cC$ is the set of members of
a $\sqsubseteq^*$-descending sequence 
$\la \bar{c}^\eta  \such \eta < \beta\ra$,
 and their $=^*$ equivalent elements,
 we also write 
$\M(\bar{c}^\eta \such \eta<\beta)$ for $\M(\cC)$.} 
\end{definition}

We write $\set(\bar{a})=\bigcup\{a_n \such n \in \omega\}$.

For a centred system $\cC \subseteq (\Fk)^\omega$,
the set 
$\Phi(\cC)$ is the filter $\filter(\{\set(\ba) \such \ba \in \cC\})$.
The forcing $\M(\cC)$ diagonalises $\Phi(\cC)$. 
Let $G$ be a $\M(\cC)$-generic filter over $\bV$. Then 
the generic real 
\[
\mu := \bigcup\{s \such \exists \bc \such (s,\bc) \in G\}
\]
is a pseudointersection of $\Phi(\cC)$.

\nothing{
It is well known \cite{Matet, Blasstoronto} that  
Matet forcing $\M$ can be decomposed into 
two steps  $\M= \bP * \M({\name{\cU}})$, such that  
$\bP=((\F)^\omega, \sqsupseteq^*)$ is $<\omega_1$-closed
(that  is, every descending sequence of conditions 
of countable length has a lower bound) and adds a
stable ordered-union ultrafilter $\cU$
over the set $\F^k$. In particular $\M$ is proper.
\smallskip
Let $G$ be an $\M(\bar{a}^\alpha \such \alpha<\beta)$-generic filter over $\bV$.
So, in $\bV[G]$ the semifilter generated by
$\{\set(\bar{a}^\alpha)(i)\such \alpha <\beta\}$ is 
meagre even if it was not meagre before. 
This means that our forcings 
are very specific: They destroy the some $P$-points and preserves $\cE$.
}

The following property of \mtu s 
$\cU$ will be  important for our proof:

\begin{theorem}\label{Eisworth} (Eisworth 
\cite[``$\rightarrow$'' Cor.~2.5, ``$\leftarrow$'' Theorem~4, this direction works also
with non-$P$ ultrafilters]{Eisworth}) 
Let $\cU$ be a \mtu\ over $\F$
and let $\cW$ be
a $P$-point. Then $\cW \not\geq_{RB} \Phi(\cU)$ 
if and only if $\cW$ continues to generate an ultrafilter after we force with $\M(\cU)$.
\end{theorem} 

\nothing{
\smallskip

If a $\sqsubseteq^*$-descending
sequence $\la \bar{c}_\eps \such\eps < \omega_1\ra$ has the property that 
$\{\FU(\bar{c}_\eps) \such \eps<\omega_1\}$,
generates an
ultrafilter $\cU$ over $\F^k$, then $\cU$ is a stable ordered-union ultrafilter
and $\M({\cU})=\M(\bar{c}_\eps \such \eps<\omega_1)$.
Only this type of \mtu s appears in this paper along the construction.
} 

\smallskip

We remark:

\begin{proposition}\label{destroyforever}
Suppose that $\cU$ is a  \mtu\ and $\cF$ is a filter over $\omega$
 and $\Phi(\cU)\leq_{RB}
 \cF$. Then $\M(\cU)$
 forces that for any finite-to-one function from the ground model, 
$f(\cF)$ is not an ultrafilter.
\end{proposition}

\proof Let $f \in \bV$ be finite-to-one such that 
$f(\Phi(\cU))\subseteq f(\cF)$. 
Let $G$ be $\M(\cU)$-generic over $\bV$, and let $\name{\mu}$ be a name for
the generic real $\mu$.
We show: 
\begin{multline*}
(\forall (s,\bar{a}) \in \M(\cU) )(\forall Y\in \cF )
(\exists (t,\bar{b}) \leq_{\M(\cU)} (s,\bar{a})) 
 (t,\bar{b}) \Vdash_{\M(\cU)} f [Y] \cap f[\name{\mu}] \neq\emptyset
\end{multline*}
\begin{multline*}
(\forall (s,\bar{a}) \in \M(\cU) )(\forall Y\in \cF )
(\exists (t,\bar{b}) \leq_{\M(\cU)} (s,\bar{a})) 
(t,\bar{b}) \Vdash_{\M(\cU)} f [Y] \cap (\omega\setminus f[\name{\mu}]) \neq\emptyset.
 \end{multline*}

Let $(s,\ba)$ and $Y$ be given.
Since $f[\set(\ba)] \in f(\Phi(\cU)) \subseteq f(\cF)$, 
we have $f[\set(\ba)]\cap f[Y]$ is infinite.
So there is $t\in \FU(\ba)$ such that 
$f[t] \cap f[Y] \neq \emptyset$.
It follows that 
$(s\cup t, \ba\past t) \Vdash_{\M(\cU)} 
f [Y] \cap f[\name{\mu}] \neq\emptyset$.

Now for the second property:

Next we define a colouring $h$ of $[\FU(\ba\past s)]_<^2$ by
\begin{equation*}
h(u<v) = \left\{
\begin{array}{ll}
1 & \mbox{if }f[u] < f[v]
\wedge \\
& \bigl( \max(f[u]), \min(f[v])\bigr)
\cap f[Y]\neq \emptyset,\\
0 & \mbox{else.}
\end{array}
\right.
\end{equation*}
Since $\cU$ is a \mtu\, by Theorem~\cite{taylor} there is a monochromatic $\bar{b} 
\sqsubseteq (\ba \past s) $, $\bb\in \cU$.
Since $f$ is finite-to-one, the colour is 1.
So we have
$(s, \bb)\Vdash_{\M(\cU)} f [Y] \cap (\omega\setminus f[\name{\mu}]) \neq\emptyset$.
\proofend

Now let $\name{f}$ be a name for the function $n \mapsto |\mu \cap n|$.
Then
\[\M(\cU) \Vdash \name{f} \mbox{ is finite-to-one and }\name{f}[\mu] =^* \omega\]
and the above proof breaks down. Information on $\name{g}(\cF)$ for particular filters  $\cF$ with $\Phi(\cU) \leq_{\rm RB} \cF$ and any name $\name{g}$ for a finite-to-one function is
contained in Theorem~\ref{vorstufe}.

\begin{definition}\label{Fplus}
Let $\cF$ be a filter.
$\cF^+=\{X\in [\omega]^\omega \such (\forall Y\in \cF)(
Y\cap X \neq\emptyset)\}$.
\end{definition}

For a non-principal $\cF$, $\cF^+$ coincides with
$\{ X \in \roth \such (\forall Y\in \cF)( Y\cap X \in \roth)\}$.

\section{Ramsey-theoretic computations in $\M(\cU)$-extensions}
\label{S3}

Now we consider the Ramsey space $((\F)^\omega, \sqsubseteq)$ in an
$\M(\cU)$-extension.
We use the notation $\bV^{\bP}$ for any $\bV[G]$, with a $\bP$-generic filter $G$ over $\bV$. All computations are about names and we use the forcing theorem freely and identify names often with their evaluations.
The first aim is to examine
$\cU$ in the forcing extension $\bV^{\M(\cU)}$.
Is there an extension of the destroyed \mtu\ $\cU$ to a new
\mtu?
We show that under \CH\ there are even $2^{\omega_1}$ possibilities with pairwise nnc cores.
This will be the successor step of an iteration of iterands of type
$\M(\cU_\alpha)$.
We introduce Ramsey-theoretic computations with names for elements
of $(\F)^\omega$ in order to establish the existence of names for
\mtu s. Although in the end, like in Cor.~\ref{diagonal1}, min-unbounded $\FU$-subsets of $\F$ are the elements of \mtu s, intermediate work
is better carried out with sequences $\ba \in (\F)^\omega$. 
We use letters $\ba, \bb, \ba \dots, A, B, \dots$ for elements of $(\F)^\omega$,
where capital letters are in particular used in work with sequences of sequences. Capital letters are also used for subsets of $\F$.

\begin{definition} 
  \begin{myrules} 
  \item[(1)] Let $\ba \in (\F)^\omega$ and $X \in \roth$.
    We let $\ba \rest X = \la a_n \such n \in \omega, a_n \subseteq X\ra$.
  Note, we do not take those $a_n$ with $a_n \cap X \neq \emptyset$ that are not subsets of $X$.
\item[(2)] Let $\cU \subseteq (\F)^\omega$ and $X \in \roth$.
  We use the restriction symbol also for subsets of $(\F)^\omega$ and let $\cU \rest X = \{\ba \rest X \such \ba \in \cU\})$.
  \end{myrules}
\end{definition}
  
\begin{lemma} \label{positive_remark}
  Now let $\mu$ be a name for the Matet generic real $\bigcup \{s \such \exists \ba (s, \ba) \in G\}$.
  \begin{myrules}
  \item[(1)] $\M(\cU) \Vdash \forall \ba \in \cU \ba \rest \mu \in (\F)^\omega$.
  \item[(2)] $\M(\cU) \Vdash \cU \rest \mu \mbox{ is an ordered-union filter.}$
  \end{myrules}
\end{lemma}

\proof (1) is an easy density argument. For (2), we use
\begin{equation*}
  \begin{split}
    \M(\cU) \Vdash  (\ba \wedge \bb) \rest \mu = (\ba \rest \mu) \wedge (\bb \rest \mu) \; \wedge \;
     \FU(\ba \rest \mu) = \FU(\ba) \rest \mu.
  \end{split}
\end{equation*}
\proofend

Now we restate and prove
\begin{theorem}\label{1}
Let $\cE$ be a $\P$-point and let $\cU$ be a \mtu\ such that $\Phi(\cU) \not\leq_{\rm RK} \cE$. In $\bV^{\M(\cU)}$ there
is a  \mtu\ $\cU^{\rm ext} \supseteq \cU \rest \mu$
such that $\Phi(\cU^{\rm ext}) \not\leq_{RK} \cE$.
\end{theorem}

In the remainder of this section we prove this theorem in an
iterable form, Theorem \ref{one-step}.

\begin{definition}\label{Matetadequatefamily}(See \cite[Def.~3.1]{Eisworth})
A set $\cH \subseteq (\F)^\omega$ is called a \emph{Matet-adequate family} if
the following hold:
\begin{myrules}
\item[(i)] $\cH$ is closed $\sqsubseteq^*$-upwards.
\item[(ii)]
  $\cH$ is countably closed, i.e., any $\sqsubseteq$-descending $\omega$-sequence of members of $\cH$ has a $\sqsubseteq^*$-lower bound in $\cH$.
\item[(iii)] $\cH$ has the \emph{Hindman property}: If $\ba \in \cH$ and $\FU(\ba)$ is partitioned into two pieces then there is some $\bb \sqsubseteq \ba$, $\bb \in \cH$ such that $\FU(\bb)$ is a subset of a single piece of the partition.
  \end{myrules} 
\end{definition}

\begin{definition}\label{diagonallowerbound}
  Let $\la \ba_n \such n < \omega \ra$ be $\sqsubseteq$-descending
  of elements $\ba_n \in (\F)^\omega$. A sequence
  $\bb \in (\F)^\omega$ is a \emph{diagonal lower bound of $\la \ba_n \such n < \omega \ra$} if \begin{equation}
    \label{def_diagonal1}
    (\forall s \in \FU(\bb) )((\bb \past s) \sqsubseteq \ba_{\max(s) +1}).
  \end{equation}
\end{definition}

It is equivalent to say $(\forall s \in \bb) )((\bb \past s) \sqsubseteq \ba_{\max(s) +1})$ in \eqref{def_diagonal1}. 

The Hindman property of $\cH$ together with the countable closure of $\cH$ implies the existence of diagonal lower bounds of sequences in $\cH$, see \cite[Cor.1.3]{Eisworth}, which is based on the deep theorem \cite[4.2]{Blass:ufs-hindman}.
We give an alternative proof in Lemma~\ref{selfstr}.
\nothing{Since $\FU(\bb \past s) = (\FU(\bb) \past s)$ any diagonal lower bound $\bb$ of $\la \ba_n \such n < \omega \ra$ fulfils $(\forall s \in \FU(\bb))(\FU(\bb) \past s \subseteq \FU(\ba_{\max(s) +1})$.
}

Matet-adequate families have better properties than stated in the definition;
this is similar to \cite[Theorem 4.2]{Blass:ufs-hindman}.

\begin{lemma} \label{selfstr}
  Let $\cH$ be a Matet-adequate family.
 \begin{myrules}
\item[(a)] Let $n\geq 1$.
      If $A \in \cH$ and $[\FUU(\ba)]^n_{<}$ is partitioned into two pieces then there is some $\bb \sqsubseteq \ba$, $\bb \in \cH$ such that $[\FUU_k(\bb)]^n_{<}$ is a subset of a single piece of the partition.
\item[(b)]
  $\cH$ contains for each descending sequence a diagonal lower bound.   
\end{myrules}
\end{lemma}

\proof (a) 
is proved in \cite[Theorem 4.2]{Blass:ufs-hindman}.
We write an alternative proof that also serves to show that ultraness is not used.
For simplicity we write the step from 1 to 2. We let $c \colon [\FUU_k(\ba)]^2_{<} \to r$  for some finite $r \geq 1$. We enumerate $\FUU_k(\ba)$ as $\la s_\ell \such \ell \in \omega \ra$ such that
for any $\ell$, all the $s_i$ with $\max(\supp(s_i)) < \max(\supp(s_\ell))$
have $i < \ell$. Let $\bc_{-1} = \ba$. For each $s_\ell \in \FUU_k(\ba)$ by induction hypothesis we may take a $\bc_\ell \sqsubseteq_k (\bc_{\ell-1} \past s_\ell)$ such that
\begin{align*}
  f_\ell \colon &  \FUU_k(\bc_{\ell-1} \past s_\ell) \to m\\
  &t  \mapsto   c(s_\ell,t)
\end{align*}
is monochromatic on $\FUU_k(\bc_\ell)$.
Let $\bb \sqsubseteq^*_k \bc_\ell$ for any $\ell < \omega$.
Now we take $g \colon \omega \to \omega$ with $g(0)=0$ and $g(\ell+1)> g(\ell) $ so large that
\begin{equation}\label{late}
  (\forall r \leq g(\ell))
\bigl((\bb \past \{g(\ell+1)\}) \sqsubseteq_k \bc_r\bigr).
\end{equation}
Next we let for $s \in \FUU_k(\bb)$, say $s = s_\ell$, $f(s) = c(s,t)$ for any $t \in  \FUU_k(\bc_\ell)\cap \FUU(\bb)$. By the induction hypothesis there is $\bc \sqsubseteq_k \bb$,
$\bc \in \cH$, such that $\FUU_k(\bc)$ is monochromatic under the colouring $f$.
Then we colour $s \in \FUU_k(\bc)$ with colour $j \in 2$ if 
$\max(\supp(s))$ is in $\bigcup\{[ g(2r+j), g(2r+j+1)) \such r \in \omega \}$,
  and take $\bd \sqsubseteq_k \bc$ such that $\FUU_k(\bd)$ is monochromatic
  for the latter colouring. Finally we take $\be \sqsubseteq_k \bd$ such that
  for each $r$ there is at most one $\ell$ such that $\supp(e_\ell) \cap [g(2r+j), g(2r+j+1)) \neq \emptyset$.
    The existence of such an $\be \in \cH$ follows  by \cite[Theorem 3.9]{Blass:ufs-hindman}. Then $[\FUU_k(\be)]^2_{<}$ is $c$-monochromatic.

    (b) Let $\la  \ba_n \such n \in \omega\ra$ be a $\sqsubseteq_k$-descending sequence
  of element of $\cH$ and let $\bc \in \cH$ be a $\sqsubseteq^*_k$-lower bound
  such that $\bc \sqsubseteq_k \ba_0$.
  We colour $\FUU_k(\bc)]^2_{<}$ via
\[c(u,v) = \begin{cases} 1 & \mbox{ if } v \in \FUU_k(\ba_{\max(\supp(u)) +1});\\
  0 & \mbox{ else.}\end{cases}
\]
Since $\bc$ is a $\sqsubseteq_k^*$-lower bound of $\la \ba_n \such n < \omega \ra$, any 
$\bd \sqsubseteq_k \bc$ such that $[\FUU_k(\bd)]^2_{<}$ is $c$-monochromatic
has colour 1 and hence is a diagonal lower bound.
\proofend

\begin{definition}\label{positive}
  Let $\cC \subseteq (\F)^\omega$ be centred.
  $\cC^+ = \{\ba \in (\F)^\omega \such \forall \bc \in \cC, \bc \not\perp \ba\}$
\end{definition}
We introduce an abbreviation:
 
\begin{definition}\label{avoidforMatet}
  Let $\cH \subseteq (\F)^\omega$ and let $\cE$ be a $P$-point. We say \emph{$\cH$ avoids $\cE$} if $\{\set(\ba) \such \ba \in \cH\}$ is nnc to $\cE$.
  For nnc see Definition~\ref{near_coherence}(4).
  \end{definition}

The technical core of the proof of Theorem \ref{1} is:

\begin{theorem}\label{one-step}
  After forcing with $\M(\cU)$, $(\cU \rest \mu)^+$ is a Matet-adequate family that for any finite-to-one $h$, for any $\ba \in (\cU\rest\mu)^+$,
  \[( \exists \bb^1, \bb^2 \sqsubseteq \ba)(\bb^0 \in (\cU\rest\mu)^+ \wedge \bb^1 \in (\cU\rest\mu)^+
  \wedge h[\set(\bb^1)] \cap h[(\set(\bb^2)]=\emptyset).\]
  \end{theorem}

By the latter property, $(\cU\rest\mu)^+$ avoids any ultrafilter from
the ground model, in particular $\cE$.
Once Theorem \ref{one-step} is proved, 
a routine downwards construction along $\omega_1$ (see e.g., \cite[Theorem 2.4]{Blass:ufs-hindman}) completes the proof of Theorem~\ref{1}:
Under \CH\, in $\bV^{\M(\cU)}$ there is a \mtu\ $\cU^{\rm ext} \supseteq \cU$ that avoids $\cE$.

Showing adequacy requires some technical work to evaluate the forcing.

\begin{definition}\label{pure}
$(t,\bb) \leq (s,\ba)$ is called a \emph{pure extension of $(s,\ba)$}
if $s=t$.
\end{definition}

\begin{lemma} (\cite[Lemma~2.6]{Eisworth})
$\M(\cU)$ has the pure decision property, that is, for 
any $\varphi$ in the forcing language for any 
$(s,\ba)\in \M(\cU)$, there is $\bb\in \cU$, $\bb \sqsubseteq \ba$
such that $(s,\bb)$ decides $\varphi$. \proofend
\end{lemma}

Eisworth introduced the notion of a neat condition for a name of a subset of $\omega$. 
We extend the notion of neatness for our purposes. \footnote{It can be used for names of subsets of $H(\omega)$ and names for any $\omega$-hierarchy of hereditary finite sets whose union is a subset $H(\omega)$. We define neatness by tailoring initial segments towards computations with diagonal lower bounds.}

\begin{definition}\label{neatness}
\begin{myrules}
\item[(1)] Let $\name{A}$ be a name for an infinite subset of $\F$
(that means, the weakest condition forces this).
We say $(s,\bb)$ is \emph{neat for $\name{A}$} if
\begin{equation*}\label{oneA}
\begin{split}
&(\forall t\in \FU(\bb))(\forall r \in \FU(\bb \past t))(\forall u\subseteq\max(r))\\
& (\forall r' \in \FU(\bb \past r))
  \bigl((s\cup t\cup r', ( \bb\past r'))\\
&
\mbox{ decides  $u \in \name{A}$ and the decision does not depend on }r'\bigr).
\end{split}
\end{equation*}
\nothing{\item[(2)]
Let $\la\name{A_j}\such j<\omega\ra$ be a sequence of names for infinite
subsets of $\omega$.
We say $(s,\bb)$ is \emph{neat for $\la\name{A_j}\such j<\omega\ra$} if
\begin{equation*}\label{manyA}
\begin{split}
&(\forall t\in \FU(\bb))
  (\forall r \in \FU(\bb\past t))(\forall i,j \leq\max(r))\\
  &(\forall r' \in \FU(\bb \past r))
\bigl((s\cup t\cup r', ( \bb\past r')) \\
&\mbox{ decides $i \in \name{A_j}$  and the decision does not depend on }r'\bigr).
\end{split}
\end{equation*}
}
\item[(2)]
Let $\name{h}$ be a name for a finite-to-one function such that $h(i) \leq i$.
We say $(s,\bb)$ is \emph{neat for $\name{h}$} if
\begin{equation*}\label{neatforh}
\begin{split}
&(\forall t\in \FU(\bb))
  (\forall r \in \FU(\bb\past t))(\forall i \leq\max(r))\\
  &(\forall r' \in \FU(\bb \past r))
\bigl((s\cup t\cup r', ( \bb\past r')) \\
&\mbox{ decides $\name{h}(i)$  and the decision does not depend on }r'\bigr).
\end{split}
\end{equation*}
\item[(3)]
Let $\name{c}$ be a name for a function  $\name{c} \colon \F \to \{0,1\}$.
We say $(s,\bb)$ is \emph{neat for $\name{c}$} if
\begin{equation*}\label{neatforc}
\begin{split}
&(\forall t\in \FU(\bb))
  (\forall r \in \FU(\bb\past t))(\forall u \subseteq \max(r))\\
  &(\forall r' \in \FU(\bb \past r))
\bigl((s\cup t\cup r', ( \bb\past r')) \\
&\mbox{ decides $\name{c}(u)$ and the decision does not depend on }r'\bigr).
\end{split}
\end{equation*}

\item[(4)]
  Let $\la\name{A_j}\such j<\omega\ra$ be a sequence of names for
  elements of $(\F)^\omega$.
We say $(s,\bb)$ is \emph{neat for $\la\name{A_j}\such j<\omega\ra$} if
\begin{equation*}\label{neatformanyA}
\begin{split}
&(\forall t\in \FU(\bb))
  (\forall r \in \FU(\bb\past t))(\forall j \leq\max(r))(\forall u \subseteq \max(r))\\
  &(\forall r' \in \FU(\bb \past r))
\bigl((s\cup t\cup r', ( \bb\past r')) \\
&\mbox{ decides $u \in \name{A_j}$ and the decision does not depend on }r'\bigr).
\end{split}
\end{equation*}
\end{myrules}
\end{definition}

\begin{lemma}\label{ex_neat}(\cite[Lemma~2.7, Lemma 2.8]{Eisworth})
$(s,\ba) \in \M(\cU)$. Let $(\name{X},\name{h},\name{c}, \la\name{A_j} \such j <\omega \ra)$ given such that the weakest condition forces: $\name{X}$ is a min-unbounded subset of $\F$, $\name{h}$ is a surjective weakly increasing finite-to-one function, $\name{c}$ is a name for a colouring, $\la\name{A_j} \such j <\omega \ra$ is a sequence of members  of $(\F)^\omega$.
Then there is $\bb \sqsubseteq \ba$
such that $(s,\bb) \in M(\cU)$ is neat for $\name{X}$, $\name{h}$, $\name{c}$, $\la\name{A_j}\such j<\omega\ra$.
\proofend
\end{lemma}

\begin{remark}\label{for_later} Since the proof does not use the fact that $\cU$ is a filter,
  Lemma~\ref{ex_neat} also holds for $\M(\cH)$ for any Matet-adequate family $\cH$.
\end{remark}

Now we prove Theorem~\ref{one-step}. It is obvious that the set of positive sets $(\cU\ \rest \mu)^+$ is upwards closed
in the $\sqsubseteq^*$-order.
Now we show that any $\sqsubseteq$-descending $\omega$-sequence has a  $\sqsubseteq^*$-lower bound. It is not harder to directly show that there is a diagonal lower bound.

We recall that $\min_{<,{\rm lex},\F}$ was defined in Definition~\ref{2.1}(11).
A $\bQ$-name is a set of the form $\tau = \{ \la \sigma, q\ra \such
\la \sigma, q\ra \in \tau\}$ with names $\sigma$ of lower rank.
For $x \in \bV$ we have the $\bQ$-name $\check{x} =
\{ \la \check{y},q\ra \such y \in x, q \in \bQ\}$. We drop the
$\check{x}$-sign.

The following technique is one of the cornerstones of
our forcing constructions and interesting for itself:

\begin{lemma}\label{diagonal2forF}(Existence of positive diagonal lower bounds)
  Let $\cU$ be a  \mtu, $\cE$ be a $P$-point, $\Phi(\cU) \not\leq_{\rm RB} \cE$. Let $\bQ = \M(\cU)$ and let $\mu$ be
  the name for the generic real.
Let $\name{\bar{X}}=\la \name{X_n} \such n \in \omega \ra$ be a sequence of 
$\bQ$-names for elements of $(\F)^\omega$ such that
\begin{align*}
  \bQ & \Vdash (\forall n \in \omega) (\name{X_n} \in (\cU \rest \mu)^+ \wedge
  \name{X_{n+1}} \sqsubseteq \name{X_n}).
\end{align*}
Then
\begin{equation}
  \begin{split}
    \name{D} = \bigl\{\la t,&(s,\ba)\ra \such   (s,\ba) \in \bQ \wedge (\exists k \in \omega)
    (\exists t_0 < t_1 < \dots < t_{k-1} \in [\F]^k_<)
    \\ & \bigl(t_{k-1} < t_k = t \wedge
     (s,\ba) \Vdash \mbox{``}t_0 = \min_{<_{{\rm lex}, \F}}(\name{X_{0}} \rest \mu)
    \wedge \\ &   
    \qquad \qquad \bigwedge_{i < k} (t_{i+1} = \min_{<_{{\rm lex}, \F}}((\name{X_{\max(t_i)+1}} \rest \mu) \past t_i))\mbox{''}\bigr)\bigr\}
\end{split}\end{equation}
fulfils
\begin{equation}
\bQ \Vdash \name{D} \in (\cU \rest \mu)^+ \wedge \name{D} \sqsubseteq \name{X_0} \wedge
(\forall t \in \name{D})((\name{D} \past t) \sqsubseteq \name{X_{\max(t)+1}}).
\end{equation}
\end{lemma}

\proof
Let $(s,\ba)$ be given such that
$(s,\ba)$ is neat for $\name{D}$ and for $\la \name{X_i \rest \mu} \such
i < \omega \ra$.
We show that there is $(s,\bb) \leq_{\bQ} (s,\ba)$ that forces
that $\name{D} \rest \mu \in (\cU \rest \mu)^+$.

Let $k \in \omega$ and $u \in \FU(\ba)$. We say that \emph{$u$ is good for $(s,\ba)$, $\bar{t}= (t_0, \dots, t_{k-1}) \in [\F]_<^k$, $\name{\bar{X}}$}, if
\[(s \cup u ,\ba \past u) \Vdash t_0 = \min_{<_{{\rm lex}, \F}}(\name{X_0} \rest \mu)
\wedge \bigwedge_{i < k-1} (t_{i+1} = \min_{<_{{\rm lex}, \F}}((\name{X_{\max(t_i)+1}} \rest \mu) \past t_i)).
\]

Note that
goodness requires $t_i \in \F$ and not just names for elements of $\F$.
We define a colouring of $[\FU(\ba)]^2_<$ as follows

\begin{equation}
  F(u<v) = \begin{cases} 1 & \mbox{if for any $\bar{t}$ such that}
    \\&
    u  \mbox{ is good for } (s,\ba), \bar{t}, \name{\bar{X}},\\
    & \mbox{there is a proper end extension }\bar{t'} \mbox{ of } \bar{t}\\
    &\mbox{ such that }
    u \cup v  \mbox{ is good for } (s,\ba), \bar{t'}, \name{\bar{X}};\\
    0 & \mbox{ else}.
\end{cases}
\end{equation}
By Taylor's theorem \ref{taylor} there is
$(s,\bb) \leq_\bQ (s,\ba)$ such that $F$ is monochromatic on
$[\FU(\bb)]^2_<$.
We argue that the monochromatic colour can only be 1:
It suffices to find $u <v \in \FU(\bb)$ such that $F(u,v) = 1$.
Suppose that $u \in \FU(\bb)$ is good for $(s,\ba)$, $\bar{t}=(t_0, \dots, t_{k-1})$, $\name{\bar{X}}$. If $k =0$, we let $X_{t_{k-1}+1} = X_0$. Then $(s \cup u, \bb) \Vdash (\name{X_{\max(t_{k-1}) +1}} \rest \mu) \past t_{k-1})  \in (\cU \rest \mu)^+$. So  
$(s \cup u, \bb) \Vdash \exists t \in (\name{X_{\max(t_{k-1}) +1}} \rest \mu) \past t_{k-1})$.
\nothing{By the forcing theorem there is a condition $q \leq (s\cup u, \bb)$
and
there is a name $\name{t}$ such that
$q \Vdash \name{t} =
\min_{<_{{\rm lex}, \F}}((\name{X_{\max(t_{k-1}) +1}} \rest \mu) \past t_{k-1})$.
Now we choose a generic filter $G$ that contains $q$
and go into a forcing extension $\bV[G]$.
The evaluation fulfils $\name{t}_G \in \F$.
By the forcing theorem there is a condition $q'\leq q$, 
and there are a $t' \in \F$ such that $q' \Vdash \name{t}=t'$.
However, since $(s,\ba)$ is neat for $(\name{X_{\max(t_{k-1}) +1}}\rest \mu)$, also $(s,\bb) \in G$ is neat for it. By neatness, there is $v \in (\bb \past u)$, for example we can take $v$ to be the next element of $\bb$ that starts after $\max{t'}$ such that  $(s \cup u \cup v, \bb \past v) \Vdash \name{t} = t' = \min_{<_{{\rm lex}, \F}}((\name{X_{\max{t_{k-1}} +1}} \rest \mu) \past t_{k-1})$. So $F(u,v)=1$.
 } 

Now we show
\begin{equation}\label{D_good}
  (s,\bb) \Vdash \name{D} \in (\cU\rest \mu)^+ \wedge \name{D} \sqsubseteq \name{X_0} \wedge
  (\forall t \in \FU(\name{D}))
  ((\name{D} \past t) \sqsubseteq \name{X_{\max(t) +1}}).
\end{equation}

The proof comes in four parts. First we show that $\name{D}$ is a min-unbounded subset of $\F$. Given any $k \in \omega$, the procedure above for
finding $v$, $t'$ is iterated  $k$ times, starting at stage $k = 0$ with $u = s$. Thus we find
$v_i$ and $t_i$ such that
\begin{equation*}\begin{split}
    (s \cup v_0, \bb \past v_0 ) & \Vdash t_0 = \min_{<_{{\rm lex}, \F}}(X_0 \rest \mu),\\
    (s \cup v_0 \cup v_1, \bb \past v_1 ) & \Vdash t_1 = \min_{<_{{\rm lex}, \F}}(X_{\max(t_0) +1} \rest \mu),\\
    \vdots\\
    (s \cup v_0 \cup \dots \cup v_k, \bb \past v_k ) & \Vdash t_k = \min_{<_{{\rm lex}, \F}}(X_{\max(t_{k-1}) +1} \rest \mu), \mbox{ and thus }\\
    (s \cup v_1 \cup \dots \cup v_k, \bb \past v_k ) & \Vdash
    t_0 < \dots < t_k \in \name{D}.
\end{split} \end{equation*}
Next we show $(s,\bb) \Vdash \name{D} \in (\cU \rest \mu)^+$.
Suppose for a contradiction, that $\bc \in \cU$  and $(s',\bb') \leq (s,\bb)$ and
$(s',\bb') \Vdash \FU(\bc) \cap \FU(\name{D}) = \emptyset$.
Then we take $\bd \sqsubseteq \bb', (\bc \past s')$, such that $\bd \in \cU$ and see that $(s', \bd) \Vdash (\name{D} \past s') \subseteq \FU(\bc\rest \mu)$.
Contradiction.

Since  $\bQ \Vdash X_{n+1} \sqsubseteq X_n$, by definition of $\name{D}$,
$(s,\bb) \Vdash \name{D} \subseteq X_0$.

For the last conjunctive clause in Equation~\eqref{D_good},
we work with the characterisation of diagonal lower bound that is given
immediately after Def.~\ref{diagonallowerbound}.
We suppose that $(s \cup v,\bb \past v) \Vdash t < t' \in \name{D}$.
Then by the definition of $\name{D}$
$(s \cup v,\bb \past v) \Vdash t' \in  \FU(\name{X}_{\max(t) +1} \rest \mu ; \past t)$.
 \proofend

 Recall $\set(\la a_n \such n < \omega \ra) = \bigcup \{a_n \such n < \omega\}$.

\begin{lemma}\label{successor_avoidE_forF}
  Let $(s,\ba) \Vdash_{\M(\cU)} \name{\bc} \in (\cU \rest \mu)^+ \wedge \name{h}  \mbox{ is finite-to-one, onto and monotone}$.
  Then there are  $E \in \cE$ and $(s,\bb) \leq (s,\ba)$ and $\name{\bd^j}$, $j =0,1$,  such that
  \begin{equation}
    \begin{split}
      (s,\bb) \Vdash & \bigwedge_{j=0,1} (\name{\bd}^j \in (\cU \rest \mu)^+ \wedge \name{\bd}^j \sqsubseteq \name{\bc}) \;\wedge\;\\
      & \name{h}[\set(\name{\bd^0})] \cap \name{h}[\set(\name{\bd^1}] = \emptyset \wedge \bigvee_{j = 0,1} \name{h}[E] \cap \name{h}[\set(\name{\bd^j}) = \emptyset.
  \end{split}\end{equation}
  \end{lemma}

\proof
We assume w.l.o.g.\ that $(s,\ba) $ is neat for $\bc$ and $h$. Every name for a positive set $\bc$ and every condition $p$ let us define
names $\name{\bd^j}$, $j =0,1$,
such that densely many $q$ below $p$ force
$\name{h}[\set(\bd^0)] \cap \name{h} [\set(\bd^1)] = \emptyset \wedge
\bigwedge_{j=0,1} (\bd^j \sqsubseteq \bc)$.
Again this is proved with a colouring.

Let $k < \omega$ and $u \in \FU(\ba)$.
We say that \emph{$u$ is good for $(s,\ba)$, $\bar{d}= (d_0, \dots, d_{k-1})$, $\name{\bc}$, and $\name{h}$} if
\begin{equation*}
 \begin{split} (s \cup u ,\ba \past u) \Vdash & d_0 = \min_{<_{{\rm lex}, \F}}(\name{\bc} \rest \mu)
   \wedge\\
   &\bigwedge_{i < k-1} (d_{i+1} = \min_{<_{{\rm lex}, \F}}\{d \in \FU(\name{\bc} \rest \mu)
   \such \name{h}[d_i] \cap \name{h}[d] = \emptyset\}).
 \end{split}
\end{equation*}

Again the $d_i$ are in the ground model, not names. This, though, is not important, since we do not use them as indices.
We define a colouring of $[\FU(\ba)]^2_<$ as follows

\begin{equation}
  F(u<v) = \begin{cases} 1 & \mbox{ if for any $\bar{d}$ such that}
    \\&
    u  \mbox{ is good for } (s,\ba), \bar{d}, \name{\bc}, \name{h}
    \\
    & \mbox{there is a proper end extension }\bar{d'} \mbox{ of } \bar{d}
    \\
    &\mbox{such that }
    u \cup v  \mbox{ is good for } (s,\ba), \bar{d'}, \name{\bc}, \name{h};
    \\
    0 & \mbox{ else}.
\end{cases}
\end{equation}
By Taylor's theorem \ref{taylor} there is
$(s,\bb) \leq_{\M(\cU)} (s,\ba)$ such that $F$ is monochromatic on
$[\FU(\bb)]^2_<$. Since for any $k$, $((\name{\bc} \rest \mu) \past k)$
is forced to be
in $(\cU\rest \mu)^+$, the monochromatic colour can only be 1.
Now \nothing{by the maximal principle (see \cite[Ch.~VII, Theorem 8.2]{Kunen} or
\cite[Ch.~I, Lemma~3.1]{Sh:f} ``existential completeness'')} there are names $\name{d_i}$ such that
\begin{equation}
  \begin{split}
    (s,\bb) \Vdash & \la \name{d_{2k+1}} \such k \in \omega\ra \in (\cU\rest \mu)^+ \wedge \\
    &
  \la \name{d_{2k}} \such k \in \omega\ra \in (\cU\rest \mu)^+ \wedge
  \\
  &
  \la \name{d_{2}} \such k \in \omega\ra \sqsubseteq \name{\bc} \wedge
   \\
  & \name{h}[\bigcup \{\name{d_{2k}} \such k \in \omega\}] \cap \name{h}[\bigcup \{\name{d_{2k+1}} \such k \in \omega\}] = \emptyset.
  \end{split}
\end{equation}

The first two conjunctive clauses are shown as in the proof of Equation~\eqref{D_good}. The last conjunctive clause follows from the new definition of goodness.
We let $(s,\bb) \Vdash \name{\bd}^j =
\la \name{d_{2k+j}} \such k < \omega\ra$ for $j =0,1$.
Since $\M(\cU)  \Vdash \name{h}(\cE)$ is an ultrafilter,
$(s,\bb)$ forces there are a $j = 0,1$ and an 
$E \in \cE$ such that
$\name{h}[E] \cap \name{h}[\set(\name{\bd}^j)] = \emptyset$.
\proofend

The next lemma is the most important step in the proof of
Theorem \ref{one-step}. Indeed, it includes again a proof that positive diagonal
lower bounds exist.
\begin{lemma}\label{Hindmansuccessor}
In $\bV^{\M(\cU)}$, $(\cU\rest \mu)^+$ has the Hindman property.
\end{lemma}

For the proof of this lemma, we adapt a theorem of Eisworth.

\begin{theorem} \cite[Theorem~5]{Eisworth} Let $\cF$ be an ordered-union filter generated by $< \cov{\mathcal B}$ sets and let $c$ be a partition of $\F$ into finite sets. Then there is an $\ba \in \cF^+$ such that $\FU(\ba)$ is included in one piece of the partition.
  \end{theorem}

At a crucial point in Eisworth's proof a Cohen real over an elementary submodel provides a name in a Galvin--Glazer framework.
We show that also a Matet-real and even an $\M(\cU)$-generic real can be used. 
We recall the Galvin--Glazer \cite{hindman-carbondale} technique.

\begin{definition}\label{gammaFIN}  We denote by $\gamma(\F)$ the set of min-unbounded (see Def.~\ref{past_etc}(3)) ultrafilters over $\F$.
This set is endowed with the topology generated by
  \[\{\{ \cU \in \gamma(\F) \such A \in \cU\} \such A\subseteq \F\}.\]
\end{definition}

The space $\gamma(\F)$  is a compact zero-dimensional Hausdorff space. With the named topology and the semigroup operation $\dotunion$ from Def.~\ref{2.1}(9), the semigroup $(\gamma( \F), \dotunion)$ is a right-topological semigroup. Details can be found in \cite{HindmanStrauss}.

\begin{lemma} (Ellis, \cite{Ellis:lemma})
Each compact  subsemigroup of $(\gamma(\F), \dotunion)$ contains an idempotent ultrafilter.
  \end{lemma}

Now we apply Ellis' lemma to $\{\cU \in \gamma(\F) \such \cU \supseteq \cF\}$ for a min-unbounded filter $\cF$.

  \begin{lemma}\label{Ellisappl}
    (\cite[Prop. 4.2]{Eisworth}) Let $\cF$ be a min-unbounded filter. There is an idempotent ultrafilter $\cU \in \gamma(\F)$ that extends $\cF$.
  \end{lemma}

  Now we prove Lemma \ref{Hindmansuccessor}.
  Let $p$ force that $\name{c}$ is  a partition of $\name{\bb_0} \in (\cU \rest \mu)^+$ into finitely many pieces and $\la \name{\bb_n} \such n < \omega\ra$ is  a $\sqsubseteq$-descending sequence of
  elements $\name{\bb_n} \in (\cU \rest \mu)^+$.
  By Lemma \ref{Ellisappl} there is an $\M(\cU)$-name $\name{\cU^i}$  such that
  \[ \M(\cU) \Vdash \name{\cU^i} \supseteq \filter((\cU \rest \mu) \cup \{\{\name{b_{n,m}} \such m \in \omega\}  \such n \in \omega\}) \wedge
  \name{\cU^i} \,\dotunion\, \name{\cU^i} = \name{\cU^i}.\]
  For $X \subseteq \F$ and $t \in \F$ we set
  \begin{equation*}
    X \ominus t = \{ s \such s\cup t \in X \}.
  \end{equation*}
  Since $\name{\cU^i}$ is forced to be idempotent,
  \[\M(\cU) \Vdash (\forall X \in \name{\cU^i})(\{t \such X \ominus t \in \name{\cU^i}\} \in \name{\cU^i}).
  \]
  \nothing{
    We define
 \nothing{in a tree of maximal antichains $T = \omega^\omega$,
  $\{ p_{n,t} \such t \in T_n\}$ such that for each $n$ and $t \in T_n$
  $\{ p_{n,t} \such t \in T_n\}$ is a maximal antichain below $p_{n-1, t \rest n-1}$
   with $p_{0, \emptyset} = 1_{\M(\cU)} = (\emptyset, \ba)$ for any $\ba \in \cU$,
 } 
    for $n \in \omega$ names $\name{X_n}$ and $\name{d_n}$ and $p_{n}= (s_{n},\ba_{n})$
    with the following rules:
  \begin{myrules}
  \item[(1)] $p_0 \Vdash X_0$ is the piece of the partition $c$ of $\FU(\bb)$  that is in $\cU^i$.
  \item[(2)] 
    $p_{n+1} \Vdash d_{n}$ is the $\leq_{\rm lex,\F}$-least element of
    \begin{equation*}
      \begin{split}
        \{d \in X_{n} \cap \FU(\{a_{{n}, k} \such k \in \omega\}) \cap \FU(\bb_n) & \such X_{n} \ominus d \in \cU^i \mbox{ and }\\
        &\min(d) > \max(d_{i}) \mbox{ for } i <n\}.
        \end{split}
    \end{equation*}
       \item[(3)]
  $p_{n+1} \Vdash  X_{n+1} = X_{n} \cap (X_{n} \ominus d_{n}).$
  \end{myrules}

  Now we have to ensure that the sequence $\la p_n \such n < \omega \ra$ has a lower bound.
  }
Now we use again the Milliken--Taylor trick. We assume that $(s,\ba)$ is neat for $\name{c}$, $\la \name{b_n} \such n < \omega \ra$.

Let $n \geq 1$.  We call \emph{$u \in \FU(\ba)$ good for $(s,\ba)$, $(\name{X_m},\name{d_m} \such m < n)$} 
if $(s \cup u, \ba \past u)$ forces the following statements:
\begin{myrules}
  \item[(1)] $\name{X_0}$ is the piece of the partition $c$ of $\FU(\bb)$  that is in $\name{\cU^i}$.
  \item[(2)] We let $d_{-1} = \{-1\}$. For any $0 \leq m < n$
    $\name{d_{m}}$ is the $\leq_{\rm lex,\F}$-least element of
    \begin{equation}\label{backwards}
      \begin{split}
        & \{d \in \name{X_{m}} \cap \FU(\{a_{k} \such k \in \omega\}\rest \mu) \cap \FU(\bb_{\max(d_{m-1}) + 1})  \such \\
        & \qquad \name{X_{m}} \ominus d \in \name{\cU^i} \mbox{ and }
        \min(d) > \max(\name{d_{m-1}})\}
        \end{split}
    \end{equation}
       \item[(3)] For any $0 \leq m < n-1$,
  $\name{X_{m+1}} = \name{X_{m}} \cap (\name{X_{m}} \ominus \name{d_{m}}).$
  \end{myrules}

Here we allow names. Only the natural numbers are meant to be pinned down.
 
  We colour $[\FU(\ba \past s)]^2_<$ as follows:

  \begin{equation}
  F(u<v) = \begin{cases} 1 & \mbox{ if for any $(\name{X_m},\name{d_m} \such m < n)$ such that}\\&
    u \mbox{ is good for } (s,\ba), (\name{X_m},\name{d_m} \such m < n),\\
    & \mbox{there is a proper end extension }(\name{X_m}, \name{d_m} \such m < n') \\
    &\mbox{such that }
    v \mbox{ is good for } (s,\ba),(\name{X_m}, \name{d_m} \such m < n');\\
    0 & \mbox{ else}.
\end{cases}
\end{equation}

  Then we find a monochromatic $\bb \in \cU$ with $(s,\bb) \leq (s,\ba)$.
    Since $\cU^i$ is idempotent and $\bb \in \cU \subseteq \cU^i$, the set in \eqref{backwards} is in $\cU^i$. Hence
  the monochromatic colour can only be 1.
  We let $\name{\be}$ be a name  such that $(s,\bb) \Vdash (\forall n)
  \la \name{e_0}, \dots, \name{e_{n-1}}\ra = \la \name{d_0}, \dots , \name{d_{n-1}}\ra$.
  
  The monochromaticity statement 
  \[\M(\cU) \Vdash \FU(\name{\be}) \subseteq X_0\] is proved literally as in Eisworth
  \cite[page 460]{Eisworth}.
\nothing{
  By the definition of $\Vdash \exists x \varphi(x)$ we thus have
  \begin{equation*}
    \begin{split}
      1 \Vdash_{\M(\cU)} & \mbox{``for any partition $c$ of $\FU(\bb_0)$} \\
        & \mbox{ and for any descending sequence $\bb_n$ of 
        $(\cU \rest \mu)$-positive elements} \\
      & \mbox{there is a $\be \sqsubseteq \bb_0$ that is a diagonal lower bound of $\bb_n$} \\
   &   \mbox{ and such that  $\FU(\be)$ lies entirely in one piece of the partition $c$.''}
    \end{split}
    \end{equation*}
}
By item (2) in the current definition of ``good'', the sequence $\name{\be}$ is a diagonal lower bound of $\la \name{\bb_n} \such n<\omega\ra$.
  Now we show that $\name{\be}$ is positive. For this we use the conditions on $\FU(\ba)$ in the goodness clause (2). 
  Suppose for a contradiction that $\name{\be}$ is not forced to be $(\cU \rest \mu)$-positive. Hence there is $q \in \M(\cU)$, $\bc \in \cU$  such that $q$ is neat for $\name{\be}$ and $\mu$, $q \leq (s,\bb)$, and
  \[q\Vdash_{\M(\cU)} \name{\FU(\be)} \cap \FU(\bc \rest \mu) = \emptyset.
  \]
 Since $\cU$ is a filter, we can assume $q = (t,\bc)$.
  We produce an extension $r$ of $q$ that forces the contrary.
  
 There is a minimal $m$ such that
  $q$ does not determine $e_{m}$. So $q$ determines
  $e_0=d_0$, \dots, $e_{m-1}=d_{m-1}$.
  
  Since $[\FU(\bc \past t)]^2_{<}$ has colour 1, 
  \begin{equation*}
    \begin{split}
      (t,\bc) \Vdash
      Y= &
      \{d \in \name{X_{m}} \cap \FU(\{c_{k} \such k \in \omega\} \rest \mu) \cap \FU(\bb_{\max(d_{m-1}) + 1}) \such \\
      &
      \name{X_{m}} \ominus d \in \cU^i \mbox{ and }
        \min(d) > \max(d_{m-1})\}
        \in \cU^i.
    \end{split}
  \end{equation*}
  
  Since $[\FU(\bc)]^2_{<}$ has colour 1 and since $q$ is neat for $\name{\be}$ and 
  $\mu$,
  there is $r\leq q$ of the form $(t \cup u, \bc \past u)$
  and there is $d \in \F$
\[r \Vdash  d \in \FU(\name{\be}) \cap \FU(\bc \rest\mu),\]
in contradiction to the assumption on $q$.
 \proofendof{\ref{Hindmansuccessor},\ref{one-step}, \ref{1}}

 Henceforth we drop the tildes underneath the names.

 Now we return to filters over $\omega$ and answer some instances of the
 question left open in the previous section: What happens to
 filters with $\Phi(\cU) \leq_{\rm RB} \cF$?

Mathias introduced the following notion 
under the name ``happy family'' \cite[Def.~0.1.]{mathias:happy}.
Louveau studied it in the special case of ultrafilters
\cite{Louveau:IJM}.
Todorcevic  \cite[Chapter~7]{Todorcevic:Ramsey}
uses the name  ``selective coideal'' for a happy family.

\begin{definition}\label{selectivecoideal}(See \cite[Def.~0.1.]{mathias:happy}, \cite[Def.~7.3]{Todorcevic:Ramsey})
A set $\cH \subseteq[\omega]^\omega$ is called a \emph{selective coideal/happy family} if
the following hold:
\begin{myrules}
\item[(i)] $\cI_\cH:={\mathcal  P}(\omega)\setminus \cH$ is an ideal that contains all 
singletons.
\item[(ii)] If $\la A_i\such i \in \omega\ra$ is a 
$\subseteq$-descending sequence of elements $A_i\in \cH$, then
there is $B \in \cH$ such that 
$B \subseteq A_0$ and $(\forall i \in B) B\setminus (i+1) \subseteq A_{i}$.
We call such a $B$ a \emph{diagonal lower bound of
  $\la A_i\such i \in \omega\ra$.}
\end{myrules} 
\end{definition}

We write $\cF_{\cH} = \{ \omega\setminus X\such X \in \cI_\cH\}$
for the filter that is dual to $\cI_{\cH}$.
Then $\cH$ coincides with the $\cF_\cH$-positive
sets, i.e.,
\[
\cH=\cF_\cH^+:=\{ X\in \roth \such (\forall Y\in \cF_\cH) (X\cap Y\neq\emptyset)\}.
\]

\begin{lemma}\label{ober-ultrafilter-eines-filters}
Let $\cF$ be a filter over $\omega$  and let $\cR$ be an ultrafilter 
over $\omega$. $\cR\subseteq \cF^{+}$ iff $\cR\supseteq \cF$.
\end{lemma}

So the forward implication, which will be invoked many times, uses that $\cR$ is ultra.

\begin{remark}\label{avoidsE}
  Let $\cH \subseteq \roth$ and let $\cE$ be a filter.
  $\cH$ and  $\cE$ are nnc iff for any finite-to-one function $f$ and
  $X \in \cH$ there are a $E \in \cE$ and
  a $Y \subseteq X$, such that  $Y \in \cH$ and  such that $f[E] \cap f[H] = \emptyset$.
  \end{remark}

The following theorem provides information on $\minl(\cU)$ and $\maxl(\cU)$.

  \begin{theorem}\label{vorstufe} 
  Assume \CH\ and that we force with $\M(\cU)$ for a  \mtu\ $\cU$,
  $\cR\in \{\minl(\cU),\maxl(\cU)\}$  and $\mu$ is the generic real.
  After forcing with $\M(\cU)$,
  the set of positive sets
\[
(\filter(\cR \cup \{\mu\}))^{+} = \{X \in ([\omega]^\omega)^{\bV^{\M(\cU)}} \such (\forall Y \in \cR) |X  \cap Y\cap \mu|=\omega\}
\]
is a happy family that is nowhere almost a filter in $\bV^{\M(\cU)}$
(and hence it is nnc to $\cE$) and by \cite[Prop. 011]{mathias:happy}
 is a Ramsey ultrafilter $\cR^{\rm ext} \supseteq \cR \cup \{\mu\}$ that is nnc to $\cE$,
  \end{theorem}
\proof
The theorem is proved like 
 Theorem~\ref{1}, however, it is much easier. Lemma \ref{diagonal2forF}, giving diagonal lower bounds, and Lemma \ref{successor_avoidE_forF}, showing the nnc-part, are adapted to names for $(\cR \cup \{\mu\})$-positive subsets of $\omega$. There are no new ideas.
  \proofend
\nothing{Beweis des Theorems Vorstufe}

  \nothing{
\begin{lemma}\label{diagonal2} Existence of positive names for diagonal lower bounds,
  Let $\cU$ be a \mtu, $\cE$ be a $P$-point, $\cR \in \{\min(\cU),\max(\cU)\}$,
  and $\Phi(\cU) \not\leq_{\rm RB} \cE$. Let $\bQ = \M(\cU)$ and let $\mu$ be
  the name for the generic real.
Let $\la X_n \such n \in \omega \ra$ be a sequence of 
$\bQ$-names for subsets of $\omega$ such that
\begin{align*}
  \bQ & \Vdash (\forall n \in \omega) \bigl(X_n \in (\filter(\cR \cup \{\mu\}))^+ \wedge
  X_{n+1} \subseteq X_n\bigr).
\end{align*}
Then
\begin{equation}
  \begin{split}
    \name{D} = & \{\la \check{n},(s,\ba)\ra \such 
    \exists n_0 < n_1 < \dots < n_k = n\\
    & (s,\ba) \Vdash \bigl(\check{n_0} = \min(\name{X_0} \cap \mu)
    \wedge  \bigwedge_{i < k} n_{i+1} = \min(\name{X_{n_i}} \cap \mu \cap (n_i,\infty))\bigr)\}
\end{split}\end{equation}
fulfils
\begin{equation}
  \begin{split}
    \bQ \Vdash & \name{D} \in (\filter(\cR \cup \{\mu\}))^+ \wedge \name{D} \subseteq \name{X_0} \cup \mu \wedge\\
    &
(\forall n \in \name{D})( \name{D} \setminus (n+1) \subseteq \name{X_{n}}).
\end{split}\end{equation}
\end{lemma}

\proof
We first prove that $\name{D}$ is positive.
Let $Y \in \cR$ be given and $(s,\ba)$ be given such that
$(s,\ba)$ is neat for $\name{D}$ and for $\name{\bar{X}}$.
We show that there is $(s,\bb) \leq_{\bQ} (s,\ba)$ that forces
that $\name{D} \cap Y \cap \mu$ is infinite.

Let $k \in \omega$ and $u \in \FU(\ba)$. We say that \emph{$u$ is good for $(s,\ba)$, $\bar{n}= (n_0, \dots, n_{k-1})$, $\name{\bar{X}}$}, if
\begin{equation*}
  \begin{split}
    &\bigl(u \in \FU(\ba \past s) \mbox{ and}\\
    &(s \cup u,\ba \past u) \Vdash \check{n_0} = \min(\name{X_0}\cap \mu)
\wedge \bigwedge_{i < k-1} \check{n_{i+1}} = \min(X_{\check{n_i}} \cap \mu \cap (\check{n_i},\infty))\bigr).
\end{split}\end{equation*}

We define a colouring of $[\FU(\ba \past s)]^2_<$ as follows

\begin{equation}
  F(u<v) = \begin{cases} 1 & \mbox{ if for any $\bar{n}$ such that}\\&
    u \mbox{ is good for }(s,\ba), \bar{n}, \name{\bar{X}},\\
    & \mbox{there is a proper end extension }\bar{m} \mbox{ of } \bar{n}\\
    &\mbox{such that }
    v \mbox{ is good for }(s,\ba), \bar{m}, \name{\bar{X}};\\
    0 & \mbox{ else}.
\end{cases}
\end{equation}
By Taylor's theorem \ref{taylor} there is
$(s,\bb) \leq_\bQ (s,\ba)$ such that $F$ is monochromatic on
    $[\FU(\bb)]^2_<$ and since for any $n$, $\name{X_{n +1}} \cap Y\cap \mu$ is forced to be
an infinite set the monochromatic colour can only be 1.
Now $(s,\bb)$ is as desired. For a detailed proof we refer to the proof of \ref{diagonal2forF}.
 \proofend
}
 
\nothing{
\begin{lemma}\label{successor_avoidE}
  Let $\cR \in \{\min(\cU),\max(\cU)\}$ and $\cE$ be a $P$-point such that $\Phi(\cU) \not\leq_{\rm RK} \cE$. Let $\name{Y}$ be a $\bQ$-name for a $\filter(\cR \cup \{\mu\})$-positive set and let
 $ \name{h}$ be a $\bQ$-name for a finite-to-one function and
  let $(s,\ba) \in \bQ$.
  Then there are  $E \in \cE$ and $(s,\bb) \leq (s,\ba)$ and $
  (s,\bb) \Vdash \name{Z} \in \cR^+$
  such that
  \begin{equation}
    (s,\bb) \Vdash \name{h}[E] \cap \name{h}[\name{Z} \cap \name{Y} \cap \mu] = \emptyset.
  \end{equation}
  \end{lemma}

\proof W.l.o.g we can assume that $(s,\ba)$ is neat for
$\name{Y}$ and $\name{h}$. 
An application of a colouring similar to above shows:
For any $p = (s,\ba)$,  $\name{h}$ for a finite-to-one function for any name $\name{X}$
for a set in $(\filter(\cR \cup \{\mu\}))^+$ there are $q \leq p$, $\name{X_1}$, $\name{X_2}$,
\[q \Vdash \name{X_i} \subseteq \name{X} \wedge \name{X_i} \in (\filter(\cR \cup \{\mu\}))^+ \wedge h[\name{X_0}] \cap h[\name{X_1}] = \emptyset.\]
Since $\cE$ is a filter,  there is a $j \in 2$ and an $E \in \cE$ such that
$q \Vdash \name{h}[E] \cap \name{h}[\name{X_j}\cap \mu]= \emptyset$. 
Again we refer the reader to the $\F$-version in \ref{successor_avoidE_forF}
for a detailed proof.
\proofend \proofendof{\ref{vorstufe}}
}

\begin{remark}
 We remark that  
 by an analogous proof to Mathias' \cite[Prop.~011]{mathias:happy},  under \CH\ any happy family that is nnc to $\cE$ contains a Ramsey ultrafilter
  as a subset that is nnc to $\cE$. So we see that instead of heading for a  \mtu\ $\cU^{\rm ext}\supseteq \cU \rest \mu$ that yields of course
$\minl(\cU^{\rm ext}) \supseteq \minl(\cU)\cup\{\mu\}$ and the same for the maximum projection, we could 
proceed into a different direction: By Theorem \ref{vorstufe}, we can extend  
the minimum and maximum projections $\minl(\cU)$, $\maxl(\cU)$
to new Ramsey ultrafilters in $\bV^{\M(\cU)}$
and not care whether these extensions are the minimum and the maximum of
a \mtu\ extending $\cU$.
This direction will be important in Section \ref{S6}.
\end{remark}
  
  \section{A name for a Matet-adequate family at limit stages}
\label{S4}

We define by induction on $\alpha \leq \omega_2$
a countable support iteration $\bP_\alpha = \la \bP_\beta, \M(\cU_\gamma) \such \beta \leq  \alpha, \gamma < \alpha \ra$ such that
for any $\gamma < \alpha$,

\begin{equation}\label{induc}
  \begin{split}
  &\bP_{\gamma+1} = \bP_\gamma \ast \M(\cU_\gamma) \mbox{ and }\\
    &\bP_{\gamma+1} \Vdash   \mu_\gamma = \bigcup\{s \such (s, \ba) \in G_{\M(\cU_\gamma)}\}\\
    &\bP_\gamma \Vdash \cU_\gamma \supseteq (\bigcup\{ \cU_\delta \rest \mu_\delta) \such \delta < \gamma) \\
    & \qquad \qquad \mbox{ is a \mtu\ that avoids $\cE$.}
\end{split}\end{equation}

In Theorem~\ref{1} we proved that there are extension of \mtu s in the successor steps
\begin{equation*}
  \begin{split}
    \bP_\gamma \Vdash_{\bP_\gamma} \Bigl(&\M(\cU_\gamma) \Vdash_{\M(\cU_\gamma)}  \mu_\gamma = \bigcup\{s \such (s, \ba) \in G_{\M(\cU_\gamma)}\}\\
    & \wedge \bigl(\exists \cU_{\gamma+1} \supseteq \cU_\gamma) \\
    & (\cU_{\gamma+1} \mbox{ is a \mtu\ and } \cU_{\gamma+1} \supseteq  \cU_\gamma \rest \mu_\gamma\bigr)\Bigr).
  \end{split}
\end{equation*}
This guarantees the continuation of our construction in the successor steps,
via $\cU_{\gamma+1} \supseteq \cU_\gamma \rest\mu_\gamma$. 
  
Now we consider limit steps $\alpha$.
If $\cf(\alpha) > \omega$, 
we can just take $\bP_{\alpha}  \vdash \cU_{\alpha}=\bigcup_{\gamma<\alpha}\cU_{\gamma}$ and 
the inductive hypotheses will be carried on, since in proper forcing
every real appears at a step of at most countable cofinality, with the only exception that for $\alpha= \aleph_2$ the \CH\ gets lost.
So we concentrate on the hard case, $\cf(\alpha)=\omega$.

\begin{theorem}\label{hardest}
  Suppose that $\bP_\beta, \cU_\beta$ are as in Equation \eqref{induc}, $\alpha < \omega$, $\cf(\alpha)=\omega$,
  $\bP_\alpha$ is the countable support limit of
  $\la \bP_\beta, \M(\cU_\beta) \such \beta < \alpha \ra$.
  In $\bV^{\bP_\alpha}$, the filter $\cE$ still generates a $P$-point and the set of positive sets 
\[
\bigl(\bigcup_{\gamma<\alpha} (\cU_{\gamma} \rest \mu_\gamma)\bigr)^{+}
\]
forms a Matet-adequate family such that 
for any $\ba \in (\cU\rest\mu)^+$,
  \[( \exists \bb^1, \bb^2 \sqsubseteq \ba)(\bb^0 \in (\cU\rest\mu)^+ \wedge \bb^1 \in (\cU\rest\mu)^+
  \wedge h[\set(\bb^1)] \cap h[(\set(\bb^2)]=\emptyset).\]
\end{theorem}

As in Theorem~\ref{one-step}, the latter implies avoidance of $\cE$.
Again \CH\ and a routine enumeration along $\omega_1$ gives the following.

\begin{corollary}\label{hardest1}
  Suppose that $\bP_\beta, \cU_\beta$ are as in Equation \eqref{induc} and
  $\bP_\alpha$ is the countable support limit of
  $\la \bP_\beta, \M(\cU_\beta) \such \beta < \alpha \ra$ and that \CH\ holds
  in $\bV^{\bP_\alpha}$.
Then 
\begin{equation*}\begin{split}
    \bP_{\alpha} \Vdash \exists \cU_\alpha \bigl( & \cU_\alpha \mbox{ is a  \mtu\ that avoids $\cE$, and }\\
    &
\cU_\alpha \supseteq \bigcup_{\gamma<\alpha} (\cU_{\gamma} \rest \mu_\gamma)\bigr).
\end{split}\end{equation*}
\end{corollary}

Now we prove Theorem \ref{hardest}.
Blass and Shelah \cite[Theorem 4.1]{BsSh:242} showed that in $\bV^{\bP_\alpha}$ the closure of $\cE$ under almost supersets is  a $P$-point.

For the new part, by induction we define an increasing sequence  
$\bar{R}= \la R_\gamma \such \gamma < \alpha \ra$ of relations $R_\gamma$
in $\bV^{\bP_\gamma}$ such that a property called
\begin{equation}\label{ind2}
  \mbox{``$\bP_\gamma$ is
    $R_\gamma$-preserving''}
\end{equation} 
is a notion we want to carry from $\gamma < \alpha$, $\cf(\gamma) < \omega_1$, to $\alpha$ in addition to the property \eqref{induc}
and properness in the inductive choice of
the iteration.
\nothing{\footnote{We could restrict to $\M \prec H(\chi)$
  with $N \cap \omega_1 \in \cS$ and thus get a
  $(\bar{R},\cS)$-version that would allow us to split the
  (collapses of the ) iteration
  stages into the ones in $\cS$ and the ones outside for
  some stationary set $\cS$ in $\omega_1$. As this modification for proper forcings is not well-known \cite{Sh:f, FischerFriedmanKhomskyy} e.g, we do not
  take up this additional complexity.}
}

Now we define the relation $R_\alpha$ for which we want to preserve
statements of  the form $(\forall f)(\exists \bg)(f R_\alpha \bg)$. The relation $R_\alpha$ will be a Borel relation
on the Baire space in $\bV^{\bP_\alpha}$
that contain complex parameters from the ground model, e.g. \mtu s and names for
\mtu s.

\begin{definition}\label{therelations} By induction on $\alpha \leq \omega_2$ we define the following
  relations.
  \begin{myrules}
    \item[(1)]
  We say that a $\bP_\alpha$-name $\ba$ for an element of $(\F)^\omega$ is \emph{$\alpha$-positive} if
  $1\Vdash_{\bP_\alpha} \ba \in (\bigcup\{ \cU_\gamma \rest \mu_\gamma \such \gamma < \alpha\})^+$.
  
\item[(2)] 
  Assume that $\la \cU_{\gamma}\such \gamma<\alpha\ra$
is an ascending sequence of \mtu s  $\cU_{\gamma}
\in \bV^{\bP_\gamma}$, 
such that $\bP_\gamma \Vdash \Phi(\cU_\gamma) \not\leq_{\rm RB} \cE$ and
$\forall \gamma < \delta < \alpha$, $\bP_{\delta} \Vdash \cU_\gamma \rest \mu_\gamma \subseteq \cU_{\delta}$.
We say $f R_{\alpha}\bg$
if the following holds in $\bV^{\bP_\alpha}$:
\begin{myrules}
\item[(a)] $f= (\bar{A}, h, c)$,
\item[(b)] 
$ \bar{A} =\la A_{\ell}\such \ell \in \omega\ra$
  is a $\sqsubseteq$-descending sequence of $\alpha$-positive sequences $A_\ell \in (\F)^\omega$,
\item[(c)] $h$ is finite-to-one,
\item[(d)]
  $c$ is a partition of $\FU(A_0)$.
\item[(e)] 
  For $j = 0,1$ we let $\bg^j := \la g_{2n+j} \such n \in \omega \ra$.
  Then
  \begin{myrules}
  \item[(i)] For $j =0,1$, $\bg^j$  is an $\alpha$-positive diagonal lower bound of $\bar{A}$.
  \item[(ii)]  For $j = 0,1$, $\FU(\bg^j)$ is in one piece of the partition $c$.
\item[(iii)] $h[\set(\bg^0)]\cap h[\set(\bg^1)] = \emptyset$.
  \end{myrules} 
\end{myrules} 
\end{myrules} 
\end{definition}
So $R_\alpha$ is a $\bP_\alpha$-name for a relation.

\begin{definition}\label{R_alpha_preserving}
  We say 
  $\bP_\alpha$ is $R_\alpha$-preserving if
  $\bP_\alpha$ is proper and
  \[\bP_\alpha \Vdash \forall f \in \dom(R_\alpha) \exists \bg (f R_\alpha\bg).
  \]
\end{definition}

There are two main differences to the known ``Case A'' of iteration theorem \cite[Ch XVIII]{Sh:f}, \cite{Goldstern93}:
For our $R_\alpha$, it is not the case
that for countably many tasks $f_n$, $n \in \omega$ there is one answer $\bg$, $E$
such that $\forall n f R_{\alpha, n } (\bg,E)$ where $R_{\alpha,n}$ is $R_{\alpha}$ up to mistakes before $n$.
There are $\sqsubseteq^*$-incompatible positive elements $\bg$.
Secondly,
not only the quests $f$ but also the answers $\bg$ are now from the forcing extension.
This differs from the traditional applications in the preservation of cardinal invariants,
see e.g. \cite{BJ}.
We do not write tildes below the $R_\alpha$'s.

\begin{lemma}\label{lemma_pres} Assume \CH, $\alpha< \omega_2$, $\cf(\alpha)<\omega_1$. If $\bP_\alpha$ is $R_\alpha$-preserving then
  \begin{equation*}
    \begin{split}
      \bP_\alpha \Vdash & (\bigcup \{\cU_\beta\rest\mu_\beta\such \beta < \alpha\})^+
      \mbox{  is a Matet-adequate family such that 
      }\\
      & (\forall\ba \in (\cU\rest\mu)^+)
      (\exists \bb^1, \bb^2 \sqsubseteq \ba)\\
      & (\bb^0 \in (\cU\rest\mu)^+ \wedge \bb^1 \in (\cU\rest\mu)^+
  \wedge h[\set(\bb^1)] \cap h[(\set(\bb^2)]=\emptyset).
    \end{split}
  \end{equation*}
    \end{lemma}
\proof This follows from Definition~\ref{R_alpha_preserving}.
       \proofend

Now we carry the preservation property upwards by induction.

\begin{lemma}\label{limit-preservation} Assume \CH, $\alpha< \omega_2$, and $\cf(\alpha)=\omega$. Let $\bP_\alpha$ be the countable support limit of $\bP_\beta$, $\beta < \alpha$.
  If for $\beta < \alpha$ such that $\cf(\beta)< \omega_1$, $\bP_\beta$ is $R_\beta$-preserving
  and for any $\beta < \alpha$ Equation \eqref{induc} holds then
  $\bP_\alpha$ is $R_\alpha$-preserving.
\end{lemma}

This lemma will proved with Lemma~\ref{limit}.
For definiteness, we can take $\chi  = (2^{|\bP_\alpha|})^+$.
Under
\CH, for $\alpha < \omega_2$, $|\bP_\alpha| \leq \aleph_1$
by \cite[page 96]{Sh:b}. So $\chi = (2^{\aleph_1})^+$ is sufficiently large.
The following lemma on the translation to countable elementary submodels
is well-known, see \cite[Theorem 2.11 and Ch.~XVIII]{Sh:f}.

\begin{lemma}\label{translation}
  The following are equivalent.
  \begin{myrules}
    \item[(1)]
      $\bP_\alpha$ is $R_\alpha$-preserving.
      \item[(2)]
  For all $M \prec H(\chi)$ such that $\bP_\alpha, \la \cU_\beta \such \beta<\alpha\ra, \cE \in M$ for all $\bP_\alpha$-names $f \in M$
  for all $p \in \bP_\alpha \cap M$,
  if $p \Vdash f \in \dom(R_\alpha)$ then there is an $(M,\bP_\alpha)$-generic $q \geq p$ and there is a $\bP_\alpha$-name $\bg \in M$  such that
  $q \Vdash f R_\alpha \bg$.
  \end{myrules}
  \end {lemma}

\nothing{
Compare with \cite[Ch.~XVIII, Theorem~3.6]{Sh:f}
and to the version \cite{MdSh:973}:
}

The property of Lemma \ref{translation}(2) is carried on by induction
in the following slightly stronger technical form that is suitable
for induction.

\begin{lemma}\label{limit} The induction lemma.
  Suppose \CH, $\xi < \zeta \leq\aleph_2$, and $\cf(\zeta)< \omega_1$.
  Let $M \prec H(\chi)$, $\zeta \in M$, $\bP_\zeta \in M$, $p \in \bP_\zeta \cap M$ and $q_0 \leq p \rest \xi$
  be $(M, \bP_\xi)$-generic. Let $f\in M$ be a $\bP_\zeta$-name
  such that \[
  p \Vdash_{\bP_\zeta} f= (\la A_n \such n \in \omega\ra, c, h)  \in \dom(R_\zeta).
  \]
  Then there is some $(M,\bP_\zeta)$-generic condition
  $q \leq p \rest [\xi,\eta) \cup q_0$ such that 
  $q \Vdash (\exists \bg)(fR_\zeta \bg)$.
\end{lemma}

\proof 
Moreover we get
$\dom(q_0)  \setminus \xi \subseteq \zeta \cap M$.
We go by induction on $\zeta$. For $\zeta=0$ there is
nothing to prove, for $\zeta$ successor a proof is included in the proof of Theorem~\ref{1}, namely in Lemmata \ref{diagonal2forF} , \ref{successor_avoidE_forF}, \ref{Hindmansuccessor}. So
let $\zeta$ be a limit.
We fix a strictly increasing sequence $\la \zeta_\ell \such \ell <\omega\ra$
with $\zeta_0= \xi$, $\zeta_\ell \in M$, and $\sup \zeta_\ell = \zeta$.
We also fix an enumeration $\la D_n \such n < \omega\ra$ of the dense subsets
of $\bP_\zeta$ that are elements of $M$.
We let $\cU^i$ be a $\bP_\zeta$-name in $M$ such that
\begin{equation}
  \begin{split}
    \bP_\zeta \Vdash &\cU^i \mbox{ is idempotent and }\\
    &\cU^i \supseteq \filter(\bigcup \{(\cU_\gamma \rest \mu_\gamma) \such \gamma < \zeta\} \cup \{A_n \such n \in \omega \}).
  \end{split}
\end{equation}
By \ref{Ellisappl} such a name $\cU^i$ exists.
We choose by induction on $n$, $q_n \in \bP_{\zeta_n}$,  a $\bP_{\zeta_n}$-name
for a $\bP_{\zeta_n,\zeta}$-name $d_n$ for
an element of $\F$,
a $\bP_{\zeta_n}$-name
for a $\bP_{\zeta_n,\zeta}$-name $X_n$ for
an element of $\cU^i$,
and a $\bP_{\zeta_n}$-name of a $p_n \in \bP_{\zeta_n,\zeta}$ with the following properties:
\begin{myrules1}
\item[(a)] $q_n \in \bP_{\zeta_n}$, $\dom(q_n) \setminus \xi \subseteq  M \cap
\zeta_n$, $q_{n+1}\rest \zeta_n = q_n$,

\item[(b)]
$q_n$ is $(M,\bP_{\zeta_n})$-generic,

\item[(c)]
$q_{n+1} \Vdash_{\bP_{\zeta_n}} p_{n+1} \in D_n\cap M \cap G$,

\item[(d)] $p_0 \rest \zeta_0 \geq q_0$ in $\bP_{\zeta_0}$,

\item[(e)] $q_{n+1} \Vdash_{\bP_{\zeta_{n+1}}} p_{n+1}
\rest[\zeta_n,\zeta_{n+1}) \leq p_{n} \rest[\zeta_n,\zeta_{n+1})$
(in $\bP_{\zeta_n,\zeta_{n+1}}$),

  \item[(f)]
$q_1$, $p_{1}$, $X_0$, $d_0$, $X_1$ is such that 
    \begin{equation*}\label{backwards1}
      \begin{split}
        q_{1} \Vdash_{\bP_{\zeta_1}} \Bigl(p_{1} \Vdash_{\bP_{\zeta_{1}, \zeta}} &
        \bigl(X_0 \mbox{ is the piece of the partition $c \rest \FU(A_0)$ that lies in $\cU^i$}\\ &
        \wedge d_{0} = \min_{\rm lex,\F}\{d \in X_{0} \cap \FU(A_0)  \such X_{0} \ominus d \in \cU^i\}
 \\
 &
 \wedge \;  X_{1} = X_0 \cap (X_0 \ominus d_0)\bigr)\Bigr)
      \end{split}
      \end{equation*}
Moreover $q_{1} \Vdash p_1 \Vdash d_0\in M$, $q_1 \Vdash p_{1} \in M \cap G \cap D_0$.
\item[(g)]    
$q_{n+1}$, $p_{n+1}$, $d_{n}$, $X_{n+1}$ is such that 
    \begin{equation*}\label{backwards2}
      \begin{split}
        q_{n+1} \Vdash_{\bP_{\zeta_{n+1}}} & \Bigl(p_{n+1} \in M \cap G \cap D_n \wedge \\
        &p_{n+1} \Vdash_{\bP_{\zeta_{n+1}, \zeta}}  \bigl(d_{n} = \min_{\rm lex,\F}\{d \in X_{n} \cap \FU(A_n)  \such X_{n} \ominus d \in \cU^i \mbox{ and }\\
        &\min(d) > \max(d_{n-1})\:
        \wedge \: h[d_{n-1}] \cap h[d] = \emptyset\}
 \\
 &
 \wedge \;    X_{n+1} = X_n \cap (X_n \ominus d_{n})\bigr)\Bigr).
      \end{split}
      \end{equation*}
    Since the name $d_{n+1}$ is defined from elements in $M$,
    $q_{n+1} \Vdash_{\bP_{\zeta_{n+1}}} p_{n+1} \Vdash_{\bP_{\zeta_{n+1}, \zeta}} d_{n+1} \in M$.
\end{myrules1}

Since $M \prec H(\chi)$, such a sequence exists by the induction hypothesis and the maximum principle.

    In the end we let $\bg$ such that
    \[q_n \Vdash_{\bP_{\zeta_{n}}} p_n \Vdash_{\bP_{\zeta_{n}, \zeta}} \bg \rest n = \la d_i \such i < n\ra\]
    and $q = \bigcup_{n<\omega}q_n$ and we let $\bg^j$ be a name
    such that $q \Vdash \bg^j = \la g_{2n+j} \such n < \omega \ra$.

Now it is easy to see that $q$, $\bg^j$ are as desired, i.e.,
$q \leq p$ is $(M, \bP_\zeta)$ generic and
\begin{equation}\label{done}q \Vdash f R_\alpha \bg.
\end{equation}
Let $G$ be $\bP_\alpha$-generic over $V$ with $q \in G$.
Since $M[G] \prec H(\chi)[G]$ (again, see \cite[Theorem 2.11]{Sh:f}) there are also names $\bg \in M$ as in Equation~\eqref{done}.
\proofendof{\ref{limit}, \ref{hardest}}

Putting Theorem~\ref{1} and Theorem~\ref{hardest} together yields the following.

\begin{theorem}\label{2} Let $\cE$ be a $P$-point and assume \CH .
  Then there is a countable support iteration of proper iterands $\bP = \la \bP_\alpha, \M(\cU_\beta) \such \beta < \omega_2, \alpha \leq \omega_2 \ra$ such that
  in the extension $\cE$ is a $P$-point, there at least three near-coherence classes of ultrafilters
  and there is a \mtu\ of character $\gd=\gc=\aleph_2$.
\end{theorem}

We note that Fern\'{a}ndez-Br\'eton \cite{Fernandez-Breton} recently built a
model in which, like in ours, $\cov{\meager}= \aleph_1< \gd = \aleph_2$
with a \mtu.

\section{Generalisation to $\Fk$}\label{S5}

In this section we generalise the results of
Sections~\ref{S2}, \ref{S3}, and \ref{S4} from $\F$ to $\Fk$ for $k \geq 1$.
We introduce \mtu s over $\Fk$ and Matet forcing for $\Fk$.

\begin{definition}\label{ksetting1}
  Let $k \in \omega \setminus \{0\}$ unless stated otherwise. 
  \begin{myrules}
  \item[(1)] For $p \colon \omega \to k+1$ we let $\supp(p) = \{ n \in \omega
    \such p(n) \neq 0\}$.
    \[\Fk = \{p \colon \omega \to k+1 \such \supp(p) \mbox{ is finite } \wedge
   k \in  \rge(p)\}.\]
    
\item[(2)]
  For $a, b \in \Fk$, we let $a < b$ denote $\supp(a) < \supp(b)$, i.e., $(\forall m \in \supp(a))\\( \forall n \in \supp(b))( m < n)$, see Def.~\ref{2.1}(3).
  A finite or infinite sequence $\la a_i \such i < m \leq \omega\ra$ of elements of $\Fk$ is in \emph{block-position}
  if for any $i< j<m$, $a_i< a_j$.
  The set  $(\Fk)^\omega$ is the set of $\omega$-sequences in block-position, also called block sequences.
  For $n \geq 1$, the set  $[\Fk]_{<}^n$ is the set of
  $n$-sequences in block-position over $\Fk$.

\item[(3)]
  We define a well-order (of type $\omega$) $\leq_{\rm lex, \Fk}$ on the set $\Fk$ via  
  $a <_{\rm lex,\Fk} b$ if $\max(\supp(a)) < \max(\supp(b))$ or ($\max(\supp(a)) = \max(\supp(b))$ and there is an $m$ such that $a \rest m = b \rest m$
  and $a(m) > b(m)$.
  For a non-empty set $X \subseteq \Fk$ we let  $\min_{<_{{\rm lex},\Fk}}(X)$ be the
  $\leq_{\rm lex, \Fk}$-least element of $X$, in parallel to Definition~\ref{2.1}(11).

\item[(4)] For $k \geq 1$, $a,b \in \Ftwo{[1,k]}$, we define the partial semigroup operation $+$ as follows: If  $\supp(a) < \supp(b)$, then $a +b \in \Fk$ is defined. We let $(a + b) (n) = a(n) +  b(n)$.  Otherwise $a+b$ is undefined.
Thus $a + b = a \rest \supp(a) \cup b \rest \supp(b) \cup 0 \rest (\omega \setminus (\supp(a) \cup \supp(b)))$.

\item[(5)]
  Let $B$ be a min-unbounded subset of $\Fk$, i.e., for any $n$, $B \cap \{s \in \Fk \such \{n\} < s\} \neq \emptyset$. We let
\begin{equation*}
\begin{split}
    \FUU_k(B) = & \{b_{n_0}+\dots  + b_{n_\ell} \such\\ 
    & \ell \in \omega, b_{n_i} \in B,
    b_{n_0} < \dots < b_{n_\ell}\}
\end{split}
\end{equation*}
be the partial subsemigroup of $\Fk$ generated by $B$.
We call $B$ an \emph{$\FUU_k$-set} if there is a sequence
$\ba \in (\Fk)^\omega$ such that $B = \FUU_k(\ba)$.

\item[(6)] The set $\gamma(\Fk)$ is the set of ultrafilters over $\Fk$ that
  contain all sets of the form $\{p \in \Fk \such \{n\}<p\}$,
  $n \in \omega$. 
  
\item[(7)] We lift $+$ to $(\bigcup_{k \geq 1}\gamma(\Fk))^2$ (so also mixed pairs are allowed)  via
  \[\cU \dotpl \cV =
   \Bigl\{ X \subseteq \Fk \such \bigl\{ s\such \{t \such s + t \in X\} \in \cV\bigr\} \in \cU\Bigr\}.
   \]
  
\item[(8)]  For $A \subseteq \Fk$, $s \in \Fk$ and $n \in \omega$  we let
  $(A \past s) = \{a \in A \such \max(\supp(s)) < \min(\supp(a))\}$.
  \nothing{and
  $(A \past n) = \{a \in A \such n < \min(\supp(a))\}$.}
For $\ba \in (\Fk)^\omega$, $s \in \Fk$ and $n \in \omega$  we let
$(\ba \past s) = \la a_m \such m \geq m_0\ra$, where
$m_0$ is the minimal $m$ such that $\max(\supp(s)) < \min(\supp(a_m))$.
\nothing{and
$(\ba \past n) = \la a_m \such m \geq m_0\ra$ where
$m_0$ is the minimal $m$ such that $n < \min(\supp(a_m))$.
}

\item[(9)] For min-unbounded sets $A, B \subseteq \Fk$ we let
  $B \sqsubseteq_k A$ if $B \subseteq \FUU_k(A)$.
  We say $B$ is a $k$-condensation of $A$.
  We let
  $B \sqsubseteq^*_k A$ if there is an $n \in \omega$  such that
  $(B \past \{n\}) \subseteq \FUU_k(A)$.
 For $\ba, \bb \in (\Fk)^\omega$ we let
 $\bb \sqsubseteq_k \ba$ if
 $\bb = \la b_n \such n < \omega \ra$ and $\{b_n \such n < \omega \} \subseteq \FUU_k(\{a_n \such n < \omega \})$, and $\bb \sqsubseteq^* \ba$ has the obvious meaning.

  \end{myrules}
\end{definition}

So $(\F_1,+) \cong (\F,\cup)$ via $p \mapsto \supp(p)$,
and $((\F_1)^\omega, \sqsubseteq_1) \cong ((\F)^\omega, \sqsubseteq)$
via $\la a_n \such m \in \omega \ra \mapsto \la \supp(a_n) \such n < \omega\ra$.
 For checking whether there is a min-unbounded $C \sqsubseteq_k A, B$,
a piecewise construction
recursively along $\omega$ with the aid of $<_{\rm {lex}, \Fk}$
as in the proof of Lemma \ref{diagonal2forF} yields $\sqsubseteq$-weakest
(i.e., largest)
common lower bounds. Since this is not entirely obvious in the $\Fk$-setting,
we state it explicitly.

\begin{lemma}\label{luxury2}
  Assume that $\ba, \bb \in (\Fk)^\omega $ and that there is
  a $\bc \in (\Fk)^\omega$ such that $\bc \sqsubseteq_k \ba, \bb$.
    Then there is a $\sqsubseteq_k$-largest witness $\bc=\la c_n \such n < \omega\ra$ to this, that can be found
  in the following recursive way:
  \begin{eqnarray*}
    c_0 &= & \min_{<_{{\rm lex},\Fk}}( \FUU_k(\ba) \cap \FUU_k(\bb)),
    \\
    c_{n+1} &= & \min_{<_{{\rm lex},\Fk}} (\FUU_k(\ba \past c_n) \cap \FUU_k(\bb \past c_n)).
    \end{eqnarray*}
\end{lemma}
The proof is like Lemma \ref{luxury1}.

We write $\bc = \ba \wedge \bb$ for the largest $\bc \sqsubseteq_k, \ba, \bb$.

Hindman's theorem is generalised to $(\Fk)^\omega$:

  \begin{theorem}\label{hindman_k}(Hindman \cite[Cor. 3.3]{hindman:sums}) Let $k \geq 1$ and  $\ba \in (\Fk)^\omega$. For every finite colouring of $\FUU_k(\ba)$ there is a infinite block sequence $\bb \sqsubseteq_k \ba$ such that the elements of $\FUU_k(\bb)$ are monochromatic.
    \end{theorem}
  
\begin{definition}\label{kadequatefamily}
A set $\cH \subseteq (\Fk)^\omega$ is called a \emph{adequate family} if
the following hold:
\begin{myrules}
\item[(i)] $\cH$ is closed $\sqsubseteq_k^*$-upwards.
\item[(ii)] $\cH$ is stable, i.e.,
any $\sqsubseteq_k$-descending $\omega$-sequence of members of $\cH$ has a $\sqsubseteq^*$ lower bound in $\cH$.
\item[(iii)] $\cH$ has the \emph{Hindman property}: If $\ba \in \cH$ and $\FUU_k(\ba)$ is partitioned into finitely many pieces then there is some $\bb \sqsubseteq_k \ba$, $\bb \in \cH$ such that $\FUU_k(\bb)$ is a subset of a single piece of the partition.
  \end{myrules} 
\end{definition}

Again Lemma~\ref{selfstr} applies.

\begin{definition}\label{ksetting2}
  \begin{myrules}
  \item[(1)]
    For $\ba, \bb \in (\Fk)^\omega$ we write $\ba \not\perp \bb$
    and say that $\ba, \bb$ are compatible, if there is $\bc \sqsubseteq_k \ba,\bb$. \nothing{Indeed, there is a $\sqsubseteq$-largest such $\bc$, which is denoted by $\ba \wedge \bb$.}
 \item[(2)]
A set $\cC \subseteq (\Fk)^\omega$ is called \emph{centred}, if 
for any finite $C \subseteq \cC$ there is  $\bar{a} \in \cC$ 
that is a generalised condensation of any $\bar{c}\in C$
and if $\cC$ is closed under finite alterations 
i.e., if $\bd\in\cC$ and $\bd =^* \be$ then $\be\in \cC$.
  \item[(3)]
  A non-principal filter $\cF$ on $\Fk$ is said to be an {\em ordered-union} filter if it has a basis of sets of the form $\FUU_k(\bar{d})$ for $\bar{d}\in (\Fk)^\omega$.
  \item[(4)]
    For a  $\sqsubseteq_k$-descending $\omega$-sequence of members of $(\Fk)^\omega$
    $\la \ba_n \such n \in \omega \ra$ and $\bb \in (\Fk)^\omega$ we say: $\bb$ is a \emph{diagonal lower bound} if $\bb \sqsubseteq \ba_0$ and
    \[ (\forall s \in \FUU_k(\bb)) (\bb \past s) \sqsubseteq_k \ba_{\max(\supp(s)) +1}.\]

\item[(5)]
Let $\mu$ be an uncountable cardinal. An ordered-union filter 
is said to be {\em $(<\mu)$-stable} if,
 whenever it contains $\FUU_k(\bb_\alpha)$ for $\ba_\alpha \in (\F)^\omega$, 
$\alpha<\kappa$, for some $\kappa<\mu$, then it also contains 
some $\FUU_k(\bb)$ for some $\bb$ such that for each $\alpha$
there is $n_\alpha$ with $(\bb \past \{n_\alpha\}) \sqsubseteq_k \bb_\alpha)$
for $\alpha<\kappa$. Such an $\bb$ is called a \emph{lower bound} of
$\{\bb_\alpha \such \alpha<\kappa\}$. For ``$<\omega_1$-stable'' we say ``stable''.
\item[(6)]
A stable ordered-union ultrafilter over $\Fk$ is  called a {\em \kmtu (over $\Fk$)}. 
\item[(7)] An ultrafilter is called \emph{idempotent} if
  $\cU\dotpl \cU=\cU$.
\end{myrules}
\end{definition}

In particular,
any \kmtu\ is Matet-adequate and Lemma~\ref{selfstr} holds.

\begin{definition}\label{kMatetforcing}
  Let $k \geq 1$.  Conditions {\em Matet forcing (over $\Fk$)}, $\M_k$,
  are pairs $(s, \bar{c})$ such that
$s\in \Fk$ and $\bar{c} \in (\Fk)^\omega$ and $\supp(s)<\supp(c_0)$.
The forcing order is $(t, \bar{d}) \leq_k (s,\bar{c})$ 
(recall the stronger condition is the smaller one) 
if the following holds:
\begin{myrules}
\item[(a)] $s=t$ or there are $n$, $i_0 < \dots < i_n \in \omega$, and  $j_r \in \{0, \dots, k-1\}$ for $r \leq n$ with at least one $j_r = 0$ such that $t$ is a sum in the sense of Def.~\ref{ksetting1}(5) of the form
  \[t = s + c_{i_0} + \dots + c_{i_n}.\]
    \item[(b)]
$\bd \sqsubseteq_k \bc$ (see Def.~\ref{ksetting1}(9)).
\end{myrules}
In the case of $s=t$ we call $(t,\bb)$ a \emph{pure extension of $(s,\ba)=p$}.
$s$ is called the \emph{trunk of $p$} and $\ba$ is called the \emph{pure part of $p$}.
\end{definition}

Since the trunk and the pure part are in block position and since the pure parts are block-sequences, 
$s \rest \supp(s) = t \rest\supp(s)$ is equivalent to $s \rest (\max(\supp(s))+1)  = t \rest (\max(\supp(s))+1)$.

\begin{definition}\label{kMatetwithH}
Given a Matet-adequate family $\cH \subseteq (\Fk)^\omega$,
the notion of forcing $\M_k(\cH)$ 
consists of all pairs $(s,\bar{a})$
 such that $\bar{a}\in \cH$.
The forcing order is the same as in the Matet forcing.
\nothing{
  If we have that $s=t$, then we call 
$(t,B)$ a pure extension of $(s,A)$ and write $(t,B) \leq_{pr} (s,A)$.
In the special case that $\cC$ is the set of members of
a $\sqsubseteq^*$-descending sequence 
$\la \bar{c}^\eta  \such \eta < \beta\ra$,
 and their $=^*$ equivalent elements,
 we also write 
$\M(\bar{c}^\eta \such \eta<\beta)$ for $\M(\cC)$.} 
\end{definition}

\begin{definition}\label{corek}
Let $\cH$ be a subset of ${\mathscr P}(\Fk)$, $X \subseteq \Fk$, $s \in \Fk$.
\begin{myrules}
\item[(1)] 
  For $1 \leq i \leq k$ we let $\set_i(s) = s^{-1}[\{i\}]$,
  $\set_i[X] = \bigcup\{\set_i(s) \such s \in X\}$. 
$\set[X] = \bigcup\{\supp(s) \such s \in X\}$.

  \item[(2)]   For $1 \leq i \leq k$ we let $\min_i(s) = \min(s^{-1}[\{i\}])$,
    $\min_i[X] = \{\min_i(s) \such s \in X\}$. We have the analogous
    notions for $\max$.

  \item[(3)]
  The {\em core of $\cH$} is the set
$\Phi(\cH) \subseteq \roth$ such that 
\[
X \in \Phi(\cH) \mbox{ iff }  (\exists Y \in \cH)
(\set[Y]  \subseteq X).
\]
\nothing{\item[(2)]
We write
and \[\set_j(\bar{a})=\bigcup\{a_n^{-1}[\{j\}] \such n \in \omega\},\]
for $j = 1, \dots k$.
}
\item[(4)] For $j = 1, \dots, k$,
the {\em core of $\cH$ at colour $j$} is the set
$\Phi_j(\cH) \subseteq \roth$ such that 
\[
X \in \Phi_j(\cH) \mbox{ iff }  (\exists Y \in \cH)
(\set_j[Y]  \subseteq X).
\]
\item[(5)] 
The {\em minimum projection of $\cH$ at colour $j$} is the set
\[\minl_j(\cH)=
\bigl\{\min_i[Y] \such Y\in \cH
\bigr\},\] 
and analogously we define  $\maxl_j(\cH)$, the \emph{maximum projection at colour $j$}.
\nothing{
\item[(4)] For $\cH \subseteq (\Fk)^\omega$ we let $\Phi(\cH) = \Phi(\{\rge(\ba) \such \ba \in \cH\})$
  and $\minl_j(\cH)=
\bigl\{ \{{\rm min}\{a_n^{-1}[\{j\}] \such n \in \omega\} \such \ba\in \cH
\bigr\}$.
}
\item[(6)] Let $\cE$ be an ultrafilter over $\omega$. We say \emph{$\cH$ avoids $\cE$} if $\Phi(\cH)$ is nnc to $\cE$.
  
\item[(7)] The $\M_k(\cH)$-generic function
  from $\omega$ to $k+1$ is
  \[\mu = \bigcup\{s \rest (\max(\supp(s))+1) \such (\exists \ba)( (s,\ba) \in G)\}.
  \]
A name for $\mu$ is 
\[\name{\mu} = \bigcup\{ \la (s \rest (\max(\supp(s))+1), (s,\ba)\ra \such
(s,\ba) \in \M_k(\cH)\}.
\]
Usually we do not write the tildes.
  \item[(8)]
 In addition we define the generic $i$-fibres for $i =1, \dots, k$ by letting
  \[\mu_i := \mu^{-1}[\{i\}].\]
\end{myrules}
\end{definition}

Again we have
$\{\min(y^{-1}[\{i\}]) \such y \in \FUU_k(\ba)\}= \{\min(a_n^{-1}[\{i\}]) \such n \in \omega\}$
and
$ \bigcup \{b^{-1}[\{i\}] \such b \in \FU(\ba)\} =   \bigcup \{a_n^{-1}[\{i\}] \such n \in \omega\}$.
Sometimes we identify each $\FUU_k$-set $\FUU_k(\ba)$ with its generating $\ba \in (\Fk)^\omega$.
Especially in the forcing notion $\M(\cH)$ we think of $\cH$ as a subset of $(\Fk)^\omega$
even when we insert an \kmtu\ $\cU$ for $\cH$.

\begin{remark}\label{diagonalisers}
In the case of a centred adequate family $\cC \subseteq (\Fk)^\omega$, a density argument shows
for any $j = 1, \dots, k$,
\[\M_k(\cC)  \Vdash
\mu_j \in \roth \: \wedge \; 
(\forall X \in \Phi_j(\cC)) (\mu_j \subseteq^* X).\]
\end{remark}

The proof of the following theorem is carried out as in the original.
It makes frequent use of  the good properties from Lemma~\ref{selfstr}.

\begin{theorem}\label{eisworthfork}
  (Generalisation of 
\cite[``$\leftarrow$'' Theorem~4, ``$\rightarrow$'' Cor.~2.5, this direction works also
with non-$P$ ultrafilters]{Eisworth}) 
Let $\cU$ be a \kmtu\ over $\Fk$
and let $\cW$ be
a $P$-point. $\cW \not\geq_{RB} \Phi(\cU)$ 
iff $\cW$ continues to generate an ultrafilter after we force with $\M_k(\cU)$.
\end{theorem} 

\nothing{
\smallskip

\begin{corollary}
  Under \CH, for every $\bar{a}\in (\Fk)^\omega$ there is a \kmtu\  $\cU$ such that $\FUU_k(\bar{a}) \in \cU$.

  Under \CH, for every Matet-adequate family $\cH$ there is a \kmtu\  $\cU$ such that $\cU \subseteq \cH$.
 \end{corollary}
} 

Now here is the $\Fk$-analogon to Theorem~\ref{1}.

\begin{theorem}\label{theorem1fork}
Let $\cE$ be a $P$-point and let $\cU$ be a \kmtu\ over $\Fk$ such that $\Phi(\cU) \not\leq_{\rm RK} \cE$. In $\bV^{\M_k(\cU)}$ there is a  \kmtu\ $\cU^{\rm ext} \supseteq \cU \rest \mu$ with $\Phi(\cU^{\rm ext}) \not\leq_{RK} \cE$.
\end{theorem}

We introduce some notions that allow us to adapt the proof of Theorem~\ref{1}
to a proof of Theorem~\ref{theorem1fork}.

\begin{definition}\label{positivek}
  Let $\cC \subseteq (\Fk)^\omega$ be centred.
  $\cC^+ = \{\ba \in (\Fk)^\omega \such \forall \bc \in \cC, \bc \not\perp \ba\}$
\end{definition}

\begin{definition}\label{restrictionfork} Let $f \colon \omega \to (k+1)$ and let $\ba \in (\Fk)^\omega$ and let $\cU \subseteq (\Fk)^\omega$.
    \begin{myrules}\item[(1)]
      We let $\ba \rest f = \la a_n \such a_n = f \rest \supp(a_n), n \in \omega \ra$.
      \item[(2)]
    We let $\cU \rest f = \{\ba \rest f \such \ba \in \cU\}$.
  \end{myrules}\end{definition}

Density arguments show:
\begin{lemma}\label{density} Let $\cU$ be a \kmtu\ over $\Fk$ and let
  $\mu$ be the generic function, see Def.~\ref{corek}(7). Then 
    \[1 \Vdash_{\M_k(\cU)} \mu \colon \omega \to (k+1).\]
    \[1 \Vdash_{\M_k(\cU)} (\forall \ba \in \cU)\ba \rest \mu \in (\Fk)^\omega.\]
\end{lemma}

Also Theorem \ref{one-step} is generalised, and the following now is proved literally as there.
  \begin{theorem}\label{one-stepk}
  After forcing with $\M_k(\cU)$, $(\cU \rest \mu)^+$ is a Matet-adequate family for any finite-to-one $h$, for any $\ba \in (\cU\rest\mu)^+$,
  \[( \exists \bb^1, \bb^2 \sqsubseteq \ba)(\bb^0 \in (\cU\rest\mu)^+ \wedge \bb^1 \in (\cU\rest\mu)^+
  \wedge h[\set(\bb^1)] \cap h[(\set(\bb^2)]=\emptyset).\]
\end{theorem}

  Thus Theorem~\ref{theorem1fork} is proved.
  There is no problem in generalising the iteration theory from Section \ref{S4} and hence we arrive at the following result: 
  
\begin{theorem}\label{3} Assume \CH\ and let $k \geq 1$.
  Then there is a countable support iteration of proper iterands $\bP = \la \bP_\alpha, \M_k(\cU_{\beta}) \such \beta < \omega_2, \alpha \leq \omega_2 \ra$ such that
  in the extension $\cE$ is a $P$-point and
  there is a \kmtu\ over $\Fk$ of character $\gc=\aleph_2=\gd$.
  \end{theorem}

The rest of the section is not used in the main theorem.

Now we are concerned with  the number of near-coherence classes among the ultrafilters that contain $\Phi(\cU)$ as a subset.

The following generalises \cite[Theorem 38]{blass-topap2009} and
adds a new aspect:

  \begin{theorem}\label{topap-generalised2}
    For any union-ultrafilter $\cU$ over $\Fk$,  all the ultrafilters among the projections \[\minl_i(\cU), i = 1,2 \dots, k\]
    are nearly coherent to $\minl_k(\cU)$
    and that class is different from the class  of $ \maxl_k(\cU)$,
    which is the class of any of the ultrafilters among
    \[\maxl_j(\cU), j = 1, \dots, k.
    \]
 \end{theorem}

  \begin{proof} We let $M \subseteq \{1, \dots k-1\}$ be maximal such that
    for any $m \in M$   the set of  $X \in \cU$
    such that $\set_m[X]$ is infinite
  is $\subseteq$-dense in $\cU$. Then exactly for $m \in M\cup\{k\}$, $\maxl_m(\cU)$ and $\minl_m(\cU)$ are ultrafilters. Recall, ultrafilter means non-principal ultrafilter.
  by \cite[Proposition~3.9]{Blass:ufs-hindman}.
  
  In \cite[Section 6]{blass-topap2009} Blass proves that the minimum and the maximum class are different.
  \nothing{
In the first step we consider the case of $\minl_j(\cU)$, $\maxl_j(\cU)$.
As in \cite[Section 6]{blass-topap2009},
we show that any minimum ultrafilter is not nearly coherent to any maximum ultrafilter for arbitrary pairs. 
We assume that $f$ is finite-to-one and $f(\maxl_i(\cU))=f(\minl_j(\cU))$
for some pair  $i \neq j\in \{1,\dots,k\}$.
W.l.o.g. we assume that $f$ is non-decreasing and surjective. We let $I_n = f^{-1}[\{n\}]$.
We let $E$ be 
\begin{equation*}
  E=\{s \in \Fk\such |\{r \in\omega \such I_r \cap s^{-1}[\{i,j\}]\neq\emptyset\}| \mbox{ is even.}\}
\end{equation*}
Since $\cU$ is an ultrafilter over $\Fk$, we  have $E\in \cU$
or $\Fk\setminus E\in \cU$.
Since $\cU$ is a union ultrafilter, we have $E\in \cU$. Hence there is
$\ba\subseteq \Fk$
such that $\FUU_k(\ba)=A\subseteq E$. Since the block-sequences
such that in any block 
for each $i$, $j$ the $i$-minimum lies before the $j$-maximum
and such that no $I_n$ contains both $\min(s^{-1}[\{i\}])$ and
$\max(s^{-1}[\{j\}])$ are $\strk$-dense,
we now assume this for $A$ ($\oplus$)

Since we assumed near coherence we have
\[
f[\{\max(s^{-1}[\{i\}] \such s \in A\}] \cap f[\{\min(s^{-1}[\{j\}] \such s \in A\}] \neq \emptyset.
\]
So there is an interval $I_r$ that meets
$\{\max(s^{-1}[\{i\}]) \such s \in A\}$ and \\
$\{\min(s^{-1}[\{j\}]) \such s \in A\}$, say
$\max(s^{-1}[\{i\}]), \min(t^{-1}[\{j\}]) \in I_r$.
Then by ($\ast$) we have $s<t$.
All of $s^{-1}[\{i,j\}]$ lies before
$\min(t^{-1}([\{j\}])\in I_r$ and hence in intervals $I_{r'}$ with $r' \leq r$.
All of $t^{-1}[\{i,j\}]$ lies after
$\max(s^{-1}([\{i\}])\in I_r$ and hence in intervals $I_{r'}$ with $r' \geq r$.
Hence
$\{r' \such I_{r'} \cap (s\cup t)^{-1}[\{i,j\}] \neq\emptyset\}$ is the 
union of $\{r'\leq r \such I_{r'} \cap s^{-1}[\{i,j\}] \neq\emptyset\}$
and $\{r'\geq r \such I_{r'} \cap t^{-1}[\{i,j\}] \neq\emptyset\}$
and the latter two intersect just in $\{r\}$.
Thus $|\{r \such  I_r \cap (s \cup t)^{-1}[\{i,j\}] \neq \emptyset \}|$ is
odd and contradicts the fact that
$s\cup t\in \FUU_k(B)\subseteq E$.

We show that $\minl_i(\cU)$ is nearly coherent to $\minl_j(\cU)$.
We let
$\ba \in (\Fk)^\omega$ such that $\FUU_k(\ba) \in \cU$.
Then there is some finite-to-one $h$ such that 
\[(\forall n) h(\min(a_n^{-1}[\{1\}]) = \dots = h(\min(a_{n}^{-1}[\{k\}])).
\]
Now any set in $\cU$ is of the form $\FUU_k(\bb)$ for some $\bb $ compatible with $\ba$, we can assume $\bb \sqsubseteq_k\ba$.
Then
\[(\exists^\infty n) h(\min(b_n^{-1}[\{1\}]) = \dots = h(\min(b_{n}^{-1}[\{k\}])).
\]
Thus $\minl_i(\cU)$ and $\minl_j(\cU)$ are nearly coherent non-meagre filters.} 
\end{proof}

\begin{theorem}\label{nc_cores}
  Let $\cU$ be a \kmtu\ over $\Fk$. Let $M$ be as above.
 The cores  $\Phi_j(\cU)$, $j \in M \cup \{k\}$, are all nearly coherent.
  \end{theorem}

\proof
Let $A = \FUU_k(\ba) \in \cU$. We define $I_0 = [0, \max(\supp(a_0)))$,
  \\
  $I_{n+1} = [\max(\supp(a_n)), \max(\supp(a_{n+1})))$.
Then we take a finite-to-one function $h_{\ba}$ that is
constant on $I_n$ for $n \in \omega$.
Then for any $i,j \in M $,
\[(\forall n \in \omega)(h_{\ba}[\set_i(a_n)] \cap h_{\ba}[\set_j(a_n)] \neq \emptyset).
\]
  Since $\cU$ is centred, for
  any $\FUU(\bb), \FUU(\bc) \in \cU$ there is $\bd \sqsubseteq_k \ba,\bb,\bc$,
  $\FUU(\bd) \in \cU$. Hence we have for
  $\bd = \la d_\ell \such \ell < \omega\ra$ for any $\ell\in \omega$ that there are natural numbers
  $n, m ,r \geq 1$ and natural numbers $i_1< \dots < i_n$, $j_1< \dots < j_m$, $k_1< \dots < k_r$,  such that 
  \[d_\ell = a_{i_1} + \dots + a_{i_n} = b_{j_1} + \dots + b_{j_m} = c_{k_1} + \dots + c_{k_r}.\]
  Then 
  \[
  h_{\ba}[\set_i(d_\ell)]
  \cap h_{\ba}[\set_j(d_\ell)] \mbox{ has size at least $\max(n,m,r)$.}
  \]
  Hence  \[
  h_{\ba}[\bigcup\{\set_i(d_\ell) \such \ell \in \omega\}]
  \cap h_{\ba}[\bigcup\{\set_j(d_\ell)\such \ell \in \omega\}] \mbox{ is infinite.}
  \]
Thus also the superset
  \[
  h_{\ba}[\bigcup\{ \set_i(b_\ell) \such \ell \in \omega\}]
  \cap h_{\ba}[\bigcup\{\set_j(c_\ell)\such \ell \in \omega\}] \mbox{ is infinite.}
  \]
  \proofend

  \begin{remark} Theorem~\ref{nc_cores} can also be proved with the
    help of Theorem~\ref{eisworthfork} and the following folklore
    result.
    \end{remark}
  
    \begin{proposition}\label{folklore}
      Any forcing that diagonalises two ncc filters adds a dominating real.
    \end{proposition}

\begin{proof} Let $\cF$ and $\cF' \in \bV$ be nnc filters and let
    $x, y \in \bV[G]$ be such that $(\forall F \in \cF)( x \subseteq^* F)$,
    $(\forall F' \in \cF') (y \subseteq^* F')$.
    Then by the next lemma,
    \[n \mapsto \max\bigl\{\min(x \setminus (n+1)), \min(y \setminus
    (n+1))\bigr\}\]
    is a dominating function.

\begin{lemma}\label{factaboutdom}(Proof of \cite[Theorem 3.2]{ncf2})
Let $\cV$, $\cW$ be non-principal filters over $\omega$.
$\cV$ is nearly coherent to $\cW$ iff
\[ \{\max(\nex(\cdot,X),\nex(\cdot,Y)) \such X \in \cV, Y\in \cW\}\]
is not a $\leq^*$-dominating family.
\end{lemma}

We give a proof for the direction we use:
\nothing{
Assume that $\cV$ and $\cW$ are nearly coherent. Let $f$ be a finite-to-one function such that  $f(\cV) \cup f(\cW)$ is a filter base.
W.l.o.g.\ we assume that there is a strictly increasing sequence of natural numbers $\la \pi_i \such i < \omega\ra$ such that $\pi_0 = 0$ and $f(n) =
i$ for $n \in [\pi_i,\pi_{i+1})$. We let $h(i) = \pi_{n+1}$ for $i \in [\pi_n,\pi_{n+1})$. We show that
 no member of $F$ dominates $h$. 
 Let $X \in \cV$, $Y \in \cV$. Then $f[X] \cap f[Y]$ is infinite.
 Let $i = f(x) = f(y)$ for some $x \in X$, $y \in Y$. Then $x,y \in [\pi_i,\pi_{i+1})$, and we assume that $x,y$ are minimal in $X \cap [\pi_i,\pi_{i+1})$
     and in $Y \cap [\pi_i,\pi_{i+1})$ respectively. Thus $h(\pi_i) = \pi_{i+1} >
       \max(x,y) = \max(\nex(\pi_i,X),\nex(\pi_i, Y))$.
       Hence we have $h \not\leq^* \max(\nex(\cdot,X),\nex(\cdot,Y))$.
}       
Assume that $\cV$ and $\cW$ are not nearly coherent. Let $h \in {}^\omega \omega$ be given, w.l.o.g.\ we assume that $h$ is strictly increasing and $h(0)>0$.
       We let $\tilde{h}$ be the iterate of $h$: $\tilde{h}(0) = 0$,
       $\tilde{h}(n+1) = h(\tilde{h}(n))$.
       We let $f_e(n) = i$  for $n \in [\tilde{h}(2i), \tilde{h}(2i+2))$ and we let
         $\tilde{h}(-1) =0$ and
         $f_o(n)= i $ for  $n \in [\tilde{h}(2i-1), \tilde{h}(2i+1))$.
           Then there are $V_e, V_o \in \cV$ and $W_e, W_0 \in \cW$ such that
           $f_e[V_e] \cap f_e[W_e] = \emptyset$ and
           $f_o[V_o] \cap f_o[W_o] = \emptyset$. Since $\cV$ and $\cW$ are filters, we can assume $V_e = V_o$ and $W_e = W_o$.
           Then for any $n  \in \omega$ we
           have $\max(\nex(n,V_e),\nex(n,W_e))\geq h(n)$.
\end{proof}
\nothing{other proof sketch           
    An interval partition $\bar{I}$ dominates an
    interval partition $\bar{J}$ if $\forall^\infty n \exists k J_k \subseteq I_n$.  Let $\bar{I}$ be an interval partition in $\bV[G]$ such that $\forall n (I_n \cap x \neq \emptyset \wedge I_n \cap y \neq \emptyset)$. Suppose for a contradiction, that there is an interval partition $\bar{J}$ in $\bV$ that is not dominated by $\bar{I}$.
    We let for $j = 0,1$, $\ell \in \omega$, $h_j(\ell) = n$ if $\ell \in J_{2n+j} \cup J_{2n+1+j}$.
    Since $\bar{J}$ is not dominated by $\bar{I}$, there is at least one $j$ such that $h_j(\cF) \cup h_j(\cF')$ is a filter subbase, and this contradicts the assumption on $\cF$ and $\cF'$. Hence $\bar{I}$ is a dominating interval partition. According to Blass' \cite[Proposition 2.10]{Blasshandbook}, $\bar{I}$ can be transformed into $\leq^*$-dominating real.
\proofend
}

Now we finish the alternative proof of Theorem~\ref{nc_cores}: $\M_k(\cU)$ preserves $\cE$ according to Theorem~\ref{eisworthfork} and hence does not add a dominating real. However $\M_k(\cU)$ diagonalises $\Phi_j(\cU)$ for $j = 1, \dots, k$ by Remark~\ref{diagonalisers}. By Proposition~\ref{folklore} all the $\Phi_j(\cU)$, $j =1, \dots, k$, are nearly coherent.

\section{A Ramsey subspace of the space of $\Fk$-sequences}
  \label{S6}

  Now we change from the Hindman $\Fk$ space to the Gowers $\Fk$ space, i.e.,
  we incorporate the  Tetris operation and use the resulting generalised
  type of condensation as partial order.
  Then we localise the
  space $((\Fk)^\omega, \strk)$  by choosing $2k$ pairwise nnc Ramsey ultrafilters $\cR_{i,{\min}}$, $\cR_{i,{\max}}$, $i = 1, \dots, k$,
  that will serve as representatives of names of near-coherence classes and get
  a subspace $((\Fk)^\omega(\bar{\cR}), \strk)$, see Def.~\ref{PP}.
  We show that this space has good Ramsey theoretic properties, as
in Theorem~\ref{gowersown}, and we introduce \gmtu s.

We tell more about our plans: 
In Section \ref{S7_neu} we define a variant of Matet forcing with a non-centred reservoir
  that is given by the sequences $\ba \in (\Fk)^\omega(\bar{\cR})$.
  The resulting forcings are in contrast to the
  partial order $\M_k(\cU)$ with a \kmtu\  $\cU$ (or to $\GM_k(\cU)$ with a \gmtu\ $\cU$) not
  $\sigma$-centred. Nevertheless they still fulfil Axiom A and hence are proper.

In order to prove that the new forcings are proper, we insert this section whose main result,  Theorem \ref{gowersown1},
  is a common strengthening of Gowers' theorem and of Blass' theorem \cite[Theorem 2.2]{Blass:ufs-hindman}.
  
  \begin{definition}\label{collection_on_Tetris}
 Let $k \in \omega \setminus \{0\}$ unless stated otherwise. 
  \begin{myrules}
\item[(1)]
For any $j \geq 2$ we define on $\Ftwo{j}$ the \emph{Tetris operation}:
$T_j \colon \Ftwo{j} \to \Ftwo{j-1}$ by $T_j(p)(n) = \max\{p(n)-1, 0\}$ and let
$T= \bigcup \{T_j \such 2 \leq j \leq k\}$.

For $\ba \in (\Ftwo{j})^\omega$, we write $T[\ba]$ for $\la T(a_n) \such n < \omega \ra$. For $X \subseteq \Ftwo{j}$, we write $T[X]$ for $\{T(x)
\such x \in X\}$.
\item[(2)]
  Let $B$ be a min-unbounded subset of $\Fk$. We let
\begin{equation*}
\begin{split}
    \TFU_k(B) = & \{ T^{(j_0)}(b_{n_0}) +\dots  + T^{(j_\ell)}(b_{n_\ell}) \such\\ 
    & \ell \in \omega, b_{n_i} \in B,\\
&
    b_{n_0} < \dots < b_{n_\ell}, j_i \in k, (\exists i \leq \ell)( j_i = 0)\}
\end{split}
\end{equation*}
be the partial subsemigroup of $\Fk$ generated by $B$.
Here $T^{(0)} $ is the identity, $T^{(1)} = T$, $T^{(j+1)} = T \circ T^{(j)}$.

\item[(3)] For $\bb \in (\Fk)^\omega$ we let 
$\TFU_k(\bb) = \TFU_k(\{b_n \such n \in \omega\})$.
We call $B$ a \emph{$\TFU_k$-set} if $B = \TFU_k(\bb)$ for some
$\bb \in (\Fk)^\omega$.

 \item[(4)]
For $j \geq 2$,  we define $\Tl \colon \gamma(\Ftwo{j}) \to \gamma(\Ftwo{j-1})$ via
   \[\Tl(\cU) 
  = \bigl\{ X \subseteq \Ftwo{j-1} \such \{ s\such T(s) \in X\} \in \cU\bigr\}. 
  \]
\item[(5)] For min-unbounded sets $A, B \subseteq \Fk$ we let
  $B \strk A$ if $B \subseteq \TFU_k(A)$. Note that here $k$ is not allowed to drop. We say $B$ is a $k$-Tetris-condensation of $A$.
  We write
  $B \strkast A$ if there is an $n \in \omega$  such that
  $(B \past n) \subseteq \TFU_k(A)$.
\item[(6)] We introduce the Tetris condensation order $\strk$ on sequences.
  For $\ba, \bb \in (\Fk)^\omega$ we let
 $\bb \strk \ba$ if
  $\bb = \la b_n \such n < \omega \ra$ and $\{b_n \such n < \omega \} \subseteq \TFU_k(\{a_n \such n < \omega \})$, and $\bb \strkast \ba$ has the obvious meaning.
  \end{myrules}
  \end{definition}

  From the following lemma we will use only the part of the statement that
  the described procedure produces a lower bound. 
  
\begin{lemma}\label{luxury3}
  Assume that $\ba, \bb \in (\Fk)^\omega $ and that there is
  a $\bc \in (\Fk)^\omega$ such that $\bc \strk \ba, \bb$.
    Then there is a $\strk$-largest witness $\bc$ to this, that can be gotten
  in the following recursive way:
  \begin{eqnarray*}
    c_0 &= & \min_{<_{{\rm lex},\Fk}}( \TFU_k(\ba) \cap \TFU_k(\bb)),
    \\
    c_{n+1} &= & \min_{<_{{\rm lex},\Fk}} (\TFU_k(\ba \past c_n) \cap \TFU_k(\bb \past c_n)).
    \end{eqnarray*}
\end{lemma}

The proof that $\bc$ is the largest lower bound requires additional analysis in the second subcase of the third case in the proof of Proposition~\ref{luxury1}, due to the Tetris operation.
We skip it.

An important and deep property of the space $((\Fk)^\omega, \strk)$ is the following:

\begin{theorem}\label{gowers}(Gowers, see e.g. \cite{Gowers92}, \cite[Theorem 2.22]{Todorcevic:Ramsey}) Let $k \geq 1$ and  $\ba \in (\Fk)^\omega$. For every finite colouring of $\TFU_k(\ba)$ there is a infinite block sequence $\bb \strk \ba$ such that the elements of $\TFU_k(\bb)$ are monochromatic.
    \end{theorem}

The following proposition gives a first picture:

\begin{proposition}\label{ordered_stairs_dense}
  For each $\ba \in (\Fk)^\omega$ there is a $\bb \strk \bb$
  such that for each $s \in \TFU_k(\bb)$ we have
  \begin{multline}\label{ordered}
      \min_1(s) < \min_2(s) < \dots <  \min_k(s) < 
        \max_k(s)
        <  \max_{k-1}(s) < \dots < \max_1(s).
  \end{multline}
  \end{proposition}

\proof For $s \in \TFU_k(\ba)$ we let
$c(s) = 1$ if $s$ fulfils \eqref{ordered}
and $c(s) = 0$ otherwise.  By Gowers' Theorem there is $\bb \strk \ba$, such that $\TFU_k(\bb)$ is monochromatic for $c$. We show that the colour can only be 1.
We take $s_i \in \TFU_k(\bb)$, $s_1< s_2 < \dots < s_{2k-1}$. Then 
\begin{equation}
  \begin{split}
    s = & T^{(k - 1)}(s_1) +  T^{(k-2)}(s_2) + \dots + T (s_{k-1}) + s_k \\
        & + T (s_{k+1}) + \dots + T^{(k-1)} (s_{2k-1}) \in  \TFU_k(\bb)
  \end{split}
\end{equation}
fulfils $c(s) = 1$.
  \proofend

\begin{definition}\label{Tetris_ksetting2}
  \hfill
  \begin{myrules}
  \item[(1)]
  A non-principal filter $\cF$ over $\Fk$ is said to be an {\em Tetris-ordered-union} filter if it has a basis of sets of the form $\TFU_k(\bar{d})$ for $\bar{d}\in (\Fk)^\omega$.
  \item[(2)]
    For a  $\strk$-descending $\omega$-sequence of members of $(\Fk)^\omega$
    $\la \ba_n \such n \in \omega \ra$ and $\bb \in (\Fk)^\omega$ we say: $\bb$ is a \emph{diagonal lower bound} if $\bb \strk \ba_0$ and
    \[ (\forall s \in \TFU_k(\bb)) (\bb \past s) \strk \ba_{\max(\supp(s)) +1}.\]

\item[(3)]
Let $\mu$ be an uncountable cardinal. A Tetris-ordered-union filter 
is said to be {\em $(<\mu)$-stable} if,
 whenever it contains $\TFU_k(\bb_\alpha)$ for $\bb_\alpha \in (\Fk)^\omega$, 
$\alpha<\kappa$, for some $\kappa<\mu$, then it also contains 
some $\TFU_k(\bb)$ for some $\bb\in (\Fk)^\omega$ such that for each $\alpha$
there is $n_\alpha$ with $(\bb \past \{n_\alpha\}) \strk \bb_\alpha$
for $\alpha<\kappa$. Such an $\bb$ is called a \emph{lower bound} of
$\{\bb_\alpha \such \alpha<\kappa\}$. For ``$<\omega_1$-stable'' we say ``stable''.
\item[(4)]
A stable Tetris-ordered-union ultrafilter is  called a {\em \gmtu}. 
\end{myrules}
\end{definition}

\begin{definition}\label{Gowersadequatefamily}
A set $\cH \subseteq (\Fk)^\omega$ is called a \emph{Gowers-adequate family} if
the following hold:
\begin{myrules}
\item[(i)] $\cH$ is closed $\strkast$-upwards.
\item[(ii)] $\cH$ is stable, i.e.,
any $\strk$-descending $\omega$-sequence of members of $\cH$ has a $\strkast$-lower bound in $\cH$.
\item[(iii)] $\cH$ has the \emph{Gowers property}: If $\ba \in \cH$ and $\TFU_k(\ba)$ is partitioned into finitely many pieces then there is some $\bb \strk \ba$, $\bb \in \cH$, such that $\TFU_k(\bb)$ is a subset of a single piece of the partition.
  \end{myrules} 
\end{definition}
In particular,
any \gmtu\ is Gowers-adequate.
Gowers-adequate families have better properties than stated in the definition;
the analogue to Lemma~\ref{selfstr} hold for them.

\begin{definition}\label{schreibweisen}
  In the general case of $i  \in \{1, \dots, k\}$, $x \in \{{\min},{\max}\}$, we let the variable $x$ range over $\{{\min}, {\max}\}$, and we also write $x_i(a)$ for $\min_i(a)$, $\max_i(a)$.
We call $i \in \{1, \dots, k\}$ sometimes a coordinate,  a value
or  a level.
\end{definition}

We introduce a new space, called
$((\Fk)^\omega(\bar{\cR}), \strk)$, the one advertised in the Equation~\eqref{gist} in the introduction.

  \begin{definition}\label{PP}
Let $k\geq 1$ and let
    \[\bar{\cR} = \la \cR_{i,x} \such i = 1, \dots, k, x = {\min},{\max}\ra
    \]
    be a $\{1, \dots ,k\} \times \{{\min},{\max}\}$-sequence of pairwise nnc Ramsey ultrafilters.
  We define 
\begin{equation}\label{wit}
  \begin{split}
    (\Fk)^\omega(\bar{\cR}) = & \{\ba \in (\Fk)^\omega \such \forall (X_{i} \in \cR_{i,{\min}}, Y_i \in \cR_{i,{\max}})_{1 \leq i \leq k}\\
    &(\exists^{\rm min-unb} s \in \TFU_k(\ba))
    \bigwedge_{1 \leq i \leq k}(\min_i(s) \in X_i \wedge \max_i(s) \in Y_i)\}.
  \end{split}\end{equation}
  \end{definition}
  Here we write  $(\exists^{\rm min-unb} s)$ as a abbreviation
  for ``the set of witnesses $s$ in min-unbounded''.
\begin{remark}
  We remark that Equation\eqref{wit}
is equivalent to
\begin{equation}
  \begin{split}
    (\Fk)^\omega(\bar{\cR}) = & \{\ba \in (\Fk)^\omega \such \forall (X_{i} \in \cR_{i,{\min}}, Y_i \in \cR_{i,{\max}})_{1 \leq i \leq k}\\
    &(\exists \bb \strk \ba)(\forall n)
    \bigwedge_{1 \leq i \leq k}(\min_i(b_n) \in X_i \wedge \max_i(b_n) \in Y_i)\}.
  \end{split}
\end{equation}

However, it speaks about the $b_n$ and \emph{not} about the elements of 
$\TFU_k(\bb)$.
\end{remark}

First we recall \cite[Lemma 1.1]{Blass:ufs-hindman}:

\begin{lemma}\label{Blass1.1} Let $\cR_{1,{\min}}$, $\cR_{1,{\max}}$ be nnc Ramsey ultrafilters,
$X_1 \in \cR_{1,{\min}}$, $Y_1 \in \cR_{1,{\max}}$ and let
$(B_n)_n$ be a descending sequence of elements of $\cR_{1,{\min}}$, $(D_n)_n$ be a descending sequence of elements of $\cR_{1,{\max}}$.
Then there are $X_1^1 \subseteq X_1$, $Y^1_1 \subseteq Y_1$, $X_1^1 \in \cR_{1,{\min}}$, $Y_1^1 \in \cR_{1,{\max}}$ such that
\begin{myrules}
\item[(1)] $(\forall x \in X_1^1)( \forall y \in Y_1^1) (y<x \rightarrow x \in B_y)$.
\item[(2)] $(\forall y \in Y_1^1 )(\forall x \in X_1^1) (x<y \rightarrow y \in D_x)$.
\item[(3)] In the increasing enumeration of $X_1^1 \cup Y_1^1$ the elements of $X_1^1$ and $Y_1^1$ alternate.
\end{myrules}
\end{lemma}

\begin{lemma}\label{not_centred} There are $\sqsubseteq_k^*$-incompatible elements in
  $(\Fk)^\omega(\bar{\cR})$. 
\end{lemma}

\begin{proof}
  For notational simplicity we assume that $k=1$ and
  the Ramsey ultrafilters are $\cR_{\min}$ $\cR_{\max}$. For transferring the situation to the general case, just replace the value 1 by $k$ and work
  only with the top ultrafilters $\cR_{k,\min}$, $\cR_{k,\max}$.
  Lemma~\ref{Blass1.1}, applied to $X_1= Y_1 = \omega$,
  $A_n = B_n = [2n,\omega)$ yields sets
  $X \in \cR_{{\min}}$ and $Y \in \cR_{{\max}}$ such that
in the  increasing enumeration
 $\la x_n \such n < \omega \ra$ of $X \cup Y$ the even numbers just give $X$ 
 and the odd number just give $Y$ and that the $x_n$ have distance at least $2$.
Then we take a 
natural number, call $\ell_n$, such that $x_{2n} < \ell_n < x_{2n+1}$.
  Then \begin{equation*}\begin{split}
      \ba = &\la \{x_{2n},x_{2n+1}\} \such n \in \omega \ra, \\
      \bb = & \la \{x_{2n},\ell_{n}, x_{2n+1}\} \such n \in \omega \ra
  \end{split}\end{equation*}
  have
  \begin{equation} \label{allowed}
      X = \{\min_1(a_n) \such n \in \omega\} \wedge Y = \{\max_1(a_n)
    \such n \in \omega \},
  \end{equation}
  and hence $\ba \in (\Ftwo{1})(\cR_{\min}, \cR_{\max})$,
  and the same holds for $\bb$.
Obviously, the conditions $\ba$, $\bb \in (\Ftwo{1})^\omega(\bar{\cR})$ are
  $\strtwo{1}$-incompatible.
\end{proof}

By the way, our only possibility to verify that some sequence is
in $(\Fk)^\omega(\bar{\cR})$ is to construct one that has ultrafilter
projections as in Equation~\eqref{allowed}.

By applying e.g.\ the trick in \cite[Cor.~2.3]{Blass:ufs-hindman} one sees that also
below every condition there are incompatible conditions.

\begin{lemma}\label{for_Thomas_the_doubter}
  For any ultrafilter $\cU$ over $\Fk$, the projections
  $ \minl_k(\cU)$, $\maxl_k(\cU)$ generate (by adding supersets) ultrafilters,
  and for $1 \leq i < k$, either
  each of the projections
  $ \minl_i(\cU)$, $\maxl_i(\cU)$ generates
a filter or each of them  contains the empty set.
  \end{lemma}

\begin{proof} They are filters:
  $\minl_i(\cU) = \{\min_i[X] \such X \in \cU\}$. Since $\cU$ is a filter,
  the projection $\minl_i(\cU)$ is closed under finite intersections.
  For $i < k$ it is possible that there is an $X \in \cU$ such that
  $\min_i[X] = \max_i[X]= \emptyset$.
  By definition of $\Fk$, the top coordinate $\minl_k(\cU)$
  does not contain  the empty set. The family
 $\minl_k(\cU)$  generates an ultrafilter: Let $A \subseteq \omega$. Then
  $\Fk= \{a \in \Fk \such \min_k(a) \in A\} \cup \{a \in \Fk \such \min_k(a) \not\in A\}$ is a partition of $\Fk$, and one part lies in $\cU$.
  The same holds for the maximum projection.
\end{proof}

\begin{definition}\label{projections}
  For an ultrafilter $\cU \in \gamma(\Fk)$ we say \\
  \emph{$\cU$ projects to
  $\bar{\cR}$}, if
  \[(\forall i \in \{1,\dots, k\}) (\minl_i(\cU) = \cR_{i,\min} \wedge
         \maxl_i(\cU) = \cR_{i,\max}).
         \]
\end{definition}

\begin{definition}
  Let $\ba \in (\Fk)^\omega$. $\set(\ba) = \bigcup\{\supp(a_n) \such n < \omega\}$.
\end{definition}

In the following Theorem~\ref{gowersown}, the case of $k=1$ is due to Blass' \cite[Theorem 2.2]{Blass:ufs-hindman}, and the unlocalised case of $k \geq 1$ is proved in Gowers' theorem \cite[2.22]{Todorcevic:Ramsey}. In a sense the result can be seen as an instance of the programme outlined
by Todorcevic in \cite[Remark 7.27]{Todorcevic:Ramsey}.
The author is thankful to Andreas Blass for his careful explanation of his proof of \cite[Theorem 2.1]{Blass:ufs-hindman}.

\begin{theorem}\label{gowersown}
Assume \CH.  For any $k$, $\bar{\cR}$
  as in Def.~\ref{PP} the following hold:
\begin{myrules}
\item[(1)]
 Let $\ba \in (\Fk)^\omega(\bar{\cR})$ and let $c $ be a colouring of
    $\TFU_k(\ba)$ into finitely many colours. Then there is a $\bb \strk \ba$,
 such that $\bb \in (\Fk)^\omega(\bar{\cR})$ and $\TFU_k(\bb)$ is $c$-monochromatic.

\item[(2)] Let $\ba_n \in (\Fk)^\omega(\bar{\cR})$, $n \in \omega$, be
  $\strk$-decending. Then there is a $\strkast$-lower bound in 
$(\Fk)^\omega(\bar{\cR})$. 

\item[(3)]
If  $\ba \in (\Fk)^\omega(\bar{\cR})$ and $\cR_{0,{\min}}$, $\cR_{0,{\max}}$ are two additional nnc ultrafilters, that are nnc
to each of the $\cR_{i,x}$, $i = 1, \dots, k$, $x = {\min}, {\max}$,
then there is a \gmtu\ $\cU_k$ over $\Fk$ such that
$\TFU_k(\ba) \in \cU_k$ and for $i =1, \dots , k$,
$\minl_i(\cU) = \cR_{i,{\min}}$, $\maxl_i(\cU) = \cR_{i,{\max}}$ and
$\Phi(\cU_k)$  is nnc 
$\cR_{0,x}$, $x = {\min}, {\max}$.
\end{myrules}
\end{theorem}

Since parts (1) and (2) are absolute, we get:

\begin{theorem}\label{gowersown1}\hfill
Let $k \geq 1$ and let $\bar{\cR}=\la \cR_{i,\min}, \cR_{i,\max}
  \such i = 1,\dots k\ra$ be a sequence  of pairwise non nearly coherent Ramsey ultrafilters.
  \begin{myrules}
\item[(1)]
  Any $\strk$-descending $\omega$-sequence of elements of $(\Fk)^\omega(\bar{\cR})$
  has a $\strkast$-lower bound in $(\Fk)^\omega(\bar{\cR})$.
\item[(2)]
Let $n \in \omega \setminus \{0\}$ and $\ba \in (\Fk)^\omega(\bar{\cR})$ and let $c $ be a colouring of
  $[\TFU_k(\ba)]^n_{<}$ into finitely many colours. Then there is a $\bb \sqsubseteq_k \ba$, $\bb \in (\Fk)^\omega(\bar{\cR})$ such that $[\TFU_k(\bb)]^n_{<}$ is $c$-monochromatic.
\end{myrules}
\end{theorem}

We prove Theorem \ref{gowersown}.

\proof
We prove (1), (2) and (3) simultaneously by induction on $k$.
\nothing{(3) at $k-1$ is essential for (1) at $k$. (1) and (2) at $k$  are used for
(3) at $k$.
Also for (2) at $k$ we use only (1), (2) and (3) at $k-1$.
So there is no decided order between (1) and (2) at $k$.}
 The case of $k = 1$ is given by Blass and will be repeated below for
$k = 1$.

It suffices to prove (1) of the theorem for a partition 
of  $\Fk$ instead of $\TFU_k(\ba)$
since $T^{j_0}(a_{n_0}) + \dots + T^{j_r}(a_{n_r})\in \TFU_k(\ba)$ corresponds
to $(n_0, k-j_0), \dots ,(n_r, k-j_r) \in \Fk$.

We collect some technical facts separately, before we prove
Theorem~\ref{gowersown}.

First we give a general pattern for proving the existence of a particular ultrafilter over $\Ftwo{i}$ via compactness.

  \begin{lemma}\label{thanks_to_Andreas}
   Let $\cR_{\min}$, $\cR_{\max}$ be two nnc Ramsey ultrafilters.
   Let \\
   $\ba \in (\Ftwo{1})^\omega(\cR_{\min}, \cR_{\max})$. The compact Hausdorff space
  \[
  \{\cU \in \gamma(\Ftwo{1}) \such \minl(\cU) = \cR_{\min}
  \wedge \maxl(\cU) = \cR_{\max} \wedge \TFU_1(\ba)\in \cU\}
  \]
  is not empty.
\end{lemma}

\begin{proof}
  Since $\minl$ and $\maxl$ are continuous function, the space is compact.
  We show that there is an ultrafilter $\cU$ over
  $\Ftwo{1}$ in the named space.
  For this it suffices to show
  that
  for any $X \in \cR_{\min}$ and any $Y \in \cR_{\max}$
  \begin{equation*}\tag*{R(X,Y)}
    \begin{split} 
  &\mbox{There are min-unboundedly many $s \in \TFU_1(\ba)$}\\
    &\mbox{such that
  $\min_1(s)\in X$ and $\max_1(s) \in Y$.}
\end{split}
\end{equation*}  

We argue why this suffices:
Then the open and closed set
\[C(X,Y) = \{\cU \in \gamma(\Ftwo{1}) \such
(\exists A \in \cU)  \mbox{($A$ is a min-unb.\ set of witnesses to }R(X,Y))\}\]
is not empty.

  Since $\cR_{i,x}$ are ultrafilters, the intersection of 
  finitely many sets of the form $C(X,Y)$ is again a set of the
  $C(X',Y')$ with $X'$ being the intersection of the finitely many $X$'s and $Y'$ being the intersection of the finitely many $Y$'s.

  Thus any finite intersection of sets of the form $C(X,Y)$ is not empty, and by compactness of $\gamma(\Ftwo{i})$, there is an ultrafilter
  \[\cU \in \bigcap\{ C(X,Y) \such X \in \cR_{k,\min},  Y\in \cR_{k,\max}\} \subseteq \gamma(\Ftwo{i}),\]
  which means 
  \[\minl_i(\cU) = \cR_{k,{\min}} \wedge \maxl_1(\cU) = \cR_{k,{\max}}.\]
  
  To fulfil requirement $R(X,Y)$:
  We verify that for any $n \in \omega$
  there is a block\\ $s\in (\TFU_1(\ba) \past \{n\})$ such that
  $\min_1(s) \in X$ and $\max_1(s) \in Y$. Why?
  Since $\ba \in (\Fk)^\omega(\bar{\cR})$,
  and $X \in \cR_{\min}$ and $Y \in \cR_{\max}$ there
  are min-unboundedly many  $s \in \TFU_1(\ba)$ such that
$\min_1(s) \in X
 \wedge \max_1(s) \in Y$.
   \end{proof}

Starting from now, we just write the requirements and show how to fulfil them.

For proving statement (1) at $k$ on the bases of (1) and (3) at $k-1$,  
we rework and extend the proof of Gowers' original theorem as given by Todorcevic \cite[pp.~35--36]{Todorcevic:Ramsey}.

\begin{lemma} \label{fitting} Let $\cR_{i,x}$, $i = 1, \dots, k$, $x = {\min},{\max}$
be pairwise nnc Ramsey ultrafilters.
For any $\ba \in (\Fk)^\omega$ the following are equivalent:
\begin{myrules}
\item[(1)] $\ba \in (\Fk)^\omega(\bar{\cR})$.
\item[(2)] There is an ultrafilter $\cU \in \gamma(\Fk)$ such that 
\[\TFU_k(\ba) \in \cU \wedge 
(\forall i \in\{ 1, \dots k\}) 
(\minl_i(\cU)= \cR_{i,{\min}} \wedge \maxl_i(\cU) = \cR_{i,\max}).
\]
\end{myrules}
\end{lemma}

\begin{proof} (1) implies (2): By the compactness method as
explained in the previous lemma it suffices to show that for any
$(X_i \in \cR_{i,{\min}} Y_i \in \cR_{i,{\max}} )$
\begin{equation}\label{witnesses}
(\exists^{\rm min-unb} s \in \TFU_k(\ba))(\forall i \in \{1, \dots, k\}) (\min_i(s) \in X_i \wedge \max_i(s) \in Y_i).
\end{equation}
However, this is guaranteed by $\ba \in (\Fk)^\omega(\cR)$.

(2) implies (1). Let $\cU$ be as in (2). Given sets $X_i \in \cR_{i,{\min}}$,  $Y_i \in \cR_{i,{\max}}$, $i = 1, \dots, k$, there is a set $S \in \cU$ such that 
\[
S \subseteq \TFU_k(\ba) \wedge (\forall i \in \{1, \dots, k\})
(\min_i[S] \subseteq X_i \wedge \max_i[S] \subseteq Y_i).
\]
Since $S \in \cU$, $S$ is min-unbounded. Any $s \in S$ fulfils
$\min_i(s) \in X_i \wedge \max_i(s) \in Y_i$.
\end{proof}

We name the set of ultrafilters with given projections: 

\begin{definition}\label{gammaspace_new}
For $1 \leq i  \leq k$ and a sequence $\la \cR_{j,x} \such 1 \leq j \leq i, x \in \{{\min},{\max}\}\ra$ or pairwise nnc Ramsey ultrafilters we let
\begin{equation}\label{positive_crucial}
  \begin{split}
   &\gamma^{\rm good} (\Ftwo{i}, \bar{\cR}) = 
    \{\cU \in \gamma(\Ftwo{i})  \such
    (\forall j \in \{1,\dots, i\})\\
& (\minl_j(\cU) = \cR_{j,{\min}} \wedge \maxl_j(\cU) = \cR_{j,{\max}})\}.
      \end{split}
\end{equation}
\end{definition}

Cave: Sometimes we use ultrafilter sequences that carry shifted indices.

\begin{lemma}\label{dot{+}}
For any $i$, $\bar{\cR}$ as above,  the space $\gamma^{\rm good}(\Ftwo{i},\bar{\cR})$ is closed under $\dotpl$.
 \end{lemma}

\begin{proof}
This is proved in \cite[page 99, second but last paragraph]{Blass:ufs-hindman}.
\end{proof}

In Lemma~\ref{closure_under_dot+} we shall look at ``mixed additions.''
First we state some simple observations.

\begin{lemma}\label{simple}\hfill
  \begin{myrules}
  \item[(1)] For any $X \subseteq \Ftwo{i}$, for any $1\leq j < j' \leq i$
    we have \[\min_j[X] \cap \min_{j'}[X] = \emptyset,\]
    since for each $s \in \Ftwo{i}$, $\min_j(s) \neq \min_{j'}(s)$.
  \item[(2)] Given $t \in \Ftwo{j}$, there is a particular (in $<_{\Ftwo{j+1},lex}$ largest)
    $s \in \Ftwo{j+1}$ with $T(s) = t$:
    We let $\supp(t) = \supp(s)$ and $t(n) = s(n) +1$ for $n \in \supp(s)$. We write for this $t$ now $s^\sharp$. 
  \item[(3)]   We can lift the sharp operation by letting
    $X^\sharp := \{s^\sharp \such s \in X\}$. 
  Given $\cU \in \gamma(\Ftwo{j})$, there is a particular ultrafilter
  $\cV \in \gamma(\Ftwo{j+1})$ with $\Tl(\cV) = \cU$,
  namely \[\cV= \{X^\sharp \such X \in \cU\}.\]
  We write $\cV = \lift(\cU)$.
  Then $\minl_{j'+1}(\cV) = \minl_{j'}(\cU)$ for $j' = 1, \dots, j$, and
  $\minl_{1}(\cV) = \{\emptyset\}$ and analogously for $\maxl$.
\nothing{
\item[(4)] Only on  level  1, the Tetris operation is not
  invertible. This means that that for $i > 1$,
  $\{s^\sharp \such s \in \TFU_{i-1}(T^{(k-i+1)}[\ba]\}$ looks on its levels $2, \dots, i$ like 
  of $\TFU_i(T^{(k-i)}[\ba])$ and for each $s \in \TFU_{i-1}(T^{(k-i+1)}[\ba])$ there
  is at least one modification $t$ of $s^\sharp$ that adds finitely many $n$
  from $\omega \setminus \supp(s^\sharp)$ into $\supp(t)$ and sets
  $t(n) = 1$ such that $t \in \TFU_i(T^{(k-i)}[\ba])$

\item[(5)] $\TFU_i(T^{(k-i)}[\ba]) \in \cU \in \gamma(\Ftwo{i})$ and
  and $\TFU_{i-1}(T^{(k-i+1)}[\ba]) \in \cV \in \gamma(\Ftwo{i-1})$
  and $\Tl({\cU}) = \cV$ implies that $\TFU_i(T^{(k-i)}[\ba]) \in \cU \dotpl \cV$.
  }
  \end{myrules}
\end{lemma}

Now we prove a technical lemma about the sums of good ultrafilters.
Prima facie the ultrafilters do not fit, e.g., $\min_{k}(\cX) = \cR_{k,{\min}}$
and $\min_{k-1}(\cV) = \cR_{k,{\min}}$. However, this mixed sum 
is defined on purpose.

\begin{lemma}\label{closure_under_dot+}
  Let $\bar{\cR} = \la \cR_{i,x} \such 1 \leq i \leq k, x = {\min}, {\max}\ra$.
  Let $\cX \in \gamma^{\rm good}(\Ftwo{k},\bar{\cR})$ and
  $\cV \in \gamma^{\rm good}(\Ftwo{k-1},\la \cR_{2,x}, \dots, \cR_{k,x}, x = {\min},{\max}\ra)$.
  Then \\
  $\cW=\cX \dotpl \cV \in   \gamma^{\rm good}(\Ftwo{k},\bar{\cR})$.
  \end{lemma}

\begin{proof}
We have $\cW \in \gamma(\Ftwo{k},\bar{\cR})$
iff for $v = 1, \dots, k$, $x = {\min}, {\max}$, for any $X \in \cW$, 
\begin{equation}\label{rightprojection}
  \{x_v(s) \such s \in X \}\in \cR_{v,x}.
\end{equation}
We fix some $X \in \cW$. Then $\{s \such \{t \such s+t \in X\} \in \cV\} \in \cX$.
For any $s \in \Ftwo{k}$, $t \in \Ftwo{k-1}$,
such that $s < t$ and $s +t \in X$, and for any
 $v = 1, \dots, k$, $x = \min,\max$,
$x_v(s+t)$ is either $x_v(s)$ or $x_v(t)$.
Either the set of $(s,t)$ on which $x_v$ is the first  is in
\[\cX \times \cV= \{Z \subseteq \Ftwo{k} \times \Ftwo{k-1} \such
\{s \such \{t \such (s,t) \in Z\} \in \cV\} \in \cX\},
\]
or its complement is in $\cX \times \cV$.

We start with the high values of $v$:
For the values $v = k$, we have for any $s \in \Ftwo{k}$, $t \in \Ftwo{k-1}$,
$(s+t)^{-1}[\{v\}] = s^{-1}[\{v\}]$ and hence its minimum and/or maximum
is from $x_v(s)$. Thus
for $X \in \cX \dotpl \cV$, $ x= \min, \max$ we have
\begin {equation}\label{6.8neu}
\{x_v(s+t) \such s+t \in X\} \in \cR_{v,x},
\end{equation}
as requested in Equation~\eqref{rightprojection}.

Now we proceed by downwards induction on $v$.

The next value is $v = k-1$ and $x = \min$ or $x = \max$. 
For $s \in \Ftwo{k}$, $t \in \Ftwo{k-1}$, $s< t$, $s +t \in X$, 
we have $x_v(s+t) = x_v(s)$ or
$x_v(s+t) = x_v(t)$,
and there is a set $Z \in \cX \times \cV$ such that for any $(s,t) \in Z$
the decision is the same.
We claim that this decision can only be for the $s$-parts and hence
$\{x_v(s+t) \such s+t \in X\}= \{x_v(s) \such s+t \in X\} \in \cR_{k-1,x}$.
This is seen as follows:
Since any $s +t \in X$ is an element of $\Ftwo{k}$, we have
\begin{equation}\label{disjoint}
  \{x_{v+1}(s+t) \such s+t \in X\} \cap \{x_v(s+t) \such s+t \in X\} = \emptyset.
\end{equation}
For $X \in \cX \dotpl \cV$,  
the left-hand set on the left-hand side of Equation~\eqref{disjoint}, i.e.,
the set $\{x_{k}(s+t) \such s+t \in X\}$,
is already in the ultrafilter $\cR_{k,x}$, by the analysis of the coordinate
$k$ in Equation \ref{6.8neu}.
Hence by Equation \eqref{disjoint},
\[
\{x_{k-1}(s+t) \such s+t \in X\} \not \in \cR_{k,x}
\]
and so $Z\in \cX \times \cV$ did not decide for the $t$-part.
Since there must be a decision on $Z$ as to whether $x_{k-1}(s+t) = x_{k-1}(s)$
or $x_{k-1}(s+t) = x_{k-1}(t)$, we have $x_{k-1}(s+t) = x_{k-1}(s)$ for $s \in Z$ and hence 
\[\{x_{k-1}(s+t) \such s+t \in X\} \in \cR_{k-1,x}.\]

In the next step we verify Equation \eqref{rightprojection} for $v = k-2$
with the same technique, shifted by $-1$ in the $v$-component, and proceed by downwards induction down to $v=1$.
\end{proof}

\begin{lemma}\label{liftofThom}(See \cite[Lemma 2.23]{Todorcevic:Ramsey})
Let $k \geq 2$.  Let $2 \leq j \leq k$. Then $\Tl \colon \gamma(\Ftwo{j}) \to \gamma(\Ftwo{j-1})$ is a continuous onto $\dotpl$-homomorphism.
\end{lemma}

So we use freely $T(\cU \dotpl \cV) = T(\cU) \dotpl T(\cV)$.

The following is a strengthening of \cite[Lemma 2.24]{Todorcevic:Ramsey}.
It is stated in more generality than we need for the proof of statement (1) at $k$ so that
we can use it also for the proof of statement (3) at $k$.

\begin{lemma} \label{ladderofidempotents}
  Let $\ba \in (\Fk)^\omega(\bar{\cR})$.
Let $\cU_{k-1}$  be a \gmtu\ over $\Ftwo{k-1}$, such that $\TFU_{k-1}(T[\ba]) \in \cU_{k-1}$ and
$\minl_i(\cU_{k-1}) = \cR_{1+i, {\min}}$, $i = 1, \dots, k-1$, and same for $\maxl$ and
such that $\cR_{1,{\min}}$, $\cR_{1, \max}$ are Ramsey ultrafilters that are  nnc $\Phi(\cU_{k-1})$.
In addition we let $\cS_1$, $ \cS_2$
 be two ultrafilters that are nnc $\Phi(\cU_{k-1})$ and nnc to $\cR_{1,{\min}}$, $\cR_{1,{\max}}$.
\footnote{ The \gmtu\ $\cU_{k-1}$ is given by the
  by the induction hypothesis,
 (2) at $k-1$ with shifted indices in the ultrafilters.}
Let $h$ be a finite-to-one function. The set
\begin{equation}
  \begin{split}
S_{k} := & \{\cX \in \gamma^{\rm good}(\Fk,\bar{\cR}) \such \Tl(\cX)
= \cU_{k-1} \wedge \\
& \bigwedge_n \TFU_k(\ba) \in \cX \wedge\\
& \bigwedge_{j=1,2}(\exists E_j \in \cS_j)( \exists X \in \cX )(h[\set(X)] \cap h[E_j] = \emptyset)
\}
\end{split}
\end{equation}
    is a closed non-empty subset of
    $\gamma(\Ftwo{k})$.
\end{lemma}

\begin{proof} In case $k = 1$, just let $\cU_{0}$ be the
  principal ultrafilter $\{f\}$ with $f(n) = 0$ for any $n$, and we leave out $B$ from the requirements below. Let $\cU_{k-1}$ be given by the induction hypothesis, be as in (2) at $k-1$ for $\cR_{1+i, x}$, $i = 1, \dots k-1$, $x = {\min}, {\max}$, such that
$\cR_{1,x}$ is nnc $\Phi(\cU_{k-1})$ and such that $\TFU_{k-1}(T[\ba]) \in \cU_{k-1}$. Let $h$ be a finite-to-one function.

  Let $B \in \cU_{k-1}$, let for $1 \leq i \leq k$ sets
  $X_i \in \cR_{i,\min}$, $Y_i \in \cR_{i,\max}$ be given. The requirement\\
$R(X_1, \dots, X_k, Y_k, \dots, Y_1, B, h)$ says
\begin{equation*}\begin{split}
    &  (\exists^{\rm min-unb} s)(s_n \in \TFU_{k}(\ba_n))(\forall i = 1, \dots, k) \\
    & (\min_i(s) \in X_i \wedge \max_i(s) \in Y_i \wedge T(s) \in B \wedge\\
& \bigwedge_{j=1,2}
(\exists E_j \in \cS_j)(h[\set(s)] \cap h[E_i]) = \emptyset).
\end{split}\end{equation*}
Again the requirements are closed under finite conjunctions,\footnote{Note that $h$ is fixed. In order to look at all $h$ simultaneously, it is best to work with two finite-to-one functions $h_1$, $h_2$. For any finitely many finite-to-one functions $h^j$, $j < J$, there are two finite-to-one functions $h_1$, $h_2$ such that for any $n \in \omega$, for any $j < J$, $(h^{j})^{-1}[\{n\}] \subseteq
h_1^{-1}[\{n\}]$ or $\subseteq h_2^{-1}[\{n\}]$. Such a compactness argument
is not used here.}
  hence
it suffices that any single requirement can be fulfilled.

We take $\bb \in \cU_{k-1}\cap B$ such that
$(\forall n)(\forall i = 1, \dots, k-1) (\min_i(b_n) \in X_{i+1} \wedge
\max_i(b_n) \in Y_{i+1})$.
Then look at where the sets $X_1$, $Y_1$ lie outside $\set(\bb)$.
We take $X_1^1 \subseteq X_1$, $X_1^1 \in \cR_{1,{\min}}$ and
$Y_1^1 \subseteq Y_1$, $Y_1^1 \in \cR_{1,{\max}}$ and $\bc \strk \bb$, $\bc \in \cU_{k-1}$ such that $X_1^1$ , $Y_1^1$, $\set(\bc)$ are pairwise disjoint.
We can do this since $\cR_{i,x}$ nnc $\Phi(\cU)$.
\nothing{
Since $\cR_{1,x}$ are $Q$-points we find $X_1^2 \subseteq X_1^1$, $X_1^2 \in \cR_{1,{\min}}$ and $Y_1^2 \subseteq Y_1^1$, $Y_1^2 \in \cR_{2,{\min}}$
such that each interval [\min(\supp(c_n)), \min(\supp(c_{n+1})))
is met by at most one point of $X_1^2$, $Y_1^2$. 
}
Let $r_0\in \omega$ be given.
We take $x_0 \in X^1_1$, $ x_0 > r_0$, $r(0) \geq r_0$  such that $x_0 < \supp(c_{r(0})$ and we take $y_0 \in Y^1_1$
that lies past $c_{r(0)}$ and let $r_1= y$ and such that
\[(\exists d= d_{r(0)} \in \TFU_k(\ba)(\min_1(d) = x_0 \wedge T(d) = c_{r(0)}) \wedge \max_1(d) = y_0).\]
Such at $d $ exists since $\ba \in (\Fk)^\omega(\bar{\cR})$.
We choose the next triple
$(x_1, c_{r(1)}, y_1)$ with $r_1$ in the place of $r_0$, and so on.
Since $\cS_1$, $\cS_2$ are ultrafilters, $\cS_1 \cap \cS_2$ is not meager
and hence contains a set $E_1 = E_2$ such that $h[E_1]$ that is disjoint from infinitely many of the sets $h[(r_m, r_{m+1}]]$, $m \in \omega$.
We choose such an $m$.
Then $s=d_{r(m)}$ fulfils the requirement.
Thus we showed that $S_{k}$ is not empty.
\end{proof}

\begin{remark} The version with the starting set $\TFU_k(\ba)$ and  the extra ultrafilters $\cS_1$, $\cS_2$
  will be used for the proof of (3) at $k$. For (1) at $k$, we can drop the
  $\cS_j$ and replace $\TFU_k(\ba) $ by the full  $\Fk$.
  \end{remark}
\medskip
 
Now we are ready to resume the thread of Todorcevic's proof of \cite[Lemma 2.24]{Todorcevic:Ramsey}.

\begin{definition}\label{minimalideals}
  For 
  $\cU, \cV \in \bigcup_{j=1}^k \gamma(\Ftwo{j})$, we let
  $\cV \leq \cU$ if $\cU \dotpl \cV =\cV \dotpl \cU = \cV= \cV \dotpl \cV$. So $\cV$ is idempotent and $\cV$ is stronger and richer than $\cU$.
\end{definition}

The relation $\leq$ is a preorder on $\bigcup_{j=1}^k \gamma(\Ftwo{j})$.

\begin{lemma}\label{ladderofidempotents2} Let $\cU_{k-1}$ be as in Lemma~\ref{ladderofidempotents}. There is an idempotent ultrafilter $\cU^{\rm idem} \in \gamma^{\rm good}(\Fk, \bar{\cR})$ such that $\cU^{\rm idem} \leq \cU_{k-1}$ and $\Tl(\cU^{\rm idem}) = \cU_{k-1}$.
\end{lemma}

\begin{proof}
We let $S_k =\{\cX \in \gamma^{\rm good}(\Fk,\bar{\cR}) \such 
\Tl(\cX) = \cU_{k-1}\}$. By Lemma \ref{ladderofidempotents}, $S_k \neq \emptyset$.
By Lemma~\ref{closure_under_dot+} the set
$\{\cX \dotpl \cU_{k-1} \such \cX \in S_k\}$ is
a subset of $\gamma^{\rm good}(\Fk,\bar{\cR})$.
The set
$\{\cX \dotpl \cU_{k-1} \such \cX \in S_k\}$ is a non-empty closed subsemigroup,
since the sum $\cV \dotpl \cU_{k-1} \dotpl \cW \dotpl \cU_{k-1} $ belongs to
$S_k \dotpl \cU_{k-1}$.
This is seen as follows:
We have
    \[ \Tl(\cV \dotpl \cU_{k-1} \dotpl \cW) = \cU_{k-1} +\Tl\cU_{k-1} + \cU_{k-1} = \cU_{k-1}, 
    \]
where the latter equation holds by induction hypothesis.
Then we pick an idempotent ultrafilter $\cW  \in S_k \dotpl\cU_{k-1}$, say $\cW = \cV \dotpl \cU_{k-1}$ with $\cV \in S_k$.  Let $\cU^{\rm idem} = \cU_{k-1} \dotpl \cV \dotpl \cU_{k-1}$.
    Then $\Tl\cU^{\rm idem} = \cU_{k-1}$ and
    $\cU^{\rm idem} = \cU_{k-1} \dotpl \cW$,
being  a sum of idempotent ultrafilters, is  idempotent.

Moreover $\cU^{\rm idem} \dotpl\cU_{k-1} = \cU_{k-1} \dotpl \cU^{\rm idem} = \cU^{\rm idem}$, so $\cU^{\rm idem} \leq \cU_{k-1}$.
By Lemma~\ref{closure_under_dot+} (and a relative thereof for additions from the righthand-side) we have  $\cU^{\rm idem} \in \gamma^{\rm good}(\Fk, \bar{\cR})$.
\end{proof}

   
Here is the penultimate step of the proof of statement (1) of Theorem~\ref{gowersown}:
We  
pick an idempotent ultrafilter
$\cU^{\rm idem}$ as in the previous lemma. Let $P$ be a piece of the given partition of $\Fk$ such that $P \in \cU^{\rm idem}$.
Now by induction on $n$ we build a tree $Tr$ of finite increasing sequences  $x_0, x_1, \dots , x_n$ of elements of $\Fk$
    and branching sets
 $A_0 \supseteq A_{1,x_0} \supseteq \dots \supseteq A_{n,x_0, \dots x_{n-1}} \in \cU_k$
    with the following properties for any $n \in \omega$, 
    \begin{myrules}
    \item[(a)] $A_0 = P$.
   \item[(b)] $x_n \in A_{n,x_0,x_1, \dots, x_{n-1}}$ and $T^{(k-\ell)}[A_{n,x_0, \dots, x_{n-1}}]= A_{n,x_0, \dots, x_{n-1}}^\ell$ for $1\leq \ell \leq k$. So $A_{n,x_0, \dots , x_{n-1}} = A^k_{n,x_0, \dots, x_{n-1}}$. 
    \item[(c)] $(\cU^{\rm idem} x)(T^{(k-i)}(x_n) + T^{(k-j)}(x) \in A^{\max\{i,j\}}_{n,x_0, \dots, x_{n-1}})$ for
      $1 \leq i,j \leq k$.
       \item[(d)] The tree $Tr$ has branching sets in $\cU^{\rm idem}$,
         i.e. given $(x_0, \dots ,x_n) \in Tr$ for any $n$, the set
         of immediate tree successors of  $(x_0, \dots ,x_n)$ fulfils:
    \begin{equation}\label{Grigorieffbed}
      \begin{split}
        \{x_{n+1} \such & x_{n+1} \mbox{ could serve in the tree $Tr$} 
  \\ &      \mbox{ as a prolongation of }
        (x_0, \dots, x_n)\}
        = A_{n+1,x_0, \dots, x_n} \in \cU^{\rm idem}.
    \end{split}\end{equation}
    \end{myrules}

    As shown in the proof of \cite[Theorem 2.22]{Todorcevic:Ramsey}, any branch $\bx$ of the tree $Tr$ has $\TFU_k(\bx) \subseteq P$.
    We show that there is a branch $\bc $ of $Tr$ such that $\bc \in
    (\Fk)^\omega(\bar{\cR})$.
    We let $\tilde{A}_0 = A_0$ and or $n \in \omega$ we let
    \begin{equation}\label{tildes}
\begin{split}
\tilde{A}_{n+1} = &\bigcap \{A_{m+1,x_0, \dots, x_{m}} \such m \leq n,
    (x_0, \dots, x_{m}) \in Tr,\\
 & \max(\supp(x_{m})) < n+1\}.
\end{split}    
\end{equation}
    Then $\tilde{A}_n \in \cU^{\rm idem}$, and
    they are $\subseteq$-descending.

    \smallskip
    
    Now we come to the last step.

\emph{Final Claim}\\ 
   There is a branch $\bc$ through the Grigorieff tree $Tr$
   such that $\bc \in (\Fk)^{\omega}(\bar{\cR})$.
\footnote{We remark that in the case of ideals (or coideals) over $\omega$, Grigorieff \cite[Section 1]{Grigorieff} worked with
 similar trees.
}

\begin{proof}

We recall another lemma. 

\begin{lemma}\label{projections_lemma} If $\cU_{k-1}$ is a \gmtu\ over $\Ftwo{k-1}$ with projections
and
\[ (\forall i \in \{1, \dots , k-1\})(\dot{x}_i(\cU_{k-1})=\cR_{i+1,x})
\] and $\cR_{1,{\min}}$, $\cR_{1,{\max}}$ are given and both are nnc to $\Phi(\cU_{k-1})$
and $\cV \in S_{k} \dotpl \cU_{k-1}$, $\cV \leq \cU_{k-1}$, then 
$(\minl_1, \Tl, \maxl_1)(\cV) = \cR_{1,{\min}} \times \cU_{k-1} \times \cR_{1,{\max}}$.
\end{lemma} 

\begin{proof} 
No other order is possible. Since $\cV = \cV + \cU_{k-1} = \cU_{k-1} + \cV$, we have
$X \in \cV$ iff $\{s \such t \such s+t \in X \}\in \cU_{k-1}\} \in \cV$
and 
$X \in \cV$ iff $\{s \such t \such s+t \in X \}\in \cV\} \in \cU_{k-1}$.
However, since $\cU_{k-1}$ is a \gmtu, by Lemma~\ref{ordered_stairs_dense} the
high values of $s \in X \in \cU_{k-1}$ are always in the middle.
So only the order  $\cR_{1,{\min}} \times \cU_{k-1} \times \cR_{1,{\max}}$
in the triple projection is possible.
\end{proof}

Now we return to the proof of the final claim.
From  $\tilde{A}_n$, $n \in \omega$, as given in Equation~\eqref{tildes}  we define three families of sets:
\begin{equation}
  \label{setting}\begin{split}
B_n &:= \{\min_1(s) \such s \in \tilde{A}_n\},\\
C_n &:= \{T(s) \such s \in \tilde{A}_n \wedge \min_1(s) \in B_n\},\\
D_n &:= \{\max_1(s) \such s \in \tilde{A}_n \wedge \min_1(s) \in B_n 
\wedge T(s) \in C_n\}.
\end{split}
\end{equation}
Since $\tilde{A}_n \in \cU^{\rm idem} \in \gamma^{\rm good}(\Fk,\bar{\cR})$, we have
 $B_n \in \cR_{1,{\min}}$, $C_n \in \cU_{k-1}$, $D_n \in \cR_{1,{\max}}$.

We prove a relative of Blass' Lemma~\ref{Blass1.1} for triples:
\begin{lemma}\label{finale}
  Fix $k \geq 1$. If $k = 1$, then there is no \gmtu\ $\cU_{k-1}$,
  we replace $\cU_{k-1}$ by the filter $\{0\}$, where $0 $ is a ``block'' with empty support. Suppose we are given 
  \begin{myrules}
  \item[(1)] A \gmtu\ $\cU$ over $\Ftwo{k-1}$ with projections $\cR_{i+1,x} \such 1 \leq i \leq k-1\ra$.  
  \item[(2)] $\cR_{1,x}$, Ramsey and nnc $\Phi(\cU_{k-1})$, $x = \min, \max$, and
    elements $X_1\in \cR_{1,x}$ and $Y_1 \in \cR_{1,\max}$,
  \item[(3)] a descending sequence $\tilde{A}_n$ as in Equation~\eqref{tildes}
of subsets of an idempotent ultrafilter $\cU^{\rm idem}\leq \cU_{k-1}$, $\Tl(\cU^{\rm idem}) = \cU_{k-1}$ and
$\cU^{\rm idem} \in \gamma^{\rm good}(\Fk,\bar{\cR})$ contains $P$. 
\item[(4)] $B_n$, $C_n$, $D_n$ are defined from $\la \tilde{A}_n \such n < \omega \ra$ as in Equation~\eqref{setting}.  
  \end{myrules}

 \emph{Then} there is some $\TFU_k(\bb) \in \cU_{k-1}$ and there is
$X \in \cR_{1,{\min}}$, $Y\in \cR_{1,{\max}}$ with the following properties
\begin{myrules}
\item[(a)]
the blocks of $X$, $\rge(\bb)$ and $Y$ lie as follows
 \[
 x_0 < b_0 <  y_0 < x_1 < b_1 < y_1 < \dots
\]
\item[(b)] the following three conditions on diagonal lower bounds 
hold 
\begin{eqnarray*}
(\forall n)& (x_n \in B_{y_{n-1}}),\\ 
(\forall n)& (b_n \in C_{x_n}),\\ 
(\forall n)& (y_n \in D_{\max(\supp(b_n))}).
\end{eqnarray*}
\end{myrules}
\end{lemma}

\begin{proof} \emph{First step:}\\
We take on each of the three parts diagonal lower bounds.
We take $X_{0.5} \in \cR_{1,{\min}}$ such that
\[\forall x,y \in X_{0.5} (y<x \rightarrow x \in B_y).\]

We take $Y_{0.5} \in \cR_{1,{\max}}$ such that  
\[\forall x,y \in Y_{0.5} (y<x \rightarrow x \in D_y).\]

We take  $\bb_{0.5} \in \cU_{k-1}$ such that
\[\forall n < m (b_{0.5,m} \in C_{\max(\supp(b_n^1))}).
\]

Then we take $X_1 \subseteq X_{0.5}$, $Y_1 \subseteq Y_{0.5}$, $\bb_1 \strk \bb_{0.5}$,
$X_1 \in \cR_{1,{\min}}$, $Y_1 \in \cR_{1,{\max}}$, $\bb_1 \in \cU_{k-1}$
such that $X_1$, $X_2$ and $\set(\bb_1)$ are pairwise disjoint.
Since $\cR_{1,{\min}}$ and $\cR_{1,{\max}}$ are nnc and both are nnc to
$\Phi(\cU_{k-1})$, there are such strengthenings. 
The diagonal conditions hold also for $X_1$, $\bb_1$, $Y_1$.

\emph{Second step:} \\
We start working towards the triple alternation pattern in clause (a).
We take $X_2 \subseteq X_1$, $X_2 \in \cR_{1,{\min}}$, $\bb_{1.5} \in \cU_{k-1}$,
$\bb_{1.5} \strk \bb_1$, $Y_{1.5} \in \cR_{1,{\max}}$, $Y_{1.5} \subseteq Y_1$, such that 
\begin{equation}\label{alternating_1}
  (\forall x < x' \in X_2)(\exists n)(\exists y \in Y_{1.5})(x < \supp(b_{1.5,n})< y < x').
\end{equation}
Since $\cR_{1,{\min}}$ is a $Q$-point and nnc $\Phi(\cU_{k-1})$ and nnc to $\cR_{1,{\min}}$, we find such a sparse set $X_2$ and $Y_{1.5}$ and $\bb_{1.5}$.

We take $\bb_2 \strk \bb_{1.5}$, $\bb_2 \in \cU_{k-1}$, such that 
\begin{equation}
  \label{alternating_2}
  (\forall n < m)(\exists x \in X_2)(\exists y \in Y_{1.5})(\supp(b_{2,n}) < y<x < \supp(b_{2,m}).
\end{equation}
Such a $\bb_2$ can be found via the Taylor-colouring law of $[\TFU_{k-1}(\bb_{1.5})]_<^2$
that sends $s< t $ to $1$ if there are the desired elements of $X_2$ and of $Y_{1.5}$ between $s$ and $t$ and to $0$ otherwise. We find  $\bb_2\strk \bb_{1.5}$, $\TFU_{k-1}(\bb_2) \in \cU_{k-1}$ such that
the monochromatic colour is 1.

We take $Y_2 \subseteq Y_{1.5}$, $Y_2 \in \cR_{1,{\max}}$, such that
\begin{equation*}
  (\forall y < y' \in Y_2)(\exists n)(\exists x \in X_2)(y < x<\supp(b_{2,n})< y').
\end{equation*}

\emph{Third step:}\\
We fill up a tiny bit: If between two adjacent blocks of $\bb_2$ there is
only an element of $X_2$ but no larger element of $Y_2$, we add such a point from $Y_{1.5}$ to $Y_2$
and thus get, when doing this for all situations a set $Y_3 $ such that
$Y_2 \subseteq Y_3 \subseteq Y_{1.5}$.
By Equation \eqref{alternating_2} there are such elements of $Y_{1.5}$ such that
the alternating pattern is established.

If between two adjacent points of $X_2$ there is
only an element of $Y_{1.5}$ but no element of $Y_2$ or if there is only an element
of $\bb_{1.5}$ but no element of $\bb_2$, we add, if applicable, such a point from $Y_1$ to $Y_3$
and thus get, when doing this for all situations a set $Y_4 $ such that
$Y_3 \subseteq Y_4 \subseteq Y_{1.5}$, and  we add, if applicable, such a block from $\rge(\bb_{1.5})$ to $\bb_2$
and thus get, when doing this for all situations a set $\bb_3$ such that
$\bb_2 \strk \bb_3 \strk \bb_{1.5}$, $\bb_2 \in \cU_{k-1}$.
By Equation~\eqref{alternating_1} such a filling up is possible.

We let $X_4 = X_2$, $\bb_4 = \bb_3$.
Thus we have an alternating chain of enumerating $X_4 \cup \rge(\bb_4) \cup Y_4$
as $x_0 < b_{0} < y_0 < x_1 < \dots$. 

\emph{Fourth step:}\\
 We ``attack the possible exceptions'' to the diagonalisation
requirements in clause (b).

Exceptions for $x$: Suppose that in the common increasing enumeration of
$X_4 \cup \rge(\bb_4) \cup Y_4$ there is an element $x\in X_4$  such that
for no $x ' \in X_{1.5}$, 
$y < x' < x < b_{4,k}$.
We map the last $\bb_4$-block before $x$ to $x$ by a finite-to-one
function that sends the whole interval $[\min(b_{4,n}), x]$ to $\min(b_{4,n})$.
We do this for all such $x$.
We define $f$ arbitrarily on the other points so that it will be a finite-to-one
function. Then we use that $\cR_{1,{\min}}$, $\cR_{1,{\max}}$ nnc $\Phi(\cU_{k-1})$ 
and take
$\bb_5 \strk \bb_4$, $\bb_5 \in \cU_{k-1}$ and $X_5 \subseteq X_4$, $X_5 \in \cR_{1,{\min}}$
and $Y_5 \subseteq Y_4$, $Y_4 \in \cR_{1,{\max}}$ such that
$h[Y_5] \cap h[\set(\bb_5)] = \emptyset$, $h[X_5] \cap h[\set(\bb_5)] = \emptyset$. Then for the common enumeration of $X_5$, $\bb_5$ and $Y_5$ we have: Any $x \in X_5$ that is preceded by an element of $Y_5$ or one of $\bb_5$ is shielded by an $x' \in X_1$
from them and hence fulfils the condition that corresponds to
the first equation in (b) for $X_5$, $\bb_5$, $Y_5$.
However clause (a) is already ruined.

Exceptions for $b_{5,k}$:  Suppose that in the common increasing enumeration of
$X_5 \cup \rge(\bb_5) \cup Y_5$ there is an element $b_{5,k}\in \rge(\bb_5)$
such that for no $b' \in \bb_{1}$ we have for the largest $x \in X_5$ preceding $b_{5,k}$,
$y<x < b' < b_{5,k}$. (There might be more $x$'s in the pattern, since in
the first manipulation of the foruth step we may have destroyed the alternating pattern.
It is only important to do the procedure with a smaller $y$.)
We take a finite-to-one
function that sends the whole interval $[y, \max(b_5,k)]$ to $\min(b_{k,k-1})$.
We do this for all such $b_{5,k}$.
We define $f$ arbitrarily on the other points so that it will be a finite-to-one
function. Then we use that $\cR_{1,{\min}}$, $\cR_{1,{\max}}$ nnc $\Phi(\cU_{k-1})$ and take
$\bb_6 \strk \bb_6$, $\bb_6 \in \cU_{k-1}$ and $X_6 \subseteq X_5$, $X_6 \in \cR_{1,{\min}}$
and $Y_6 \subseteq Y_5$, $Y_6 \in \cR_{1,{\max}}$ such that
$h[Y_6] \cap h[\set(\bb_6)] = \emptyset$, $h[X_6] \cap h[\set(\bb_6)] = \emptyset$.

Exceptions for $y$: Suppose that in the common increasing enumeration of
$X_6 \cup \rge(\bb_6) \cup Y_6$ there is an element $y\in Y_6$
such for no $y' \in Y_{1}$ we have that for the latest block $b_{6,n}$ before
$y$ the constellation
$x <b_{6,n} < y' < y < x'$, $x, x' \in X_6$. Again the remark in parantheses applies. 
We map the last interval $[x,y]$ before $y$ to $x$ by a finite-to-one
function.
We do this for all such $y$.
We define $f$ arbitrarily on the other points so that it will be a finite-to-one
function. Since $\cR_{1,{\min}}$, $\cR_{1,{\max}}$ are nnc $\Phi(\cU_{k-1})$ there are
$\bb_7 \strk \bb_6$, $\bb_7 \in \cU_{k-1}$ and $X_7 \subseteq X_6$, $X_7 \in \cR_{1,{\min}}$
and $Y_7 \subseteq Y_6$, $Y_7 \in \cR_{1,{\max}}$ such that
$h[Y_7] \cap h[\set(\bb_7)] = \emptyset$, $h[X_7] \cap h[\set(\bb_7)] = \emptyset$.

By the thinning out in the fourth step we ensured that 
$X_7$, $Y_7$, $\bb_7$ respect clause (b) in the following form:
\begin{equation}\label{clause(b)}
  \begin{split}
&(\forall x \in X_7)( \forall y \in Y_7) (y < x \rightarrow x \in B_y)\\
&(\forall b \in \rge(\bb_7) )(\forall x \in X_7) (x < \supp(b) \rightarrow
b \in C_x)\\
&(\forall y \in Y_7)( \forall b \in \rge(\bb_7)) (\supp(b) < y \rightarrow
y \in D_{\max(\supp(b))}).
  \end{split}
\end{equation}
However, clause (a) is ruined. Note that the diagonal laws
from Equation~\ref{clause(b)} are preserved
under further thinning out.

\emph{Fifth step:}\\
 We restore the triple alternation pattern again. We repeat the steps
two and three, and thus get $X_8\subseteq X_7$, $Y_8\subseteq Y_7$, $\bb_8\strk \bb_7$  that respect (a) and (b).
\end{proof}

We conclude the proof of part (1) of Theorem~\ref{gowersown}:
Let $X$, $\bb$, $Y$ be as in Lemma~\ref{finale}. By Equation \eqref{setting}, for each $n$, there is 
 \[
c_n \in \tilde{A}_{\max(\supp(c_{n-1}))} \mbox{ s.t. } (\min_1(c_n) = x_n \wedge T(c_n) = b_n \wedge \max_1(c_n) = y_n).
\]
By the definition of the $\tilde{A}_n$ we have $\TFU_k(\bc) \subseteq P$. 
Since
\begin{equation*}\begin{split}
  &\bb \in \cU_{k-1}  \in  \gamma^{\rm good}(\Ftwo{k-1}, \la \cR_{1+i,x}\such 1 \leq i \leq k-1, x = {\min},{\max}\ra),\\
  &\{\min_1(c_n) \such n < \omega\}  =  X  \in  \cR_{1,{\min}},\\
  & \{\max_1(c_n) \such n<\omega\}  =  Y   \in  \cR_{1,{\max}}, \mbox{ and }\\
   & \Phi(\cU_{k-1}) \mbox{ is nnc }  \cR_{1,{\min}}, \cR_{1,{\max}},
  \end{split}
\end{equation*}
we have $\bc \in (\Fk)^\omega(\bar{\cR})$. So $\bc$  is branch of the Grigorieff tree with the
desired properties.
Thus Theorem~\ref{gowersown} part (1) is proved.

\medskip

We prove part (2). Again we encounter the phenomenon (as in the proof
of the Hindman colouring in the $\M(\cU)$-extension in Lemma~\ref{Hindmansuccessor})
  that we almost proved the existence of diagonal lower bounds already during the
  proof of the colouring theorem.

\begin{lemma}\label{finale2}
\hfill
  \begin{myrules}
  \item[(1)] a $\strk$-descending sequence $\ba_n$, $n \in \omega$,
$\ba_n \in (\Fk)^{\omega}(\bar{\cR})$,
\item[(2)] 
  \begin{equation}
    \label{setting2}
    \begin{split}
B_n &:= \{\min_1(s) \such s \in \TFU(\ba_n)\}\\
C_n &:= \{T(s) \such s \in \TFU(\ba_n) \wedge \min_1(s) \in B_n\}\\
D_n &:= \{\max_1(s) \such s \in \TFU(\ba_n) \wedge \min_1(s) \in B_n 
\wedge T(s) \in C_n\}
\end{split}
\end{equation}
  \end{myrules}

  \emph{Then} there is some $\bb \strkast T(\ba_n)$, $\bb \in (\Ftwo{k-1})^\omega
  (\la \cR_{1+i, x} \such 1 \leq i \leq k-1, x={\min}, {\max}\ra)$ and there is
$X \in \cR_{1,{\min}}$, $Y\in \cR_{1,{\max}}$ with the following properties
the blocks of $\bb$ and $X$ and $Y$ lie as follows
 \[
 x_0 < b_0 <  y_0 < x_1 < b_1 < y_1 < \dots
\]
and the following three conditions on diagonal lower bounds 
hold 
\begin{equation}\label{alternate}
(\forall n) \bigl(x_n \in B_{y_{n-1}} \wedge b_n \in C_{x_n} \wedge y_n \in D_{\max(\supp(b_n))})
\end{equation}
\end{lemma}

\begin{proof} This is proved as in Lemma \ref{finale}. In the first step, we take separate
diagonal lower bounds, call them as there. For the middle part $\bb_1 \in (\Fk)^\omega(\la \cR_{1+i,x}
\such 1 \leq i \leq k-1, x = {\min},{\max}\ra)$ we choose by (3) at $k-1$ a \gmtu\ $\cU_{k-1} \in \gamma^{\rm good}(\Ftwo{k-1}\la \cR_{1+i,x}
\such 1 \leq i \leq k-1, x = {\min},{\max}\ra)$ such that $\Phi(\cU_{k-1})$ is nnc $\cR_{1,x}$ and such that $\TFU_{k-1}(\bb_1) \in \cU_{k-1}$. Then we can go on as in the proof of Lemma \ref{finale}. \end{proof}

Equation \ref{alternate}
yields 
 a diagonal bound $\bc \in (\Fk)^\omega(\bar{\cR})$ of $\ba_n$, $n \in \omega$.
This finishes the proof of statement (2) at $k$. 
 
\smallskip

Now we turn to (3) at $k$:

 Let $\cS_1$ and $\cS_2$ be two new ultrafilters, nnc to any of the
 $\la \cR_{i,x} \such i =1, \dots, k, x = {\min},{\max}\ra$.
 Using \CH, we enumerate via $\la X_\alpha,h_\alpha \such \alpha < \omega_1\ra$
 any subset $X_\alpha $ of $\Fk$ and any finite-to-one function $h_\alpha$ such that
 each pair appears cofinally often in the enumeration.
 By induction on $\alpha< \omega_1$ we choose $\ba_{\alpha} \in
 (\Fk)^\omega(\bar{\cR})$ such that
   \begin{myrules}
     \item[(1)]
 $\ba_{\alpha+1}$ decides $X_\alpha$ and
       \item[(2)] for $i =1,2$, $\exists E_i \in \cS_i$, $h_\alpha[\set(\ba_{\alpha+1})] \cap
         h_\alpha[E_i] = \emptyset$.
         \item[(3)] $(\forall \beta < \alpha)(\ba_{\beta} \strkast \ba_\alpha)$.
   \end{myrules}
   Then we let $\cU$ be generated by $\{\TFU_k(\ba_{\alpha} \such \alpha < \omega_1\}$.

   Given $\ba_\beta$, $\beta < \alpha$, we first
   take a $\strkast$-lower bound $\ba_\alpha$ of them as guaranteed by (2) at $k$. We let $h = h_{\alpha}$
   By induction hypothesis there is a \gmtu\ $\cU_{k-1}$ avoiding $\cS_i$
   and containing $\TFU_{k-1}(T[\bc_{\alpha}])$. In step 0, just drop all clauses about $\ba_\alpha$.
   Then we build $S_k$ in Lemma \ref{ladderofidempotents} (this time
   in the form with the $\cR_{1,x}$ and the $\cS_j$ and the starting set
   $\TFU_k(\ba)$, note that this lemma uses only (3) at $k-1$), i.e., the requirement pertains to
   $h_\alpha$, and $\TFU_k(\ba)$ such that 
   for any $X \in \cX \in S_k$,  $(\exists E_i \in \cS_i)_{i=1,2} (h[\set[X]] \cap h[E_i] = \emptyset)$ and $\TFU_{k}(\ba_{\alpha}) \in \cU^{\rm idem}$.
   We use, that $h$-avoiding $\cS_i$ is preserved by $\dotpl$ to find an idempotent ultrafilter $\cU^{\rm idem} \leq \cU_{k-1}$, $\cU^{\rm idem} \in S_k + \cU_{k-1}$
   such that $\TFU_k(\ba_{\alpha}) \in \cU^{\rm idem}$ and $\cU^{\rm idem} $ $h$ avoids $\cS_j$. Then we take $Y_\alpha$ to be
 $X_\alpha$, if $X_\alpha \in \cU^{\rm idem}$ or $Y= \Fk \setminus X_\alpha$ otherwise.
 Now we build all the Grigorieff tree with
 $\tilde{A}_n \in \cU^{\rm idem}$ and $\tilde{A}_n \subseteq \TFU_k(\ba_\alpha)$.
 A branch of it is taken as $\ba_{\alpha+1}$. Then $\TFU_k(\ba_{\alpha+1}) \subseteq Y_\alpha$.

 Thus we finish the proof of part (3) at $k$ of Theorem \ref{gowersown}.
\end{proof}

\begin{question}\label{baumgartnerstyle}
    We do not know a proof in the Baumgartner \cite{baumgartner:hindman} style.
The analogous question for Gowers' $\Fk$ theorem is asked in \cite{Todorcevic:Ramsey}.
 \end{question}

Now we show that under \CH\ or under \MA\ or after forcing with the pure part
of $\GM(\bar{\cR})$ there is a \gmtu.

\begin{corollary} \label{adequate}
  Let $k$, $\bar{\cR}$ be as in Def.~\ref{PP}. The family
  $(\Fk)^\omega(\bar{\cR})$ is Gowers-adequate.
\end{corollary}

\begin{corollary} \label{gmtu_existence}
  Let $k$,  $\bar{\cR}$ be as in Def.~\ref{PP}.
  Under \CH\ or under \MA($\sigma$-centred) there is an \gmtu.
\end{corollary}
\begin{proof} We perform a $\strkast$-downward construction of length
  $\gc$ with lower bound in the limit steps that are either given by \CH\ or
  provided by \MA($\sigma$-centred).
  In the successor we use Theorem~\ref{gowersown} to decide whether the next subset of $\Fk$ is in or out with the aid of Theorem \ref{gowersown}.
  \end{proof}

\begin{definition}\label{higher_core}
  Let $k \geq 2$ and let $1\leq i < k$.
  Let $\cU$ be an ultrafilter over $\Fk$. We let
  \[\Phi_{\geq i+1}(\cU) = \{\{s^{-1}[\{i+1, \dots, k\}] \such s \in X\} \such X \in \cU\},\]
  and call $\Phi_{\geq i+1}(\cU)$ the core strictly above $i$.
\end {definition}

The following generalises \cite[Theorem 38]{blass-topap2009} and shows
that in any \gmtu\ over $\Fk$ the $2k$ pairwise nnc Ramsey projections
can be found.
The additional information on the higher cores
can be seen as a justification of the inductive procedure in the proof \ref{gowersown}.

 \begin{theorem}\label{structure}
   For any \gmtu\ $\cU$ over $\Fk$ the following holds:
   \begin{myrules}
   \item[(1)]
     For each $1 \leq i  \leq k$, $\cR_{i,{\min}}$ and $\cR_{i,{\max}}$ are nnc and
     for $i< k$, $\cR_{i,{\min}}$ and $\cR_{i,{\max}}$  are nnc $\Phi_{\geq i+1}(\cU)$.
\item[(2)]
The projections $\minl_1(\cU)$, \dots , $\minl_k(\cU)$, $\maxl_k(\cU)$, \dots, $\maxl_1(\cU)$ are  pairwise non-nearly coherent Ramsey ultrafilters over $\omega$.
 \item[(3)] All cores $\Phi_i(\cU)$ are nearly coherent.
        \end{myrules}
 \end{theorem}

\proof
In \cite[Proposition~3.9]{Blass:ufs-hindman} it is shown that all minima and all maxima are
Ramsey ultrafilters. It remains to show that they are pairwise non-nearly coherent. 

We prove (1). 
In the first step we consider the case of $\minl_j(\cU)$, $\maxl_j(\cU)$.
Almost verbatim as in \cite[Section 6]{blass-topap2009} it is shown that
for any $j \in \{1, \dots, k\}$, the two Ramsey ultrafilters
$\minl_j(\cU)$ and
$\maxl_j(\cU)$ are not nearly coherent.
The reasoning uses only the closure of $\TFU_k(\bb)$ under $+$ from Def.~\ref{ksetting1}(7) and does not refer to the Tetris function.

We show that for any  $1 \leq i < \leq k$, $\minl_i(\cU)$ and $\Phi_{\geq i+1}(\cU)$ are not nearly coherent. The crucial step is the Tetris step that does not exist for the Hindman space.
For a contradiction, we assume that $f$ is finite-to-one and $f(\minl_i(\cU)) \supseteq 
f(\Phi_{\geq i+1}(\cU))$. 
W.l.o.g. we assume that $f$ is non-decreasing and surjective. We let $I_n = f^{-1}[\{n\}]$.
We colour $\Fk$ as follows. Let $c(s) = 0$ if the $k$ such that $\min(s^{-1}[\{i\}]) \in I_k$
is equal to the $k'$ such that 
$s^{-1}[\{i+1, \dots k\}]$ meets $I_{k'}$; and $1$ otherwise.
We choose a  $c$-monochromatic $\TFU_k(\ba) \in \cU$.
We show that the colour is 1. For any $r < t$, $r,t \in \rge(\ba)$   whose supports do not meet the same $I$-interval  we have for any $j \geq i+1$ that $j \not\in \rge ( T^{(k-i)}(r))$
but  $i \in \rge ( T^{(k-i)}(r))$ and hence
\[c( T^{(k-i)}(r) + t) = 1.
\]
However, by the assumed near coherence there are infinitely many $n$ such that
\[
I_n \cap \{\min(a_m^{-1}[\{i\}])\such m \in \omega\} \cap
(\bigcup\{a_\ell^{-1}[\{i+1, \dots, k\}]\such \ell \in \omega\}) \neq \emptyset.
\]
We take $n$, $m$, $\ell$ such that
\[\min(a_m^{-1}[\{i\}]) \in I_n , a_\ell^{-1}[\{j\}] \cap  I_n \neq \emptyset.
\]
Then in case $m < \ell$ we have $c(a_m + a_\ell) = 0$ and in case
$m = \ell$ we have $c(a_m) = 0$. Contradiction.

Similarly one shows the statement about $\cR_{i,{\max}}$.

(2): Statement (1) implies (2), since for $j > i$, $\minl_j(\cU), \maxl_j(\cU) \supseteq \Phi_{\geq i+1}(\cU)$.

(3)
Now we consider the cores.
We let
$\ba \in (\Fk)^\omega$ such that $\TFU_k(\ba) \in \cU$.
Then there is some finite-to-one $h$ such that 
\[(\forall n) (h[a_n^{-1}[\{1\}]] = \dots = h[a_{n}^{-1}[\{k\}]]).
\]
Any set in $\cU$ is of the form $\TFU_k(\bb)$ for some $\bb $ compatible with $\ba$, we can assume $\bb \strk\ba$.
Then by the definition of the condensation order, in each block of $\bb$ there
must be one block  of $\ba$ that does not drop by the Tetris function.
Hence 
\[(\forall n) (\forall i =1, \dots, k)
(h[b_n^{-1}[\{i\}] \cap h[b_{n}^{-1}[\{j\}] \neq \emptyset).
\]
Thus $\Phi_i(\cU)$ and $\Phi_j(\cU)$ are nearly coherent non-meagre filters.
\proofend
 
\section{Localised Gowers--Matet Forcing}
 \label{S7_neu}

 Now we harvest a fruit from the Ramsey theoretic work:
 a new family of forcing orders that allows us to arrange the number of near-coherence classes.  Since we use generalised condensations
  by including the Tetris operation, we call the new forcing Gowers--Matet forcing.
 \label{begin_forc}
\begin{definition}\label{GowersMatetforcing}
  Let $k \geq 1$.  Conditions in {\em Gowers--Matet forcing}, $\GM_k$,
  are pairs $(s, \bar{a})$ such that
$s\in \Fk$ and $\ba \in (\Fk)^\omega$ and $\supp(s)<\supp(a_0)$.
Let $p=(s,\ba) \in \GM_k$. The component $s$ is called the \emph{trunk of $p$} and
the component $\ba$ is called the \emph{pure part of $p$}.

  The forcing order is $(t, \bb) \leq_k (s,\ba)$ 
(recall the stronger condition is the smaller one) 
if the following holds:
\begin{myrules}
\item[(a)] $s=t$ or there are $n$, $i_0 < \dots < i_n \in \omega$, and  $j_r \in \{0, \dots, k-1\}$ for $r \leq n$ with at least one $j_r = 0$ such that $t$ is a sum of the form
  \[t = s + T^{(j_0)}(a_{i_0}) + \dots + T^{(j_n)}(a_{i_n}).\]
    \item[(b)]
$\bb \strk \ba$ (see Def.~\ref{collection_on_Tetris}(6)).
\end{myrules}
In the case of $s=t$ we call $(t,\bb)$ a \emph{pure extension of $(s,\ba)=p$}.
\end{definition}

\begin{definition}\label{GowersMatetwithH}
Given a Gowers--Matet-adequate family $\cH \subseteq (\Fk)^\omega$,
the notion of forcing $\GM_k(\cH)$ 
consists of all pairs $(s,\bar{a})$
 such that $\ba\in \cH$.
The forcing order is the same as in the Gowers--Matet forcing.
\end{definition}

\begin{definition}\label{GMkR}
  For $\cH = (\Fk)^\omega(\bar{\cR})$ from the previous
  section we let $\GM_k(\bar{\cR}) = \GM_k(\cH)$.
  \end{definition}

Thanks to Theorem~\ref{gowersown} we have:

\begin{lemma}\label{puredecision} $\GM_k(\bar{\cR})$ has the \emph{pure decision property}, i.e.,
  for any $\varphi \in {\mathcal L}(\in)$, $(s,\ba) \in \GM_k(\bar{\cR})$
  there is $(s,\bb) \leq_0 (s,\ba)$ such that 
$((s,\bb) \Vdash \varphi \vee (s,\bb) \Vdash \neg\varphi)$.
\end {lemma}
\begin{lemma} \label{proper}
  The forcing poset
  $(\GM_k(\bar{\cR}), \leq, (\leq_n)_{n< \omega}))$ fulfils Axiom A and hence is proper.
\end {lemma}

\proof A derivation of Axiom A for the relations $(\leq_n)_n$ from the pure decision property for
tree forcings can be found, e.g., in \cite[Section 7.1]{BJ}.
Properness alone can also be proved as in \cite[2.3, 2.4, 2.5]{BsSh:242}.
The latter proof  uses
the existence of lower bounds, i.e., statement (2) of Theorem~\ref{gowersown}, but not the pure decision property.
\proofend
\begin{definition}\label{Gowers_corek}
Let $\cH$ be a subset of $(\Fk)^\omega$.
\begin{myrules}
\item[(1)] The $\GM_k(\cH)$-generic function
  from $\omega$ to $k+1$ is
  \[\mu = \bigcup\{s \rest (\max(\supp(s))+1) \such (\exists \ba)( (s,\ba) \in G)\}.
  \]
    \item[(2)]
      Again the generic $i$-fibres for $i =1, \dots, k$
      are written as
      \[\mu_i := \mu^{-1}[\{i\}].\]
 \item[(3)]
Any \gmtu\ $\cU$ avoids $\cE$ iff $\Phi(\cU) \not\leq_{\rm RB} \cE$.
\end{myrules}
\end{definition}

Now we investigate the number of  near-coherence classes in models
our our new forcing.
  
\begin{definition}\label{fmu}
 Let $X \in \roth$. We let $f_X(n) = |X \cap n|$.
  \end{definition}

We start with a density argument for evaluating our forcings.

\begin{theorem}\label{fmu(E)=fmu(W)}
Assume \CH.  Let $\bar{\cR}$ be as in Def.~\ref{PP}.
    Let $h \colon \omega \to \omega$ be a finite-to-one function.
    Let $\cE$ and $\cW$ be ultrafilters over $\omega$ such that
    $\cW$,    $\cE \not\geq_{RB} \cR_{i,x}$ for $i =1, \dots, k$, $x ={\min},{\max}$.
    Then
    \begin{equation*}
      \begin{split}
        \GM_k(\bar{\cR}) \Vdash_{\GM_k(\bar{\cR})}
        & f_{\supp(\mu)}(\cW) = f_{\supp(\mu)}(\cE).
              \end{split}
    \end{equation*}
    \end{theorem}

\proof 
Suppose for some $E \in \cE$, $W \in \cW$,  some $(s,\ba) \in \GM_k(\bar{\cR})$,
$(s,\ba) \Vdash f_{\supp(\mu)}[W] \cap f_{\supp(\mu)}[E] \cap [n,\omega) = \emptyset$. 
By Theorem \ref{gowersown}(3), there is a \gmtu\ $\cU_k$ such that
$\TFU_k(\ba) \in \cU$ and such that $\Phi(\cU_k)$ avoids $\cE$, $\cW$.

So there is a $\bb \strk \ba$ such that
$\TFU(\bb) \in \cU_k$ and such that there are infinitely many elements of $E$ outside $\set(\bb)$ and
infinitely many elements of $W$ outside $\set(\bb)$.
Since $\TFU_k(\bb) \in \cU_{k}$, Lemma ~\ref{fitting} implies
that $\bb \in (\Fk)^\omega(\bar{\cR})$. We take some $ e \in (E \setminus \set(\bb)) $ such that $ e > \max\{ n, \max(\supp(s))\}$.
and some $w \in W \setminus \set(\bb)$ such that $w > \max\{n, \max(\supp(s))\}$.
Then we can lengthen $s$ by taking $t > s$ such that \\
$\min(\supp(t)) > w,e$
and $(s+t, \bb) \leq_{\GM_k(\bar{\cR})} (s,\ba)$.
Then $(s+ t, \bb) \Vdash f_{\supp(\mu)}(e) = f_{\supp(\mu)}(w)$.
Contradiction.
\proofend

  \begin{theorem}\label{EisworthforGMk(cR)}
  (Adaption of 
\cite[Theorem~4]{Eisworth}) 
Let $k \geq 1$ and $\bar{\cR}$ be as above and assume that $\cE$ is
a $P$-point with
$\cE \not\geq_{RB} \cR_{i,{\min}}, \cR_{j,{\max}}$ for any $i \in \{1, \dots, k\}$. Then
$\cE$ continues to generate an ultrafilter after we force with $\GM_k(\bar{\cR})$.
  \end{theorem}
  
  \proof
  By Theorem \ref{gowersown}(2) the conditions $(s,\ba)$ such that
  there is a \gmtu\ $\Phi(\cU_k)$ is nnc to $\cE$ (which means
  $\Phi(\cU_k) \not\leq_{rm RB} \cE$) and such that $\ba \in \cU_{k}$ are
  dense in $\GM_k(\bar{\cR}$.
  Since $\GM_k(\bar{\cR})$ has the pure decision property  we can carry out
  the proof literally as in \cite[Lemma 2.3 -- Lemma 2.10]{Eisworth}.
  \proofend

Now we are concerned with the second iterand. The following follows from
an easy density argument.

\begin{lemma}\label{density2}
For any $i \in \{1, \dots, k\}$, $x = {\min},{\max}$, we have
 \[ \GM_k(\bar{\cR}) \Vdash \cR_{i,x} \cup \{\mu_i\} \mbox{ is a filter subbase.}\]
  \end{lemma}

Now we come to the $k$-dimensional version of  Theorem~\ref{vorstufe}. It is crucial that the
Ramsey ultrafilter $\cR_{i,x}$ is not diagonalised. Otherwise we would add a dominating real.
For each $1\leq i\leq k$ the two Ramsey ultrafilters $\cR_{\min,i}$ and $\cR_{\max,i}$ are destroyed by the $i$-fiber $\mu_i$.

\begin{theorem}\label{onestepforMk(cR)}
  Let $k$, $\bar{\cR}$ be as in Def.~\ref{PP}.
  \begin{equation}
    \begin{split}
      \GM_k(\bar{\cR}) \Vdash &(\forall i \in \{1,\dots,k\}, x \in \{{\min},{\max}\})\\
      &(\filter((\cR_{i,x} \cup \{\mu_i\}))^+ \mbox{ is a happy family}\\
      &\mbox{that is nowhere almost a filter}).
          \end{split}
  \end{equation}
  and hence by Mathias' \cite[Prop.~011]{mathias:happy} combined with
  \cite[proof of Theorem 14]{ncf1},
\begin{equation*}
    \begin{split}
      \GM_k(\bar{\cR}) \Vdash &(\exists \la \cR_{i,x}^{\rm ext} \such i=1,\dots,k, x ={\min},{\max} \ra)
      \Bigl((\forall (i,x))\bigl( \cR_{i,x}^{\rm ext} \supseteq (\cR_{i,x} \cup \{\mu_i\}) \wedge
      \\
      &\mbox{$\cR_{i,x}^{\rm ext}$ is a Ramsey ultrafilter that is nnc to $\cE$} \bigr)\\
      & \mbox{ and the extensions are pairwise nnc} \Bigr).
  \end{split}
\end{equation*}
\end{theorem}

\proof 
Theorem~\ref{gowersown1} allows us to transfer the proof of
Theorem~\ref{vorstufe} to $\GM_k(\bar{\cR})$.
\proofendof{\ref{onestepforMk(cR)}}

We rework Section~\ref{S4}, the iteration theory for limit steps, for
iterands of the form $\GM_k(\bar{\cR})$.
We fix $k$ and $\bar{\cR}= \bar{\cR}_0$ as in Def.~\ref{PP}.
We define by induction on $\alpha \leq \omega_2$
a countable support iteration (in the sense of \cite[Def.~III, 3.1]{Sh:f}) $\bP_\alpha = \la \bP_\beta, \GM(\name{\bar{\cR}_\gamma}) \such \beta \leq \alpha, \gamma < \alpha \ra$ such that
for any $\beta < \alpha$,
\begin{equation}\label{inducforMk(cR)}
  \begin{split}
    &\bP_\beta \Vdash (\forall (i,x)) \bigl(\cR_{i,x,\beta} \supseteq
    \bigcup\{ \cR_{i,x,\gamma} \cup \{\mu_{\gamma,i}\} \such \gamma < \beta\} \\
    &\wedge \cR_{i,x,\beta}\mbox{ is a Ramsey ultrafilter that is nnc to $\cE$}\bigr)\\
    &\mbox{ and the $\cR_{i,x,\beta}$, $i=1,\dots,k$, $x = {\min},{\max}$,
    are pairwise nnc}.
  \end{split}
\end{equation}
We use for names the same letters as for the corresponding evaluated names.

Here $\bP_{\beta+1}$ forces that $\mu_\beta$ is the $\GM_k(\bar{\cR_\beta})$-generic real, and $\mu_{\beta,i}$ denotes its $i$-fibre according to Definition~\ref{corek}(8).
In Theorem~\ref{onestepforMk(cR)} we proved that there are extensions in the successor steps.
This guarantees the existence of $\bar{\cR}_{\beta+1}$ with the desired properties.

Now we consider limit  steps $\alpha$.
If $\cf(\alpha) > \omega$, 
we can just take $\bP_\alpha \Vdash \cR_{i,x,\alpha}=\bigcup_{\gamma<\alpha}\cR_{i,x,\gamma}$ and 
the inductive hypotheses will be carried on, since in proper forcing
every real appears at a step of at most countable cofinality, with the only exception that for $\alpha= \aleph_2$ the \CH\ gets lost.
So we concentrate on the hard case, $\cf(\alpha)=\omega$.
\begin{equation}\label{induc2forMk(cR)}
  \begin{split}
  \bP_\alpha \Vdash(\exists \bar{\cR_\alpha})&\bigl(\bigwedge_{i,x}\cR_{i,x,\alpha} \supseteq \bigcup\{ \cR_{i,x,\beta} \cup \{\mu_{\beta,i}\}) \such \beta < \alpha\}\\
  &\mbox{ and the $\cR_{i,x,\alpha}$, $i = 1, \dots, k$, $x ={\min}, {\max}$,
    are pairwise nnc} \\
   &\mbox{ Ramsey ultrafilters that are nnc to $\cE$}\bigr),
  \end{split}
\end{equation}
will follow from the \CH, a routine enumeration according to \cite[Prop. 011]{mathias:happy} and \cite[proof of Theorem 14]{ncf1}  the following
theorem:

\begin{theorem}\label{hardestforMk(cR)} Suppose \CH,
  $\alpha< \omega_2$, $\cf(\alpha) = \omega$, and that
$\bP_\beta, \bar{\cR_\beta}$ are as in Equation \eqref{inducforMk(cR)} and
  $\bP_\alpha$ is the countable support limit of
  $\la \bP_\beta, \GM_k(\bar{\cR_\beta}) \such \beta < \alpha \ra$.
  In $\bV^{\bP_\alpha}$, for any $i,x$, the set of positive sets 
\[
\bigl(\bigcup_{\gamma<\alpha} (\cR_{i,x,\gamma} \cup \{\mu_{\gamma,i}\})\bigr)^{+}
\]
forms a happy family that is nowhere almost a filter.
\end{theorem}

As in the proofs of Theorem~\ref{hardest} we introduce an increasing sequence  
$\la R_\gamma(\bar{\cR}_{<\gamma})\such \gamma < \alpha \ra$ of relations $R_\alpha(\bar{\cR}_{<\alpha})$
in $\bV^{\bP_\alpha}$ such that a property called
\begin{equation}\label{ind2for(cR)}
  \mbox{``$\bP_\alpha$ is
    $R_\alpha(\bar{\cR}_{<\alpha})$-preserving''}
\end{equation} 
is carried in addition to the property \eqref{inducforMk(cR)}
and properness in the inductive choice of
the iteration. \nothing{\footnote{We could restrict to $M \prec H(\chi)$
  with $N \cap \omega_1 \in \cS$ and thus get a
  $(\bar{R},\cS)$-version that would allow us to split the
  (collapses of the ) iteration
  stages into the ones in $\cS$ and the ones outside for
  some stationary set $\cS$ in $\omega_1$. As this modification for proper forcings is not well-known \cite{Sh:f, FischerFriedmanKhomskyy} e.g, we do not
  take up this additional complexity.}
}

We define a relation $R_\alpha(\bar{\cR}_{<\alpha})$
and are concerned with
statements of  the form $(\forall f)( \exists g )( f R_\alpha(\bar{\cR}_{<\alpha}) g)$ for a Borel relation $R_\alpha(\bar{\cR}_{<\alpha})$ on the Baire space in $\bV^{\bP_\alpha}$.
Here the argument $\bar{\cR}_{<\alpha}$ indicates 
the invoked parameters from the ground model: Ramsey ultrafilters in the stage 0 and names
for Ramsey ultrafilters in any later stage.

We do not write tildes below the $R_\alpha(\bar{\cR}_{<\alpha})$'s, although these
relations are 
names.

\begin{definition}\label{therelationsforMk(cR)} By induction on $\alpha \leq \omega_2$ we define the following
  relations.
  \begin{myrules}
    \item[(1)] Let $\bP_\alpha = \la \bP_\beta, \GM_k(\bar{\cR}_\gamma) \such \beta \leq \alpha, \gamma < \alpha \ra$ be defined with \eqref{inducforMk(cR)}.
      Let $i \in \{1, \dots, k\}$, $x \in \{{\min},{\max}\}$. We say that a $\bP_\alpha$-name $g$ for an element of $\roth$ is \emph{$\alpha$-$(i,x)$-positive} if  
  $1\Vdash_{\bP_\alpha} g \in (\bigcup\{ \filter(\cR_{i,x,\gamma} \cup \{\mu_{\gamma,i}\}) \such \gamma < \alpha\})^+$.
  
\item[(2)]  Let $\bP_\alpha = \la \bP_\beta, \GM_k(\bar{\cR}_\gamma) \such \beta \leq \alpha, \gamma < \alpha \ra$ be defined with \eqref{inducforMk(cR)}.
  We say $f R_{\alpha}(\bar{\cR}_{<\alpha})g$
if the following holds in $\bV^{\bP_\alpha}$:
\begin{myrules}
\item[(a)] $f= (\bar{A}, h, (i,x))$,
\item[(b)] $ \bar{A} =\la A_{\ell}\such \ell \in \omega\ra$
  is a $\subseteq$-descending sequence of $\alpha$-$(i,x)$-positive members of $\roth$,
\item[(c)] $h$ is finite-to-one, 
\item[(d)] for $j =0,1$ there are $g^j \in \roth$, with the following properties:
  \begin{myrules}
    \item[(i)] $g = g^0 \cup g^1$,
    \item[(ii)] $g^j$ is an $\alpha$-$(i,x)$-positive diagonal lower bound of $\bar{A}$ and
      \item[(iii)] $h[g^0]\cap h[g^1] = \emptyset$.
  \end{myrules}
\end{myrules} 
  \end{myrules} 
\end{definition}

\begin{definition}\label{R_alpha_preserving_forMk(cR)}
  We say 
  $\bP_\alpha$ is $R_\alpha(\bar{\cR}_{<\alpha})$-preserving if
  $\bP_\alpha$ is proper and
  \[\bP_\alpha \Vdash (\forall f \in \dom(R_\alpha(\bar{\cR}_{<\alpha})))( \exists g)( 
  f R_\alpha(\bar{\cR}_{<\alpha}) g).
  \]
\end{definition}

\begin{lemma} We assume \CH, $\alpha< \omega_2$ and $\cf(\alpha) = \omega$. If $\bP_\alpha=
 \la \bP_\beta, \GM_k(\bar{\cR}_\gamma) \such \beta \leq \alpha, \gamma < \alpha \ra$
is $R_\alpha(\bar{\cR}_{<\alpha})$ preserving then
\begin{equation*}\begin{split}
    \bP_\alpha \Vdash & (\bigcup \{\filter(\cR_{i,x,\beta} \cup\{\mu_{i,\beta}\})\such \beta < \alpha\})^+\\
    &
  \mbox{  is a happy family
    that is nowhere almost a filter}).
\end{split}\end{equation*}
\end{lemma}
\proof This follows from Definition \ref{R_alpha_preserving_forMk(cR)}. \proofend

As mentioned  a routine enumeration gives:

\begin{corollary}\label{mathias011again}
  If $\bP_\alpha \Vdash \CH$ and $\bP_\alpha$ is $R_\alpha(\bar{\cR}_{<\alpha})$-preserving then Equation~\eqref{induc2forMk(cR)} holds.
  \end{corollary}

We carry the preservation property upwards by induction.

\begin{lemma}\label{limit-preservation_forMk(cR)}
  Assume \CH, $\alpha< \omega_2$, $\cf(\alpha)= \omega$.  Let $\bP_\alpha$ be the countable support limit of $\bP_\beta$, $\beta < \alpha$.
  If for $\beta < \alpha$ such that $\cf(\beta) < \omega_1$,  $\bP_\beta$ is $R_\beta(\bar{\cR}_{<\beta})$-preserving
  and for any $\beta < \alpha$, \eqref{inducforMk(cR)} holds then
  $\bP_\alpha$ is $R_\alpha(\bar{\cR}_{<\alpha})$ preserving.
\end{lemma}

\proof As in Lemma~\ref{limit-preservation}.

This finishes the proof of Theorem~\ref{hardestforMk(cR)}.

Finally we answer Banakh's and Blass' question on the finite part near-coherence spectrum:

\begin{theorem}\label{4} Assume \CH\ and let $k \geq 1$.
  Then there is a countable support iteration of proper iterands $\bP = \la \bP_\alpha, \GM_k(\bar{\cR_{\beta}}) \such \beta < \omega_2, \alpha \leq \omega_2 \ra$ such that
  in the extension there \emph{exactly} $2k+1$ near-coherence classes of ultrafilters. Namely, one class is represented by a $P$-point of character $\omega_1$ and
  $2k$ classes represented by the names Ramsey ultrafilters $\cR_{i,x,\omega_2}$, $i = 1, \dots, k$, $x = {\min}, {\max}$.
  \end{theorem}
\proof We start with a ground model ${\bf V}$ that  
fulfils \CH. 
By \CH\ there are a $P$-point $\cE$  and a $\{1,\dots, k\}\times\{{\min},{\max}\}$-sequence of Ramsey ultrafilters $\bar{\cR}_0$, such that any two of them are
nnc. A countable support iteration
$\la \bP_\beta, \GM_k(\bar{\cR}_\gamma) \such \beta \leq \omega_2, \gamma <\omega_2 \ra$ of  proper
forcings with \eqref{inducforMk(cR)}
at any stage $\beta$ forces:
There are at least the near-coherence classes of $\cE$, $\cR_{i,x,\omega_2}$,
$i =1, \dots, k$, $x = {\min},{\max}$.

Why are there no other near-coherence classes in $\bV^{\bP_{\omega_2}}$? For a contradiction, we suppose that $\cW$ is a non-principal ultrafilter in $\bV[G]$ for some $\bP_{\omega_2}$-generic filter $G$
and $\cW$ is nnc to any $\cR_{i,x\omega_2}$, $i=1,\dots, k$, $x={\min},{\max}$ and nnc to $\cE$. Then by \cite[Lemma 5.10]{BsSh:242} there are a $p \in \bP_{\omega_2}\cap G$, say $p \in \bP_{\alpha_0}$ for an $\alpha_0 \in \omega_2$ (by \cite[Lemma 5.6]{BsSh:242}),  and an $\omega_1$-club $C \subseteq [\alpha_0,\omega_2)$ such that for each $\alpha \in C \cup \{\omega_2\}$,
$\cW \cap \bV[G\cap \bP_\alpha]$ has a $\bP_\alpha$-name (for which we again write $\cW \cap \bV[G\cap \bP_\alpha]$)  and
\[p\Vdash_{\bP_{\alpha}}  \cW \cap \bV[G \cap \bP_\alpha] \mbox{ is an ultrafilter that is nnc to $\cR_{i,x,\alpha}$ and nnc to $\cE$}.
\]

We fix $\alpha = \min(C)$.

By Theorem~\ref{fmu(E)=fmu(W)} and Theorem~\ref{EisworthforGMk(cR)} we have for any $\alpha \in C$,
\[p \Vdash_{\bP_{\alpha+1}} f_{\supp(\mu_\alpha)}(\cW \cap \bV[G\cap \bP_\alpha]) = f_{\supp(\mu_\alpha)}(\cE) \wedge \cE \mbox{ generates an ultrafilter}.
\]

Since
\[p \Vdash_{\bP_{\omega_2}} \cW \mbox{ is a filter and }\cW \supseteq \cW \cap \bV[G \cap \bP_\alpha],\]
we have
$p \Vdash_{\bP_{\omega_2}} f_{\supp(\mu_\alpha)}(\cW) = f_{\supp(\mu_\alpha)}(\cE)$.
\proofend

\begin{theorem}\label{odd_number}
Under \CH\, there is an iteration $\la \bP_\beta, \bQ_{\alpha} \such \beta \leq \omega_2, \alpha< \omega_2 \ra$ that forces exactly $2k$ near-coherence classes.
\end{theorem}
\begin{proof}
We fix a  $\diamondsuit(S^{\aleph_2}_{\aleph_1})$-sequence  $\la D_\alpha
\such \alpha \in S^{\aleph_2}_{\aleph_1}\ra$ in the ground model
or add one by  a preliminary forcing that preserves \CH.
\nothing{E.g. with $\{p \such \dom(p) \subseteq \omega_2 \times \omega_2, \rge(p) \subseteq 2, |p| < \omega_2\}$, $q \leq p$ if $q \supseteq p$, $D_\alpha = \{\beta < \alpha \such (\bigcup G)(\alpha, \beta) = 1\}$.}
At iterand $\bQ_\alpha = \GM_k(\la \cR_{i,x,\alpha} \such i = 1, \dots, k, x = {\min}, {\max}\ra)$, $\cf(\alpha) = \omega_1$, we work with only one
top ultrafilter $\cR_{k,{\min},\alpha}$ that extends the
ultrafilters $\cR_{k,{\min},\beta}$, $\beta < \alpha$.
The other top ultrafilter
$\cR_{k,{\max},\alpha}$ is varied according to a $\diamondsuit(S^{\aleph_2}_{\aleph_1})$ sequence, namely to be some Ramsey ultrafilter that is (forced to be) nnc to $D_\alpha$
that is handed down by the diamond, if the diamond
hands down a name of an ultrafilter $D_\alpha$ that is
forced to be nnc to $\cE$ and nnc to any of the
$\cR_{i,x,\alpha}$, $i = 1, \dots, k-1$, $x = {\min},{\max}$ or $(i,x) = (k,{\min})$.\footnote{Details about book-keeping functions that are defined from diamond
  sequences can for example be found in \cite[Section 2]{MdShTs:847}.}
According to Theorem~\ref{fmu(E)=fmu(W)} in stage $\bV^{\bP_{\alpha+1}}$
  we have that $D_\alpha$ is nearly coherent to $\cE$.
  If $\cW \in \bV^{\bP_{\omega_2}}$ extends $D_\alpha$,
  then in the final model $\cW$ is nearly coherent
  to $\cE$ via $f_{\mu_\alpha}$.
  At stages $\alpha$ of uncountable cofinality and below conditions where the diamond does not hand down an
  ultrafilter that is nnc to $\cE$ and nnc to any of the $\cR_{i,x,\alpha}$, $i = 1, \dots, k-1$, $x = {\min},{\max}$ or $(i,x) = (k,{\min})$,
  or at stages of countable cofinality or successor stages we choose for
  $\cR_{k,{\max},\alpha}$ any Ramsey ultrafilter that is nnc to  any of the $\cR_{i,x,\alpha}$, $i = 1, \dots, k-1$, $x = {\min},{\max}$ or $(i,x) = (k,{\min})$
  and nnc to $\cE$.  Since for any $\cP_{\omega_2}$-name $\tau$ for an ultrafilter,  there are stationarily many
  $\alpha \in S^{\aleph_2}_{\aleph_1}$ such that $\tau \cap \bV^{\bP_{\alpha}} = D_\alpha$, any ultrafilter $\cW$ in $\bP_{\omega_2}$ extension is nearly coherent to $\cE$ or nearly coherent to
  some of the 
   $\cR_{i,x,\omega_2}$, $i = 1, \dots, k-1$, $x = {\min},{\max}$ or $(i,x) = (k,{\min})$.
\end{proof}
\nothing{
\begin{proof} Let $S = \{\alpha < \omega_2 \such \cf(\alpha)=\omega_1\}$. A $\diamondsuit_S$-sequence $\la D_\alpha \such \alpha \in S\ra$, of $\cP_{<\omega_2}$-names is chosen such that
for any $\bP_{\omega_2}$-name $\tau$ of an ultrafilter, there are stationarily many
$\alpha \in S$ the forcing such that $\tau \cap \bV^{\bP_{\alpha}} = D_\alpha$.
For details, see e.g.\ \cite{MdShTs:847}.
Now at stage $\alpha \in S$, choose a Ramsey ultrafilter $\cS_\alpha $
that is nnc $D_\alpha$, $\cE$, and nnc $\cR_{i,x}$, $i = 1, \dots, k-1$, $x = {\min}, {\max}$ or $(i,x) = (k,{\min})$. Then we let 
 \[\bQ_\alpha:= \GM_k(\la \cR_{i,\min}, \cR_{i,\max} \such 1 \leq i \leq k-1\ra \concat
\la  \cR_{k,\min}, \cS_\alpha \ra)\] is proper. By the density lemma,
  it makes in an iteration with increasing $\cR_{i,x}$, as above, any ultrafilter that is nnc to $\cR_{i,\min}$, $i = 1, \dots k$  $\cR_{i,\max}$, $i =1 \dots , k-1$ and to $\cS_\alpha$
  nnc to $\cE$  in $\bV^{\bP_{\alpha+1}}$ and henceforth. 
\end{proof}
}
\begin{corollary}\label{cor4} For any $n\in \omega$, the statement ``there are exactly $n+1$ near-coherence classes of ultrafilters'' is consistent relative to {\sf ZFC}.
\end{corollary}
For $n=1$ we can take the Matet model without any localisation
\cite{Blasstoronto} or the Miller
model \cite{ncf3} or the Blass--Shelah model \cite{BsSh:242}.

\nothing{\begin{remark}\label{question_differences?}
Theorem~\ref{4} shows that actually there are many ways to force exactly
$n+1$ near coherence classes: Any $k \geq 1$, $PP$ as in Def.~\ref{PP}
with $n=|PP| \leq 2k$ can serve. We do not know any combinatorial
statement that distinguishes between two members, say with different $k$
but same $|PP|$.
\end{remark}}
The following two propositions give more information on the
decompositions of the iterands. 

\begin{proposition}\label{two_substeps}
  Fix $k$, $\bar{\cR}$ as Definition~\ref{PP}.
  We let \mbox{$\bQ_{\rm pure} =  (\Fk^\omega(\bar{\cR}), \strkast)$}  and we let
  $\name{\cU} = \{\la \ba, \check{\ba} \ra  \such \ba \in \bQ_{\rm pure}\}$.
  Then the following holds:
  \begin{myrules}
  \item[(1)] $\bQ_{\rm pure}$ is $\omega$-closed.
    \item[(2)]
      $\GM_k(\bar{\cR})$ is densely embedded into $\bQ_{\rm pure} \ast \GM_k(\name{\cU})$.
      \item[(3)]
        $\bQ_{\rm pure}$ forces  that $\name{\cU}$ is a \gmtu\ with
        $\min_i(\name{\cU}) = \cR_{i,\min}$ and $\maxl_i(\name{\cU}) = \cR_{i,\max}$
        for $i=1,\dots, k$.
      \item[(4)]
        $\bQ_{\rm pure}$ forces that 
   $ \Phi(\name{\cU})$  ($\Phi$ is defined in Def.~\ref{corek}(3)) is nnc to any filter from the ground model that is nnc $\cR_{i,x}$. $i =1,\dots,k$, $x ={\min},{\max}$.
    \end{myrules}
  \end{proposition}

\proof
(1) The forcing order $\bQ_{\rm pure}$ is $\omega$-closed by (2) of Theorem~\ref{gowersown}.

(2) We map $(s, \ba) \in \GM_k(\bar{\cR})$ to $(\ba, (s, \ba))$. This is a dense embedding. Details can be found in \cite[Proposition 3.2]{Eisworth}.

(3) 
By Theorem~\ref{gowersown} and density arguments, the first forcing, $\bQ_{\rm pure}$, forces that is generic filter
is an \gmtu. Since $\bQ_{\rm pure}$ does not add reals, only colourings and descending $\omega$-sequences from the ground model have to be considered.

Statement (4) follows from 
Theorem~\ref{fmu(E)=fmu(W)}, applied to $\bQ_{\rm pure}$.
\proofend

\begin{remark}
We sketch an alternative way to force at least $2k +1$ nameable near-cohernce classes: Force with an iteration of
iterands $\GM_k(\cU_\alpha)$, $\cU_\alpha$ a \gmtu\ in $\bV^{\bP_\alpha}$ without localisation to $\bar{\cR}$. In the successor steps one needs to prove the analogon to Theorem~\ref{one-step} about completing the destroyed \gmtu\ $\cU_\alpha$,  in 
a non-localised form. So here for the existence of $\cU_{\alpha+1} \supseteq \cU_\alpha$ we need Gower's theorem in a computation about names. We expect that this can be done similarly to the proofs in Section \ref{S3}.
However, there is no way to prevent that more nnc Ramsey ultrafilters
than the projections are created in the iteration, since the iteration could just be the $\Fk$-shadow of an iteration with \gmtu s over $\Ftwo{k+1}$.
\end{remark}

\bibliographystyle{plain}
\bibliography{../sh/lit,../sh/listb,../sh/lista}
\end{document}